\documentclass{amsart}
\usepackage{amsmath}
\usepackage{amssymb}
\usepackage{graphicx}
\usepackage[all]{xy}
\newtheorem{definition}{Definition}[section]
\newtheorem{theorem}[definition]{\bf Theorem}
\newtheorem{lemma}[definition]{\bf Lemma}
\newtheorem{corollary}[definition]{\bf Corollary}
\newtheorem{proposition}[definition]{\bf Proposition}
\newtheorem{remark}[definition]{\bf Remark}
\newtheorem{example}[definition]{\bf Example}
\newcommand{\be}{\begin{equation}}
\newcommand{\ee}{\end{equation}}

\begin{document}
\title{Representations of Some Associative Pseudoalgebras}

\author{Zhixiang Wu}
\address{School of Mathematical Sciences, Zhejiang University,
Hangzhou, 310027, P.R.China} \email{wzx@zju.edu.cn}

\thanks{This paper  is sponsored
by  NNSFC (No.11871421,No.12171129).}
 \subjclass[2010]{Primary  17B05, 18D05, 18G60}

%\date{January 1, 2001 and, in revised form, June 22, 2001.}


\keywords{Associative pseudoalgebra, Conformal Lie algebra, Schur-Weyl duality, Morita theorem}

\begin{abstract} In this paper, we generalize Schur-Weyl duality and Morita Theorem on associative algebras to those on associative $H$-pseudoalgebras. Meanwhile, we get  a plenty of associative $H$-pseudoalgebras over a cocommutative Hopf algebra $H$.
\end{abstract}

\maketitle

\section{Introduction}
In this article, $H$ is a cocommutative Hopf algebra over a field ${\bf k}$ of characteristic zero.  An $H$ module means  a left $H$ module. A finitely generated $H$ module is simply called a finite $H$ module. An $H$-pseudoalgebra $\mathcal A$ is an $H$-module with a left $H^{\otimes 2}$ linear map $\mathcal A\otimes \mathcal A\to H^{\otimes 2}\otimes_H\mathcal A$, where $\mathcal A\otimes \mathcal A$ is regarded as a left $H^{\otimes 2}$-module with $(h\otimes h')(a\otimes a')=ha\otimes h'a'$ for $h\otimes h'\in H\otimes H$ and $a\otimes a'\in \mathcal A\otimes \mathcal A$.
$H$-pseudoalgebras are  a class of Chiral algebras, which were introduced in \cite{BD}.   Lie $H$-pseudoalgebras and associative $H$-pseudoalgebras  introduced  in \cite{BDK} are $H$-pseudoalgebras with additional axioms.  A finite $H$-pseudoalgebra means  a finite $H$ module. In \cite{BDK}, the authors  classified 
 finite simple Lie $H$-pseudoalgebras and  studied representations of  Lie $H$-pseudoalgebras and associative $H$-pseudoalgebras.
 A  representation of a Lie $H$-pseudoalgebra is said to be finite if  it is a finite $H$ module. Classification of all finite irreducible representations of simple finite Lie $H$-pseudoalgebras has been  completed in \cite{BDK1, BDK2, BDK3, DM}. When $H$ is an enveloping algebra of a finite-dimensional Lie algebra, an associative $H$-pseudoalgebra  $\mathcal{A}  $ with identiy and a left symmetric $H$-pseudoalgebra  $\mathcal{L}$  with identity are isomorphic to $H\otimes A$ and $H\otimes L$ respectively, where $A$ and $L$ are $H^*$-differential associative and left symmetric algebras respectively (see \cite{R, Wu2}).
 
 From the view of $H$-pseudoalgebras, conformal Lie algebras are Lie $H$-pseudoalgebras, where $H$ is the Hopf algebra of a polynomial algebra ${\bf k}[\partial]$ with one variable $\partial$. Conversely, any ${\bf k}[\partial]$-pseudoalgebra is a conformal algebra.
 With these,  P. Kolesnikov and his coauthors obtained  many interesting results on conformal algebras, which correspond to well-known results of  ordinary algebras over a field (see \cite{K1,K2, K3}).  Since  conformal Lie algebras can be used to study the singular part of  vertex operator algebras \cite{K},  conformal algebras have been extensively received attention  from various aspects.  E. Zelmanov   studied associative  conformal algebras with idempotent elements in \cite{Ze}. B. Bakalov, V. Kac and A. Voronov introduced the cohomology of conformal Lie algebras  in \cite{BKV}. After this, Y. Su determined the cohomology groups of a general linear conformal Lie algebra  in \cite{Su}. Further, Su and his  coauthors  study representations of various conformal Lie algebras (see \cite {FHS,HS1,HS2,SY,SY1,SY2}).
 
After B. Bakalov, A. D'Andrea and V. Kac studied finite Lie $H$-pseudoalgebras in  \cite{BDK},  various $H$-pseudoalgebras were  introduced and studied   (\cite{BL, Ro,SQ1,SQ2,SQ3,Wu2,Wu3, HW}). In present paper we focus on associative $H$-pseudoalgebras. We try to describle $H$-pseudoalgebras by their representations.
 A representation of an algebra is used to describe this algebra by some linear transformations of  a  (finite-dimensional) vector space.
 Similar to  an algebra, some  $H$-pseudoalgebras can be described by $H$-pseudolinear maps of a left $H$ module, even if  some conformal algebras  cannot be embedded in a linear confromal algebras (see \cite{SY1,SY2}). If $M$ and $N$ are  left $H$-modules,  a pseudoinear map $\phi$ from $M$ to $N$ is a linear map from $M$ to $H^{\otimes 2}\otimes_H N$  such that $\phi(hm)=(1\otimes h\otimes_H1)\phi(m)$ for any $h\in H$ and $m\in M$. All $H$-pseudolinear maps constitute a left $H$-module $Chom(M,N)$, where $(h\phi)(m)=(h\otimes 1\otimes_H1)\phi(m)$.
 Then $Cend(M):=Chom(M,M)$ is an associative $H$-pseudoalgebra provided that $M$ is a finitely generated left $H$-module.  The goal of this article is to describe  
 associative $H$-pseudoalgebras by pseudolinear maps.
 
 In the following, we use $\mathcal{A},\mathcal{B},\mathcal{C},\cdots$  to denote $H$-pseudoalgebras.
 Suppose that $M$ is a left representation  of an associative $H$-pseudoalgebra $\mathcal{A}$, that is, there is an $H^{\otimes 2}$-linear map $*:\mathcal{A}\boxtimes M\to H^{\otimes 2}\otimes_HM$ such that $(a*b)*m=a*(b*m)$ for any $a,b\in \mathcal{A},m\in M$. In this case we say that $M$ is a left $\mathcal{A}$ pseudomodule.  Let $R$ be an ordinary algebra over the field ${\bf k}$. Then an $\mathcal{A}$-$R$ hybrid bimodule is both left $\mathcal A$ pseudomodule and right $R$ module  such that  $$(a*m)(1\otimes 1\otimes_Hr)=a*(mr)$$ for any $a\in \mathcal{A},m\in M$ and $r\in R$. Similarly, one can define a $R$-$\mathcal{A}$ hybrid bimodule.
 For any $\mathcal{A}$-$R$ hybrid bimodule $M$,  if it is  finite as a left $H$ module, then  we can obtain an ordinary algebra $R':=End(^{\mathcal{A}}M)$ of all endomorphisms of the left $\mathcal{A}$ pseudomodule $M$ and an associative $H$-pseudoalgebra $\mathcal{A}':=\{\phi\in Cend(M)|\phi(m)(1\otimes 1\otimes_Hr)=\phi(mr)$, $m\in M$ and $r\in R\}$.  It is obvious that $M$ is an $\mathcal{A}'$-$R'$ hybrid bimodule. Continuing this way, we can define $(\mathcal A')'$ and $(R')'$. It is easy to prove that  $(\mathcal A')'=\mathcal A'$ and $R'=(R')'$.  This is called  double centralizer property.
 The double centralizer property of bimodules over ordinary algebras plays a central role in the representation theory of algebraic Lie theory. For example, the Schur-Weyl duality between the general linear group $GL(V )$ and the symmetric group $S_n$ on the $n$-tensor space $V ^{\otimes n}$ (\cite{We, CL}) implies that the Schur algebra $S(m, n)$ has the double centralizer property with respect to $V^{\otimes n}$. The Schur-Weyl duality can be generalized to a ${\bf k}[S_n]$-$\mathcal{A}$ hybrid  bimodule  for some  associative $H$-pseudoalgebra $\mathcal A$.

The classical Morita theorem describes when two ordinary rings have an equivalent module categories.  Morita Theorem has been generalized in different settings. For example, Morita Theorem  for triangulated categories has been generalized to derived categories \cite{Ke}.  B. Keller, D. Murfet, and M. Van den Bergh obtained Morita-type Theorem  for cluster categories \cite{KMV}. D. P. Blecher and U. Kashyap generalized  Morita Theorem to dual operator algebras in \cite{BK,Ka}. We will give  Morita Theorem for two associative $H$-pseudoalgebras.
 An associative $H$-pseudoalgebra perspective on Morita theory is based on the observation that Morita theory is a theory of bi-pseudomodules, not simply left pseudomodules or right pseudomodules. 
  
 The paper is organised as follows. In Section 2, we first fix some notions of Hopf algebras and their dual algebras. Then we recall some concepts related to associative $H$-pseudoalgebras and their representations. As we know, the annihilation algebra $X\otimes _H\mathcal A$ of an $H$-pseudoalgebra $\mathcal A$ in \cite{BDK}, where $X$ is a commutative differential $H$ bimodule algebra,  is very important to understand $\mathcal A$. We will introduce another annihilation algebra $X\otimes \mathcal A$ of  $\mathcal A$, which is both an associative ordinary algebra and an associative $H$-pseudoalgebra. Morover, it is a hybrid bimodule over itself. After this we recall some notations related to unital associative $H$-pseudoalgebras  in \cite{R} and  give some new examples of unital associative $H$-pseudoalgebras. Affected by the $H$-pseudoalgebra $Cend(H\otimes V)\simeq H\otimes H\otimes End(V)$ for any finite-dimensional vector space $V$,  we try to determine all   associative $H$-pseudoalgebra structures  on  $H\otimes B\otimes A$,   where $B$ is a bialgebra and $A$ is an ordinary associative algebra.
  In the final part of Section 2, we recall the definition of pseudolinear maps. Since there is no natural definition of pseudoproducts to make $Cend(M)$ into an associative $H$-pseudoalgebra if $M$ is not a finite $H$-module, we give two submodules of $Cend(M)$ and then prove that these two submodules are associative $H$-pseudoalgebras with a natural pseudoproduct.
 
 In Section 3, we define a hybrid bimodule. After briefly discussing the basic properties of a hybrid bimodule, we focus on the hybrid bimodules, where the ring is the group algebra of a symmetric group. We generalize the Schur-Weyl duality to these hybrid bimodules. Then we try to describe associative $H$-pseudoalgebras which have a faithful simple pseudomodule. We obtain a similar density theorem for some semismple pseudomodules. Finally, we try to describe when two hybrid bimodule categories $_R{M}od^{\mathcal{A}}$ and $_{R'}{Mod}^{\mathcal{A}}$ are equivalent.
 
 In Section 4, we define bipseudomodules and a tensor product of two bipseudomodules to prepare for Morita theory of two associative $H$-pseudoalgebras. We begin from the functor $Chom(\mathcal{A}^{\mathcal{A}},-)$. Unlike the functor $Hom_R(R,-)$ which is equivalent to the identity functor, $Chom(\mathcal{A}^{\mathcal{A}},-)$ is equivalent to the functor $H\otimes-$ on the subcategory of all unital pseudomodules.  For any associative $H$-pseudoalgebra $\mathcal A$, $H\otimes\mathcal A$ is an associative $H$-pseudoalgebra with operations $h(h'\otimes a)=(hh')\otimes a$ and $$(h\otimes a)*(h'\otimes b)=\sum\limits_ihSh_i\otimes h'\otimes_Ha_{x_i}b$$
 for any $h,h'\in H$ and $a,b\in \mathcal A$, where $a*b=\sum\limits_iSh_i\otimes 1\otimes_Ha_{x_i}b$.  In fact, $H\otimes\mathcal A\simeq Cend(\mathcal A^\mathcal A)$, the conformal endomorphism $H$-pseudoalgebra of $\mathcal A$ as a right $\mathcal A$ pseudomodule over itself for some $H$-pseudoalgebra $\mathcal A$.
  Finally, we define a kind of tensor product of two  bipseudomodules and  give a conformal Morita Theorem.

\section{Preliminaries }

In this section, we recall some concept and results related to a cocommutative Hopf algebra $H$ and $H$-pseudoalgebras. For convenience of the reader, we  recall some ``folklore" results. Meanwhile, we obtain some useful new results. 

\subsection{Hopf algebras}
In this subsection, we will fix some notations related to Hopf algebras. For more details, we refer to \cite{Sw, BDK}. Throughout this paper, ${\bf k}$ is a field of characteristic zero, and $H$ is a cocommutative Hopf algebra over the field ${\bf k}$ with coproduct $\Delta$,  counit $\varepsilon$  and an antipode $S$. All the unadorned $\otimes$ mean the tensor product over ${\bf k}$. We will simplify the sum  $\Delta(h)=\sum\limits_{(h)}h_{(1)}\otimes h_{(2)}$ into $\Delta(h)= h_{(1)}\otimes h_{(2)}$. The following equations hold.
$$(id\otimes \Delta)\Delta(h)=(\Delta\otimes id)\Delta(h),\qquad \varepsilon(h_{(1)})h_{(2)}=h_{(1)}\varepsilon(h_{(2)})=h,$$
and $$S(h_{(1)})h_{(2)}=h_{(1)}S(h_{(2)})=\varepsilon(h)$$for any $h\in H$.
There is a natural structure of a right $H$-module on $H^{\otimes n}$ given by
$(f_1\otimes\cdots\otimes f_n)g=(f_1g_{(1)})\otimes\cdots\otimes (f_ng_{(n)})$ for $f_i,g\in H$, where $$g_{(1)}\otimes\cdots\otimes g_{(n)}=(\Delta\otimes id^{n-2})(\Delta\otimes id^{n-3})\cdots (\Delta\otimes id)\Delta(g):=\Delta^{n-1}(g)$$ for $n\geq 2$. In addition, $g_{(1)}\otimes\cdots\otimes S(g_{(i)})\otimes\cdots\otimes g_{(n)}$ is denoted by $g_{(1)}\otimes\cdots\otimes g_{(-i)}\otimes\cdots\otimes g_{(n)}$. 

For any cocommutative  Hopf algebra $H$, there is an increasing subspace sequence of subspaces of $H$   \begin{eqnarray}\label{se1}\cdots \subseteq F^nH\subseteq F^{n+1}H\subseteq \cdots\end{eqnarray} where $F^n(H)=0$ for $n\leq -1$, $F^0H$ is the subspace  spanned by all group-like elements of $H$, and $F^n(H)=span_{\bf k}\{h\in H|\Delta(h)\in F^0H\otimes h+h\otimes F^0H+\sum\limits_{i=1}^{i-1}F^iH\otimes F^{n-i}H\}$ for all $n\geq 1$.  We always assume that the dimension of $F^0H$ and the dimension of $F^1H$ are finite. Then the dimension of $F^nH$ is finite for each $n\geq 0$. 

The dual space $X=H^*$ of the Hopf algebra $H$ is a commutative associative algebra. We always fix a basis of $\{h_i|i\in I\}$ of $H$ and its dual basis $\{x_i|i\in I\}\subseteq X$. Namely $x_i(h_j)=\delta_{ij}$.
For an infinite-dimensional Hopf algebra $H$, we can consider only $\Delta:X\to (H\otimes H)^*$. But for such a Hopf algebra it is possible to write $\Delta(x)$ for $x\in X$ as an infinite series $\Delta(x)=\sum\limits_{(x)}x_{(1)}\otimes x_{(2)}$. We  simply write  this sum as $\Delta(x)=x_{(1)}\otimes x_{(2)}$. There is a left and right action of $H$ on $X$ given by $$\langle xf,g\rangle =\langle x,gS(f)\rangle,\qquad \langle fx,g\rangle =\langle x,S(f)g\rangle$$  respectively for $x\in X,$ and $f,g\in H$. Using this action, we have
\begin{eqnarray}\label{d1}\Delta(x)=\sum\limits_ixS(h_i)\otimes x_i=\sum\limits_i x_i\otimes S(h_i)x.\end{eqnarray}
From (\ref{d1}), we get \begin{eqnarray}\langle x_{(1)},h\rangle \langle x_{(2)},h'\rangle=\sum\limits_i \langle  x_i,h\rangle \langle Sh_ix,h'\rangle=\langle Shx,h'\rangle=\langle x , hh'\rangle.\end{eqnarray}
  In addition, there is a linear transformation $S:X\to X$ defined by $\langle S(x),f\rangle :=\langle x,S(f)\rangle$ for $x\in X$ and $f\in H$. We denote the element  $x_{(1)}\otimes \cdots \otimes Sx_{(i)}\otimes \cdots\otimes x_{(n)}\in X^{\otimes n}$ by  $x_{(1)}\otimes \cdots \otimes x_{(-i)}\otimes \cdots\otimes x_{(n)}$. 
  
  Let  $F_{-1}X=X$ and $F_nX:=(F^nH)^{\perp}=\{x\in X|\langle x,h\rangle =0,\forall h\in F^nH\}$ for $n\geq 0$.  Then
  \begin{eqnarray}\label{se2} \cdots \supseteq F_nX\supseteq F_{n+1}X \supseteq \cdots\end{eqnarray}  Let $X^\circ=\{x\in X|\langle x,h_i\rangle=0$ for all $i\in I$ but finite$\}$, namely, $X^\circ$ is a vector space spanned by $\{x_i|i\in I\}$ over ${\bf k}$.
  
 \begin{proposition}\label{pr21}(1) $X^\circ $ is a submodule of $X$ as a left or right $H$-module.
 
 (2) For any $x\in X$, $x_{(-1)}x_{(2)}=x_{(1)}x_{(-2)}=\langle x,1\rangle \varepsilon$ and $\langle x_{(1)},1\rangle \varepsilon x_{(2)}=x_{(1)}\langle x_{(2)},1\rangle \varepsilon=x$.
 \end{proposition}
 
 \begin{proof}  (1) For any $h\in H$, $\langle hx_i, h_j\rangle=\langle x_i,Shh_j\rangle$.  Since $H=\cup_{n\geq 0}F^nH$, there are minimal $n_0$ and $n_i$ such that $Sh\in F^{n_0}H$ and $h_i\in F^{n_i}H$ respectively.  
 Then there are finite many $j\in I$ such that $Shh_j\in F^{n_i}H$ as the dimension of $F^{n_i}H$ is finite and $H$ is a domain. Hence there are finite many $j\in I$ such that $\langle hx_i, h_j\rangle \neq 0$ and $hx_i\in X^\circ$. Similarly, we have $x_ih\in X^\circ$.

 (2) For any $h\in H$, $\langle x_{(-1)}x_{(2)},h\rangle=\sum\limits_i \langle S(x_i), h_{(1)}\rangle\langle Sh_i x,h_{(2)}\rangle=\sum\limits_i \langle x_i, h_{(-1)}\rangle\langle  x, h_ih_{(2)}\rangle=\langle x,$ $ \sum\limits_i \langle x_i,h_{(-1)}\rangle h_i h_{(2)}\rangle=\langle x,h_{(-1)}h_{(2)}\rangle=\langle x, 1\rangle \varepsilon(h).$
 Thus $x_{(-1)}x_{(2)}=\langle x,1\rangle \varepsilon$. In the similar way, we can prove that  $x_{(1)}x_{(-2)}=\langle x,1\rangle \varepsilon$.  
 
 Since $\langle\langle x_{(1)},1\rangle \varepsilon x_{(2)},h\rangle=\langle x_{(1)},1\rangle \langle \varepsilon, h_{(1)}\rangle \langle  x_{(2)}, h_{(2)}\rangle=\langle x_{(1)},1\rangle \langle x_{(2)},h\rangle=\sum\limits_i \langle x_i,1\rangle\langle Sh_i x,h\rangle=\langle S(\sum\limits_i\langle x_i,1\rangle h_i) x,h\rangle=\langle x,h\rangle,$  we have $\langle x_{(1)},1\rangle \varepsilon x_{(2)}=x$. 
 In the similar way, we can obtain $\langle x_{(2)},1\rangle x_{(1)} \varepsilon$ $ =x$.
 \end{proof}
 
For any left $H$ modules  $U_i$ ($1\leq i\leq n$), we use  $U_1\boxtimes U_2\boxtimes \cdots\boxtimes U_n$ to  denote the left $H^{\otimes n}$ module $U_1\otimes U_2\otimes\cdots\otimes U_n$, where the action is given by $$(h_1\otimes h_2\otimes \cdots\otimes h_n)(u_1\otimes u_2\otimes\cdots\otimes u_n):=(h_1u_1)\otimes (h_2u_2)\otimes\cdots\otimes (h_nu_n).$$ While $U_1\otimes U_2\otimes \cdots\otimes U_n$ is still a left $H$-module with the action given by $$h(u_1\otimes u_2\otimes \cdots\otimes u_n)=(h_{(1)}u_1)\otimes (h_{(2)}u_2)\otimes\cdots\otimes (h_{(n)}u_n).$$
For any left $H$-modules $U$ and $V$, we use ${}_HU\otimes V$  to denote the left $H$-module with the acton given by $h(u\otimes v)=(hu)\otimes v$ for any $h\in H$ and $u\otimes  v\in U\otimes V$. Then we have 

\begin{lemma}\label{l2.1} ${}_HU\otimes V\simeq U\otimes V$ as left $H$-modules for any left $H$-modules $U$ and $V$. Thus $U\otimes V$ is a projective (resp. injective) left $H$-module if either $U$ or $V$ is a projective (resp. injective) $H$ module.
\end{lemma}
 \begin{proof} Let $\Phi:{}_HH\otimes V\to H\otimes V$ be a linear map given by $h\otimes u\mapsto h_{(1)}\otimes h_{(2)}u$. Then $\Phi$ is a left $H$-module homomorphism. It is easy to check that the linear map  $\Psi: H\otimes V\to {}_HH\otimes V$, given by $h\otimes u\mapsto h_{(1)}\otimes h_{(-2)}u$, is an inverse of $\Phi$. Thus $\Phi$ is an isomorphism of $H$ modules. Now if $U=H^{(J)}$ for some index set $J$, then $_HU\otimes V\simeq (_HH\otimes V)^{(J)}\simeq (H\otimes V)^{(J)}\simeq H^{(J)}\otimes V=U\otimes V$. If $M$ is a submodule of $H^{(J)}$, then $M\otimes V$ is also a submodule of $H^{(J)}\otimes V$. Moreover, the isomorphism $_HH^{(J)}\otimes V\simeq H^{(J)}\otimes V$  induces an isomorphism between $_HM\otimes V$ and $M\otimes V$. Thus $_H(H^{(J)}/M)\otimes V\simeq (H^{(J)}\otimes V)/M\otimes V\simeq (H^{(J)}/M)\otimes V$.
 
Let $U$ be a projective $H$-module. Then there is a left $H$-module $U'$  and an index set $J$ such that $U\oplus U'=H^{(J)}$. Thus $(U\oplus U')\otimes V\simeq U\otimes V\oplus U'\otimes V\simeq (H\otimes V)^{(J)}$ is a free $H$-module. Hence $U\otimes V$ is a projective $H$ module. Similarly, we can prove $U\otimes V$ is projective if $V$ is projective.
 
 Let $U$ be an injective $H$-module.  Then ${}_HU\otimes V$ is isomorphic to a direct sum of ${}_HU$. Since $H$ is noetherian, ${}_HU\otimes V$ is injective. 
 Hence $U\otimes V$ is an injective $H$-module. Since $H$ is cocommutativ, $\sigma: U\otimes V \to V\otimes U; u\otimes v\mapsto v\otimes u$ is an isomorphism of $H$ modules. Thus  $U\otimes V$ is an injective $H$ module provided that $V$ is an injective $H$ module.
 \end{proof}

 \subsection{Associative $H$-pseudoalgebras,  pseudomodules and $H$-coalgebras} In this subsection, we recall the definitions  of associative  $H$-pseudoalgebras and their representations. Then we define $H$-coalgebras with counit. For more details, we refer to  \cite{BDK} and \cite{LG}.
 
 If $\mathcal{A}$ is a left $H$-module, a pseudoproduct of $\mathcal A$ is a  left  $H^{\otimes 2}$ linear map $\mu:\mathcal{A}\boxtimes \mathcal{A}\to H^{\otimes 2}\otimes_H\mathcal{A}$. For any $a,b\in\mathcal A$, the image $\mu(a\otimes b)$ is usually denoted by $a*b$. 
 An $H$-pseudoalgebra  is a  left $H$ module $\mathcal{A}$ with a  pseudoproduct.  An $H$-pseudoalgebra $\mathcal{A}$ is called an associative $H$-pseudoalgebra with pseudoproduct $\mu$ if 
 \begin{eqnarray}\label{alw}\mu(\mu(a\otimes b)\otimes c)=\mu(a\otimes\mu(b\otimes c)),\ \text{namely,} \ (a*b)*c=a*(b*c)\end{eqnarray}
  for any $a,b,c\in A$. In  this definition, we use the following rule
 \begin{eqnarray} \label{rule} \mu((f\otimes_Ha)\otimes(g\otimes b)):=(f\Delta^{n-1}\otimes g\Delta^{m-1}\otimes_H1)\mu(a\otimes b),\end{eqnarray}
 for $f\otimes_Ha\in H^{\otimes n}\otimes_HA$ and $g\otimes_Hb\in H^{\otimes m}\otimes_HA$.
 Throughout  this paper, we use  rule (\ref{rule}) to deal with pseudoproducts. 
 
 Note that $H^{\otimes n}$ in $H^{\otimes n}\otimes_H\mathcal{A}$ means that  it is a left $H^{\otimes n}$ module and right $H$ module. Hence  $H^{\otimes n}\otimes_H\mathcal{A}$ is a left  $H^{\otimes n}$ module.

 \begin{definition}If $\mathcal{A}$ is an associative $H$-pseudoalgebra,  a left (resp. right) $\mathcal A$ pseudomodule is a left $H$ module $M$ together  with an $H^{\otimes 2}$ linear map $\rho: \mathcal{A}\boxtimes M\to H^{\otimes 2}\otimes_HM$ (resp.   $\mu: M\boxtimes \mathcal A\to H^{\otimes 2}\otimes_H M$) such that 
 $$\rho(a\otimes\rho(b\otimes m))=\rho((a*b)\otimes m)\  (resp.\   \mu(\mu(m\otimes a)\otimes b)=\mu(m\otimes (a*b)))$$ for all $a,b\in\mathcal{ A}$ and $m\in M$. For any $a\in\mathcal A$ and $m\in M$, the image $\rho(a\otimes m)$ (resp. $\mu(m\otimes a)$) is uaually simplified as $a*m$ (resp. $m*a$). We write $M={}^{\mathcal A}M$ (resp. $M=M^\mathcal A$) to indicate that $M$ is a left $\mathcal A$ (resp. right) pseudomodule.
   \end{definition}
 
\begin{definition}If $\mathcal A$ and $\mathcal B$ are two associative $H$-pseudoalgebras,  an $\mathcal A$-$\mathcal B$ bi-pseudomodule $M$ is both 
a left $\mathcal A$ pseudomodule and a right $\mathcal B$ pseudomodule such that $$(a*m)*b=a*(m*b)$$ for any $a\in \mathcal A$, $m\in M$ and $b\in \mathcal B$. We write  $M={}^\mathcal AM^\mathcal B$ to denote that $M$ is an $\mathcal A$-$\mathcal B$ bi-pseudomodule.
\end{definition}

Let $M_i$ ($i=1,2$) be two right $\mathcal A$ pseudomodules. A left $H$-module homomorphism $\phi:M_1\to M_2$ is called a homomorphism from 
 $M_1$ to $M_2$ of  $\mathcal A$ pseudomodules provided that $$\phi(m*a)=\phi(m)*a$$ for any $m\in M_1$ and $a\in\mathcal A$. Similarly, one can define a homomorphism of two left $\mathcal A$ pseudomodules and a homomorphism  of two $\mathcal A$-$\mathcal B$ bi-pseudomodules.
The category of all left $\mathcal A$ (resp. right $\mathcal B$) pseudomodules is denoted by $^\mathcal A\mathcal{M}od$ (resp. $\mathcal{M}od^\mathcal B$). The category of all  $\mathcal A$-$\mathcal B$ bi-pseudomodules is denoted by $^\mathcal A\mathcal{M}od^\mathcal B$.

For an associative $H$-pseudoalgebra $\mathcal{A}$, we use  $\mathcal{A}^{op}$ to  denote an associative $H$-pseudoalgebra with the same  $\mathcal{A}$ as left $H$-modules, and pseudoproduct $*'$  of $\mathcal{A}^{op}$  given by $$a*' b:=((12)\otimes_H1)b*a$$ for any $a,b\in \mathcal{A}$, where $((12)\otimes_H1)\sum\limits_jf_j\otimes g_j\otimes_Ha_j=\sum\limits_jg_j\otimes f_j\otimes_Ha_j$. This associative $H$-pseudoalgebra $\mathcal{A}^{op}$ is called the opposite associative $H$-pseudoalgebra of $\mathcal{A}$. If $M$ is a right $\mathcal A$ pseudomodule, then $M$ is a left $\mathcal A^{op}$ pseudomodule via  $a*m=((12)\otimes_H1)m*a$, where  $m\in M$ and $a\in\mathcal A$, and vice versa.
 
 Let $A,B,C$ be three left $H$-modules and $\mu:A\boxtimes B\to H^{\otimes 2}\otimes_HC$ be a left  $H^{\otimes 2}$-linear map. For any $a\in A,b\in B$, the image $\mu(a\otimes b)$ can be unique expressed  as $\sum\limits_{i\in I} Sh_i\otimes 1\otimes c_i$, where  $c_i$ is usually  denoted  by $\mu(a_{x_i}b)$, or simply $a_{x_i}b$. This means that $\mu(a\otimes b)=\sum\limits_{i\in I}Sh_i\otimes 1\otimes_H a_{x_i}b$. 
 Define $\mu(a_xb)$, or simply $a_xb$,   by $$a_{x}b:=\sum\limits_{i\in I}\langle x, h_i\rangle a_{x_i}b$$ for any $x\in X$.
 Sometimes $\sum\limits_{i\in I}$ is simplified as $\sum\limits_i$, or $\sum$. With these, we get that $(a*b)*c=a*(b*c)$ if and only if
 \begin{eqnarray}\label{ass1} (a_xb)_yc=a_{x_{(2)}}(b_{yx_{(-1)}}c), \text { or equivalently,} \qquad a_x(b_yc)=(a_{x_{(2)}}b)_{yx_{(1)}}c
 \end{eqnarray}
 for all $x,y\in H^*$.
 
 Generalized Lie pseudo-bialgebras in \cite{BL}, Liu and Guo introduce associative $H$-coalgebras in \cite{LG}.  In the following, we define $H$-coalgebras with unital and give some examples of $H$-coalgebras. In the last subsection of this section,  we will prove that the conformal dual of a coassociative $H$-coalgebra with a unital is an associative  $H$-pseudoalgebra with unit.
 
 \begin{definition} An $H$-coalgebra is a left $H$-module  $C$ together with  a left $H$-module homomorphism $\Delta_H:C\to C\otimes C$. The map $\Delta_H$ is called the coproduct  of this $H$-coalgebra. The coproduct $\Delta_H$ is said to be coassociative if \begin{eqnarray}\label{26}(\Delta_H\otimes 1)\Delta_H=(1\otimes \Delta_H)\Delta_H.\end{eqnarray}
A coassociative $H$-coalgebra $C$ is an $H$-coalgebra with a coassociative coproduct $\Delta_H$. As coalgebras over a field, $\Delta_H(a)$ is denoted by $a_{(1)}\otimes a_{(2)}$ for any $a\in C$

In addition, if there is a left $H$-module homomorphism $\varepsilon_H^l:C\to { H}$ such that  $$x\varepsilon_H^l(a_{(1)})_{(1)}\otimes \varepsilon^l_H(a_{(1)})_{(2)}a_{(2)}=\delta_{1,x}x\otimes a$$
for any $a\in C$ and any $x\in H^*$,  then $\varepsilon_H^l$ is called a left counit of the $H$-coalgebra $C$. If there is a left $H$-module homomorphism $\varepsilon_H^r:C\to H$ such that 
$\varepsilon_H^r(a_{(2)})a_{(1)}=a$ and 
$$\langle x,\varepsilon_H^r(a_{(2)})_{(1)}\rangle (\varepsilon_H^r(a_{(2)})_{(2)} y)\otimes a_{(1)}=\delta_{1,y}x\otimes a$$
for any $a\in C$ and $x,y\in H^*$, then $\varepsilon_A^r$ is called a right counit of the $H$-coalgebra $C$. If $\varepsilon_H^l=\varepsilon_H^r$, then we simply denote them  by $\varepsilon_H$, which is called a unit of the $H$-coalgebra $C$.
\end{definition}

\begin{remark} If $H={\bf k}$, then the counit of an $H$-coalgebra is the usual counit of this coalgebra.
\end{remark}

\begin{proposition}\label{lem54}For any two  coassociative $H$-coalgebras $A$ and $A'$,  $A\otimes A'$ is a coassociative H-coalgebra with a coproduct 
$\Delta_H=(1\otimes \tau\otimes 1)(\Delta_H\otimes \Delta_H)$, where $\tau$ is the flip map given by $\tau(a\otimes a'):=a'\otimes a$ for any $a\otimes a'\in A\otimes A'$.
\end{proposition}

\begin{proof} For any $h\in H$, $a\in A$ and $a'\in A'$, $\Delta_H(h_{(1)}a\otimes h_{(2)}a')=(1\otimes \tau\otimes 1)((h_{(1)}\otimes h_{(2)})\Delta(a)\otimes (h_{(3)}\otimes h_{(4)})\Delta_H(a'))=
h((1\otimes\tau\otimes 1)\Delta_H(a)\otimes \Delta(a')).$ Hence $(1\otimes \tau\otimes 1)(\Delta_H\otimes \Delta):A\otimes A'\to (A\otimes A')\otimes(A\otimes A')$ is a left $H$-module homomorphism. It is routine to verify that this coproduct is coassociative.
\end{proof}

It is obvious that $H$ is a left $H$-coalgebra with the coproduct $\Delta$ of $H$. For any family of  $H$-coalgebras  $C_{\alpha}$ ($\alpha\in \Omega$),  the coproduct $\oplus_{\alpha\in \Omega}C_{\alpha}$ of these  left $H$ modules $C_\alpha$  is also an $H$-coalgebra. A morphism $f$ from an $H$-coalgebra $C_1$ to an $H$-coalgebra $C_2$
is a left $H$-module homomorphism such that $\Delta_H(f(a))=(f\otimes f)\Delta_H(a)$ for any $a\in C_1$.

\begin{example} (1) Let $C$ be an $H$-coalgebra. Then $H\otimes C$ and $_HH\otimes C $ are $H$-coalgebras via $\Delta_H(h\otimes a)=(h_{(1)}\otimes a_{(1)})\otimes (h_{(2)}\otimes a_{(2)})$  for any $h\otimes a\in H\otimes C$. The isomorphism $\Phi$, where $\Phi(h\otimes a)=h_{(1)}\otimes h_{(-2)}a$,  is also an isomorphism of $H$-coalgebras.

(2) The tensor algebra $T(H)=\oplus_{n\geq 0} H^{\otimes n}$ is an $H$-coalgebra, where $\Delta_H(h_1\otimes h_2\otimes\cdots\otimes h_n)=(h_{1(1)}\otimes h_{2(1)}\otimes\cdots\otimes h_{n(1)})\otimes (h_{1(2)}\otimes h_{2(2)}\otimes\cdots\otimes h_{n(2)})$ and $\Delta_H(1)=1\otimes 1$. Furthermore, $T(C)=\oplus_{n\geq 0}C^{\otimes n}$ is an $H$-coalgebra for any $H$-coalgebra $C$.
\end{example}

 \subsection{Annihilation algebras}
 
 An ordinary  algebra  $A$ is called a left (resp. right) $H$-differential algebra if it is a left (resp. right) $H$-module satisfying 
\begin{eqnarray}h(a b)=\sum (h_{(1)}a)( h_{(2)}b) \  \ (resp. \ (ab)h=\sum (ah_{(1)})(bh_{(2)}))
\end{eqnarray}
for $h\in H$ and $a,b\in A$. Suppose that $Y$ is an associative commutative algebra and an $H$-bimodule such that it is a left and right $H$-differential algebra. Let $\mathcal{A}$ be an associative $H$-pseudoalgebra.
Then the annihilation algebra $\mathcal{A}_Y(\mathcal{A}):= Y\otimes_H\mathcal{A}$ of $\mathcal{A}$  is defined in [BDK], whose product is given by 
$$(x\otimes_Ha)(y\otimes_Hb):=\sum\limits_i(xf_i)(yg_i)\otimes_Hc_i$$
where $a*b=\sum\limits_if_i\otimes g_i\otimes_Hc_i$.
Next, we will prove that $Y\otimes \mathcal{ A}$ is either an associative $H$-pseudoalgebra, or  an ordinary associative algebra. Thus $\mathcal{A}(Y\otimes A)=Y\otimes_HY\otimes\mathcal A$ is an associative algebra.

\begin{proposition} \label{prop2.9}Let $\mathcal{A}$ be an associative $H$-pseudoalgebra and $Y$ be an $H$ bimodule commutative associative  algebra such that it is both left and right $H$-differential algebra. Then $Y\otimes \mathcal{A}$ is an associative $H$-pseudoalgebra
with operations
$$h(x\otimes a)=\sum xh_{(-1)}\otimes h_{(2)}a$$
$$(x\otimes a)*(y\otimes b)=\sum\limits_i f_{i(1)}\otimes g_{i(1)}\otimes_H(xf_{i(2)})(yg_{i(2)})\otimes c_i$$
for any $h\in H$ and $x\otimes a,b\otimes y\in Y\otimes\mathcal A$, where  $a*b=\sum\limits_i f_i\otimes g_i\otimes_Hc_i$. 

Let $\mathbb{A}_Y(\mathcal A):=Y\otimes\mathcal A$. Then it is  a left $H$-differential algebra with the operation given by
$$(x\otimes a)(y\otimes b)=\sum\limits_i(xf_i)(yg_i)\otimes c_i$$
where $a*b=\sum\limits_i f_i\otimes g_i\otimes_Hc_i$. In addition,  there is an algebra homomorphism $\varphi:\mathcal{A}_Y(Y\otimes \mathcal A)\to \mathbb{A}_Y(\mathcal A)$ given by $\varphi(x\otimes_Hy\otimes a)=xy\otimes a$. Moreover, 
\begin{eqnarray}\label{h1}(h(x\otimes a))(y\otimes b)=(x\otimes a)(h(y\otimes b))\end{eqnarray} and
\begin{eqnarray}\label{h2}(1\otimes 1\otimes_H(x\otimes a))((y\otimes b)*(z\otimes c))=((x\otimes a)(y\otimes b))*(z\otimes c)\end{eqnarray}
for any $h\in H$, and $x\otimes a,y\otimes b,z\otimes c\in Y\otimes \mathcal A$.
\end{proposition}

\begin{proof}It is easy to check that $Y\otimes \mathcal A$ is an $H$-pseudoalgebra. Then $\mathcal{A}_Y(Y\otimes \mathcal A)=Y\otimes_HY\otimes \mathcal A$ is a left $H$-differential algebra.
Moreover,  for $a*b=\sum\limits_i f_i\otimes g_i\otimes_Hc_i$, $(x\otimes_Hx'\otimes a)(y\otimes_Hy'\otimes b)=\sum\limits_i (xf_{i(1)})(yg_{i(1)})\otimes_H(x'f_{i(2)})(y'g_{i(2)})\otimes c_i$. Then $\varphi(x\otimes_Hx'\otimes a)\varphi(y\otimes_Hy'\otimes b)=\sum\limits_i((xx')f_{i(1)})((yy')g_{i(1)})\otimes c_i=\varphi(
\sum\limits_i (xf_{i(1)})(yg_{i(1)})\otimes_H(x'f_{i(2)})(y'g_{i(2)})\otimes c_i).$
 Thus  the map $\varphi: x\otimes_Hy\otimes a\mapsto xy\otimes a$ induces a  homomorphism of algebras.
 
For any $h\in H$ and $x\otimes a,y\otimes b\in Y\otimes \mathcal A$,  we have $(h(x\otimes a))(y\otimes b)=(xh_{(-1)}\otimes h_{(2)}a)*(y\otimes b)=\sum\limits_i(xh_{(-1)}h_{(2)}Sh_i\cdot y)\otimes a_{x_i}b=\varepsilon(h)(x\otimes a)(y\otimes b)$. Similarly, we can prove that $(x\otimes a)(h(y\otimes b))=\varepsilon(h)(x\otimes a)(y\otimes b)$. Thus (\ref{h1}) holds.
 
Let  $x\otimes a,y\otimes b,z\otimes c\in Y\otimes \mathcal A$. Then $(1\otimes 1\otimes_H(x\otimes a))((y\otimes b)*(z\otimes c))=(1\otimes 1\otimes_H(x\otimes a))\sum\limits_i(Sh_{i(1)}\otimes 1)\otimes_H(ySh_{i(2)}\cdot z\otimes b_{x_i}c)=\sum\limits_{i,j}Sh_{i(1)}\otimes 1\otimes_H(xSh_j\cdot ySh_{i(2)}\cdot z\otimes a_{x_j}(b_{x_i}c))$
 and $((x\otimes a)(y\otimes b))*(z\otimes c)=\sum\limits_i(xSh_i\cdot y\otimes a_{x_i}b)*(z\otimes c)=\sum\limits_{i,j}Sh_{j(1)}\otimes 1\otimes _H(xSh_iSh_{j(2)}\cdot ySh_{j(3)}\cdot z\otimes (a_{x_i}b)_{x_j}c).$ Since $(a*b)*c=a*(b*c)$, $\sum\limits_{i,j}Sh_j\otimes Sh_{i(1)}\otimes Sh_{i(2)}\otimes 1\otimes_Ha_{x_j}(b_{x_i}c)=\sum\limits_{i,j}Sh_iSh_{j(1)}\otimes Sh_{j(2)}\otimes Sh_{j(3)}\otimes 1\otimes _H(a_{x_i}b)_{x_j}c$. Thus $(1\otimes 1\otimes_H(x\otimes a))((y\otimes b)*(z\otimes c))=\sum\limits_{i,j}Sh_{j(2)}\otimes 1\otimes_H(xSh_iSh_{j(1)}\cdot ySh_{j(3)}\cdot z\otimes (a_{x_i}b)_{x_j}c$. Because $H$ is cocommutative, $(1\otimes 1\otimes_H(x\otimes a))((y\otimes b)*(z\otimes c))=((x\otimes a)(y\otimes b))*(z\otimes c)$.
\end{proof}

 \subsection{Unital associative $H$-pseudomodule}
Let $\mathcal A$ be an associative $H$-pseudoalgebra. Then an element $1\in \mathcal A$ is called an identity of
$\mathcal A$ if $1*1=1\otimes 1\otimes _H1$, and $1_1a=a$ for any $a\in \mathcal A$. We always use $1$ to denote the identity of $\mathcal A$.
We say  $\mathcal A$ a unital associative $H$-pseudoalgebra if it has an identity.  A right (resp. left) $\mathcal A$ pseudomodule $M^{\mathcal A}$ (resp. ${}^{\mathcal A}M$) over a unital associative $H$-pseudoalgebra $\mathcal A$ is said to be unital  if $m_11=m$ for all $m\in M^{\mathcal A}$ (resp. $1_1m=m$ for all $m\in {}^{\mathcal A}M$). Thus a unital associative $H$-pseudoalgebra $\mathcal A$ is a left unital $\mathcal A$ pseudomodule.

There is a unital associative $H$-pseudoalgebra $\mathcal A$ such that $a_11\neq a$ for some $a\in \mathcal A$.

\begin{example} Let $\mathcal A=He\oplus He'$ be a free left $H$-module with basis $\{e,e'\}$. Define $\mu:\mathcal A\boxtimes\mathcal A\to H^{\otimes 2}\otimes_H\mathcal A$ via
$$e*e=1\otimes 1\otimes_He, e*e'=1\otimes 1\otimes_He', e'*e=e'*e'=0.$$
Then $\mathcal A$ is a unital associative $H$-pseudoalgebra with identity $e$. However, $0=e'_1e\neq e'$.
\end{example}

The identity of a unital associative $H$-pseudoalgebra is not unique.

\begin{example} Let $\mathcal A=He\oplus He'$ be a free left $H$-module with basis $\{e,e'\}$. Define $\mu:\mathcal A\boxtimes \mathcal A\to H^{\otimes 2}\otimes_H\mathcal A$ via
$$e*e=1\otimes 1\otimes_He, e*e'=1\otimes 1\otimes_He', e'*e=1\otimes 1\otimes_He, e'*e'=1\otimes 1\otimes_He'.$$
Then $\mathcal A$ is a unital associative $H$-pseudoalgebra with identities $e$ and $e'$. However, $e=e'_1e\neq e'$.
\end{example}

There is a unital associative $H$-pseudoalgebra $\mathcal A$ with identity $e$ such that $e_1a=a_1e=a$ for any $a\in \mathcal A$.

\begin{example} Let $\mathcal A=He\oplus He'$ be a free left $H$-module with basis $\{e,e'\}$. Define $\mu:\mathcal A\boxtimes \mathcal A\to H^{\otimes 2}\otimes_H\mathcal A$ via
$$e*e=1\otimes 1\otimes_He, e*e'=1\otimes 1\otimes_He', e'*e=1\otimes g\otimes_He', e'*e'=k\otimes 1\otimes_He',$$
where $g\in H$ is a group-like element and $k\in {\bf k}$.
Then $\mathcal A$ is a unital associative $H$-pseudoalgebra with identity $e$. If $g=1$ and $k=1$, then  $e_1e'=e'_1e=e'$ and $e_1a=a_1e$ for all $a\in \mathcal A$.
\end{example}

 \subsection{Associative $H$-pseudoalgebra $H(B,A)$}
If $B$ is a sub-Hopf algebra of $H$ and $A$ is an associative algebra, then $H\otimes B\otimes A$ can be endowed various associative $H$-pseudoalgebra structures ([Wu 2]). Next, we assume that  $A$ is an ordinary associative algebra such that $A^3\neq 0$,  $B$ is a bialgebra,  $\alpha,\beta:B^{\otimes 2}\to H$ and $\gamma:B^{\otimes 2}\to B$ are linear maps. Let $H(B,A):=H\otimes B\otimes A$, which is  a left $H$-module with the action $h_1(h_2\otimes b\otimes a):=(h_1h_2)\otimes b\otimes a$ for $h_1,h_2\in H$, $a\in A$ and $b\in B$.   
 
\begin{proposition}\label{lem23} With the above assumption, $H(B,A)$ is an associative $H$-pseudoalgebra with pseudoproduct
\begin{eqnarray}(1\otimes b_1\otimes a_1)*(1\otimes b_2\otimes a_2)=(\alpha(b_1,b_2)\otimes \beta(b_1,b_2))\otimes_H(1\otimes\gamma(b_1,b_2)\otimes a_1a_2)
\end{eqnarray}
if and only if 
\begin{eqnarray}\label{eq24}\sum \alpha(b_1,b_2)\alpha(\gamma(b_1,b_2),b_3)_{(1)}\otimes \beta(b_1,b_2)\alpha(\gamma(b_1,b_2),b_3)_{(2)}\otimes \beta(\gamma(b_1,b_2),b_3)\nonumber \\
=
\sum \alpha(b_1,\gamma(b_2,b_3))\otimes \alpha(b_2,b_3)\beta(b_1,\gamma(b_2,b_3))_{(1)}\otimes_H \beta(b_2,b_3)\beta(b_1,\gamma(b_2,b_3))_{(2)}
\end{eqnarray}for any $b_1,b_2,b_3\in B$
\end{proposition}

\begin{proof}It follows by a direct computation and using the associative law (\ref{alw}). We leave it to readers.\end{proof}

\begin{proposition} For any  left $A$ module $M$  and  $\phi,\psi:B\to H$ two linear maps, 
$H(M)=H\otimes M$ is a left pseudomodule over the pseudoalgebra $H(B,A)$  defined in Proposition \ref{lem23} with the action
\begin{eqnarray}(1\otimes b\otimes a)*(1\otimes m)=(\phi(b)\otimes \psi(b))\otimes_H(1\otimes am)
\end{eqnarray}if and only if  
\begin{eqnarray}\label{E12}\alpha(b_1,b_2)\phi(\gamma(b_1,b_2))_{(1)}\otimes\beta(b_1,b_2)\phi(\gamma(b_1,b_2))_{(2)}\otimes \psi(\gamma(b_1,b_2)) \nonumber\\
=\phi(b_1)\otimes \phi(b_2)\psi(b_1)_{(1)}\otimes \psi(b_2)\psi(b_1)_{(2)}.
\end{eqnarray}
\end{proposition}
\begin{proof} It is a straight verification. We omit the details.
\end{proof}

The next proposition indicates that there  exists   $\alpha$ and $\beta$ satisfying (\ref{eq24}) and (\ref{E12})  in some case.

\begin{proposition}\label{prop2.15}Let $B$ be a cocommutative bialgebra and $\pi:B\to H$ be a bialgebra homomorphism. Set $H(B,\pi, A):=H\otimes B\otimes A$ for any associative algebra $A$.
Then $H(B,\pi,A)$ becomes an associative $H$-pseudoalgebra with each of the following pseudoproducts
\begin{eqnarray}\label{5}(f\otimes a\otimes c)*(g\otimes b\otimes d):=(f\otimes g\pi(a_{(1)}))\otimes_H(1\otimes ba_{(2)}\otimes cd)\\
\label{9}(f\otimes a\otimes c)*(g\otimes b\otimes d):=(f\pi(b_{(1)})\otimes g)\otimes_H(1\otimes ab_{(2)}\otimes cd)\\
\label{6}(f\otimes a\otimes c)*(g\otimes b\otimes d):=(f\pi(ab_{(-1)})\otimes g)\otimes_H(1\otimes b_{(2)}\otimes cd)\\
\label{7}(f\otimes a\otimes c)*(g\otimes b\otimes d):=(f\otimes g\pi(ba_{(-1)}))\otimes_H(1\otimes a_{(2)}\otimes cd)\\
\label{8}(f\otimes a\otimes c)*(g\otimes b\otimes d):=(f\otimes g\pi(a_{(-1)}))\otimes_H(1\otimes a_{(2)}b\otimes cd)\\
\label{10}(f\otimes a\otimes c)*(g\otimes b\otimes d):=(f\pi(b_{(-1)})\otimes g)\otimes_H(1\otimes b_{(2)}a\otimes cd).
\end{eqnarray}
We denote these associative $H$-pseudoalgebras  by $\mathcal A_1,\mathcal A_2,\cdots,\mathcal A_6$ respectively. 
\end{proposition}
\begin{proof} The proof follows from a direct computation.  We omit the detail.
\end{proof}
These $H$-pseudoalgebra $\mathcal A_2,\mathcal A_4,\mathcal A_6$ are opposite $H$-pseudoalgebras of $\mathcal A_1,\mathcal A_3,\mathcal A_5$ respectively. 
If $V$ is an $A$ bimodule, then $H\otimes V$ is a left $\mathcal A_1$ pseudomodule with
$$(f\otimes a\otimes b)*(h\otimes v):=(f\otimes h\pi(a))\otimes_H(1\otimes bv),$$
and $H\otimes V$ is a right $\mathcal A_2$ pseudomodule with  
$$(h\otimes v)*(f\otimes a\otimes b):=( h\pi(a)\otimes f)\otimes_H(1\otimes vb).$$
Similarly,  $H\otimes V$  is  a left (resp. right)  $\mathcal A_3$ (resp. $\mathcal A_4$) pseudomodule  with the action given by 
$$(f\otimes a\otimes b)*(h\otimes v):=(f\pi(a_{(1)})\otimes\pi(a_{(2)}))\otimes_H(1\otimes bv),$$
$$ (resp.\ \  (h\otimes v)*(f\otimes a\otimes b):=(h\pi(a_{(1)})\otimes h\pi(a_{(2)}))\otimes_H(1\otimes vb)).$$
In particular, for any associative algebra $A$,  $H\otimes H\otimes A$ is an associative $H$ pseudoalgebra with  $h'(h\otimes g\otimes a)=h'h\otimes g\otimes a$ and $(1\otimes g\otimes a)*(1\otimes g'\otimes b):=(g'_{(-1)}\otimes 1)\otimes _H(1\otimes g'_{(2)}g\otimes ab)$.

 \subsection{Pseudolinear mapping}
In this subsection, we recall the definition of pseudolinear maps. Using pseudolinear maps, we obtain two kinds of associative
$H$-pseudoalgebras $Cend^{\bf r}(V)$ and $Cend^{\circ}(V)$ from arbitrary left $H$-module $V$.

 An $H$-pseudolinear map $\phi$ from a left $H$ module $V$ to a left $H$ module $W$ is a linear map $\phi:V\to H^{\otimes2}\otimes_HW$ such that $\phi(hv)=((1\otimes h)\otimes_H1)\phi(v)$ for $h\in H$ and $v\in V$.  The space of all $H$-pseudolinear maps from $V$ to $W$ is denoted by $Chom(V,W)$. In particular,  
$Chom(V,V)$ is denoted by $Cend(V)$. It is well-known that $Chom(V,W)$ is a left $H$-module, where $h\phi(v)=((h\otimes 1)\otimes_H1)\phi(v)$ for any $h\in H$, $\phi\in Chom(V,W)$ and $v\in V$.  $\phi(v)$ is also  denoted by $\phi*v$ in the sequel. We also obtain a left $H^{\otimes 2}$ linear map $*: Chom(V,W)\boxtimes V\to H^{\otimes 2}\otimes_HW, \phi\otimes v\mapsto \phi(v)=\phi*v.$ Thus, if $V$ is a finitely generated left $H$-module, then $Cend(V)\otimes V\to H^{\otimes 2}\otimes_HV, \phi\otimes v\mapsto \phi*v$ is a left $Cend(V)$ pseudomodule.

It has been proved  that $Chom(V,W)$ is isomorphic to $Hom_H(V, H\otimes W)$ as left $H$-modules. If $V=H\otimes V_0$ is a  free $H$-module, then $Chom(V,{\bf k})\simeq H\otimes V_0^*$ (see \cite{Wu2}).
If $\varepsilon_H^l$ ($\varepsilon_H^r$)  is a left (resp. right) counit of an $H$-coalgebra  $C$, then we define $\varepsilon^{cl}:C\to H^{\otimes 2}\otimes_H{\bf k}$ (resp. $\varepsilon^{cr}:C\to H^{\otimes 2}\otimes_H{\bf k}$) via $\varepsilon^{cl}*a=1\otimes \varepsilon_H^l(a)\otimes _H1=\sum\limits_i\langle \varepsilon_H^l(a), x_i\rangle \otimes h_i\otimes_H1$ (resp. $\varepsilon^{cr}*a=1\otimes \varepsilon^r_H(a)\otimes_H1$). Thus $\varepsilon^{cl}_{x_i}a=\langle \varepsilon_H^l(a),x_i\rangle$ (resp.
$\varepsilon^{cr}_{x_i}a=\langle \varepsilon_H^r(a), x_i\rangle )$. 

\begin{proposition}\label{lem53}Let $C$ be a coassociative $H$-coalgebra with a left counit $\varepsilon^{cl}$ (resp. right counit $\varepsilon^{cr}$). Then $Chom(C,{\bf k})$ is an associative $H$-pseudoalgebra with
$\phi*\psi=\sum\limits_iSh_i\otimes1\otimes_H(\phi_{x_i}\psi)$, where $(\phi_x\psi)_y(a)=\phi_{x_{(2)}}(a_{(1)})\psi_{yx_{(-1)}}(a_{(2)})$ for any $a\in C$, $x,y\in H^*$. Moreover,
$\varepsilon^{cl}*\varepsilon^{cl}=1\otimes 1\otimes_H\varepsilon^{cl}$ (resp. $\varepsilon^{cr}*\varepsilon^{cr}=1\otimes 1\otimes_H\varepsilon^{cr}$) and $\varepsilon^{cl}_1\phi=\phi$ (resp. $\phi_1\varepsilon^{cr}=\phi$) for any $\phi\in Chom(C,{\bf k})$.
\end{proposition}

\begin{proof} It has been proven in \cite{LG} that $Chom(C,{\bf k})$ is an associative $H$-pseudoalgebra.
Since $\phi_xha=h_{(2)}(\phi_{h_{(-1)}x}a)$ (see \cite[(9.19)]{BDK}), we have $h(\phi_xa)=\phi_{h_{(1)}x}(h_{(2)}a)$. Thus,
$$\begin{array}{lll}(\varepsilon^{cl}_x\phi)_ya&=&\varepsilon^{cl}_{x_{(2)}}(a_{(1)})\phi_{yx_{(-1)}}(a_{(2)})\\
&=&\sum\limits_i \varepsilon^{cl}_{x_{i}}(a_{(1)})\phi_{y(S(x)h_i)}(a_{(2)})\\
&=&\sum\limits_i\phi_{y(S(x)h_i)}(\langle \varepsilon_H^l(a_{(1)}),x_i\rangle a_{(2)})\\
&=&
\phi_{y(S(x)\varepsilon_H^l(a_{(1)}))}(a_{(2)})\\
&=&\delta_{1,S(x)}\phi_{yS(x)}(a).\end{array}$$ Consequently,  $\varepsilon^{cl}_1\phi=\phi$  and $\varepsilon^{cl}*\varepsilon^{cl}=1\otimes 1\otimes_H\varepsilon^{cl}$.  Similarly, we can prove $\varepsilon^{cr}*\varepsilon^{cr}=1\otimes 1\otimes_H\varepsilon^{cr}$ and $\phi_1\varepsilon^{cr}=\phi$ for any $\phi\in Chom(C,{\bf k})$.
\end{proof}

\begin{example}Let $C$ be a finite-dimensional coalgebra over the field ${\bf k}$. Then $C$ is an $H$-coalgebra with a trivial action by $H$. Thus $H\otimes C$ is an $H$-coalgebra which is a finitely generated free left $H$-module. Moreover,  $Chom(H\otimes C,{\bf k})\simeq Hom_H(H\otimes C,H)\simeq H\otimes C^*$.  Let $\{e_1,e_2,\cdots,e_n\}$ be a basis of $C$  and $\Delta(e_i)=\sum\limits_{j=1}^{n}a_{ij}\otimes e_j$. Assume that  $\{e_1^*,e_2^*,\cdots,e_n^*\}$ is the dual basis of $\{e_1,e_2,\cdots,e_n\}$.  Define $(h\otimes e_i^*)*(h'\otimes e_j):=h\otimes h'\otimes_H\delta_{ij}$. Then $Chom(H\otimes C,{\bf k})$ is a free $H$-module with a basis 
$\{1\otimes e_i^*|1\leq i\leq n\}$. If $a_{ij}=\sum\limits_{l=1}^na_{ij}^le_l$, then $((1\otimes e_i^*)_x(1\otimes e_j^*))_y(1\otimes e_p)=(1\otimes e_i^*)_{x_{(2)}}(1\otimes a_{pq})((1\otimes e_j^*)_{yx_{(-1)}}(1\otimes e_q))=\delta_{1x}\delta_{1y}\delta_{jp}a_{pq}^i=\delta_{1x}\delta_{1y}(e_i^*e_j^*)(e_p)$.
Thus the   pseudoproduct of $H\otimes C^*$ is given by $$(h\otimes a)*(h'\otimes b)=h\otimes h'\otimes_Hab$$ by  Proposition \ref{lem53}.  Hence $Chom(H\otimes C,{\bf k})$ is the current associative $H$-pseudoalgebra of $C^*$.
\end{example}

In addition to the above associative $H$-pseudoalgebras provided by $H$-coalgebras, there is  an associative $H$-pseudoalgebra $Cend(V)$ defined  in \cite{BDK}, where 
$V$ is  a finitely generated left $H$-module.  Let us briefly recall the definition of $Cend(V)$ as follows.
If $V_0$ is a
finite-dimensional vector space over the field ${\bf k}$, then $Chom(V,W)\simeq Hom_H(H\otimes V_0,H\otimes W)\simeq H\otimes V^*_0\otimes W\simeq V\otimes W$
as left $H$-modules. For any finitely generated left $H$-module $V$,
$Cend(V)$ is an associative $H$-pseudoalgebra with a pseudoproduct $``*"$ such that $(\psi*\phi)(v)=\psi(\phi(v))$ for all $\psi,\phi\in Cend(V)$ and $v\in V$, where $\psi*\varphi=\sum\limits_iSh_i\otimes 1\otimes_H\psi_{x_i}\varphi$ and $(\psi_x\varphi)_yu:=\sum\limits_{i}\psi_{x_i}(\varphi_{y(h_iS(x))}u)$ for any $u\in V$. However,
there is no a natural  pseudoproduct $*$ such that  $(\psi*\phi)(v)=\psi(\phi(v))$ for all $\psi,\phi\in Cend(V)$ and $v\in V$ if $V$ is not a finitely generated left $H$-module.
If $V$ is not finite, then $\phi_{x_i}\psi$ is not equal to zero for infinitely many $i\in I$ and therefore $\sum\limits_{i\in I}Sh_i\otimes 1\otimes_H\phi_{x_i}\psi$ is not in $H^{\otimes 2}\otimes_HCend(V)$. To overcome this obstacle, we introduce two submodules of $Cend(V)$ such that $\phi_{x_i}\psi=0$ for all $i\in I$ but finite if $\phi,\psi$ are in these two submodules. To introduce these submodules, let us recall the definition of  $ker_n(\phi)$ in \cite{BDK}.

 For any $\phi\in Chom(V,W)$, let $ker_n(\phi):=\{v\in V|\phi_xv=0$ for all $x\in F_nX\}$, where $F_nX$ is defined in (\ref{se2}).
Using the increasing sequence $F^nH\subseteq F^{n+1}H$, we introduce the following notations. 
For any left $H$-modules $V,W$ and any integer $n$, let $$Chom^n(V,W):=\{\varphi\in Chom(V,W)|\varphi*v\in F^nH\otimes1\otimes_HW,\forall v\in V\}.$$ This means that  $Chom^n(V,W)=\{\varphi\in Chom(V,W)|ker_n(\varphi)=V\}$. Let $$Chom^{\mathbf r}(V,W):=\cup_{n\geq 0}Chom^n(V,W).$$ We call the  pseudolinear maps in $Chom^{\mathbf r}(V,W)$ {rational pseudolinear maps}. Usually, $Chom^{\mathbf r}(V,V)$ is denoted by $Cend^\mathbf r(V)$.

\begin{proposition}For any two left $H$ free modules $V=H\otimes V_0$ and $W=H\otimes W_0$,  $Chom^\mathbf r(V,W)$ $\supseteq H\otimes H\otimes Hom(V_0,W_0)$.
If $V_0$ is finite-dimensional, then $Chom^{\mathbf r}(V,W)=H\otimes H\otimes Hom(V_0,W_0)$.
\end{proposition}

\begin{proof}
Since $V=H\otimes V_0$, we have $Chom(V,W)\simeq Hom_H(H\otimes V_0, H\otimes W)\simeq Hom_{\bf k}(V_0, H\otimes W)$. In fact, define $\Phi:Hom_H(H\otimes V_0,H\otimes W)\to Hom_{\bf k}(V_0, H\otimes W)$ via $\Phi(\alpha)(v_0)=\alpha(1\otimes v_0)$ and define
$\Psi:Hom_{\bf k}(V_0,H\otimes W)\to Hom_H(H\otimes V_0,H\otimes W)$ via $\Psi(\beta)(1\otimes v_0)=\beta(v_0)$. Then $\Psi(\Phi(\alpha))(1\otimes v_0)=\Phi(\alpha)(v_0)=
\alpha(1\otimes v_0)$ and $\Phi(\Psi(\beta))(v_0)=\Psi(\beta)(1\otimes v_0)=\beta(v_0)$. Hence $Hom_H(H\otimes V_0, H\otimes W)\simeq Hom_{\bf k}(V_0,H\otimes W)=Hom_{\bf k}(V_0, \cup_{n\geq 0}F^nH\otimes W)\supseteq \cup_{n\geq 0}Hom_{\bf k}(V_0, F^nH\otimes W)\simeq \cup_{n\geq 0}F^nH\otimes Hom_{\bf k}(V_0,W)$. Further we assume that $W=H\otimes W_0$ for some vector space $W_0$. Then there is an injective linear map from $H\otimes Hom_{\bf k}(V_0,W_0)$
 to $Hom_{\bf k}(V_0,H\otimes W_0)$. Hence $H\otimes H\otimes Hom_{\bf k}(V_0,W_0)$ can be regarded as an subspace of $Chom^{\bf r}(H\otimes V_0,H\otimes W_0)$ up to above isomorphisms and inclusions.
 
 From the above arguments, we obtain that $Chom^{\bf r}(V,W)= H\otimes H\otimes Hom(V_0,W_0)$ if $V_0$ is finite-dimensional.
 \end{proof}
 
  In particular, $Cend^{\bf r}(H\otimes V_0)\supseteq H\otimes H\otimes End(V_0)$. From Proposition \ref{prop2.15}, we know that $ H\otimes H\otimes End(V_0)$ is an associative $H$-pseudoalgebra. Next, we prove that $Cend^{\bf r}(V)$ is also an associative $H$-pseudoalgebra.
 
 \begin{proposition} \label{prop21}For any three left $H$-modules $U, V,W $,  there is a unique ${H^{\otimes 2}}$ linear map
$\mu:Chom^{\bf r}(V, W)\boxtimes  Chom^{\bf r}(U, V) \to H^{\otimes 2}\otimes _HChom^{\bf r}(U,$ $ W)$, denoted as $\mu(\alpha,\beta)=\alpha*\beta$, such that
$$(\alpha*\beta)*u=\alpha*(\beta*u)$$ in $H^{\otimes 3}\otimes_H W$ for any $\alpha\in Chom^{\bf r}(V, W), \beta\in Chom^{\bf r}(U, V),$ and $u\in U.$ In particular, $Cend^{\bf r}(U)$ is an  associative $H$-pseudoalgebra. 
\end{proposition}

\begin{proof}Let $\{h_i\}$ and $\{x_i\}$ be dual bases of $H$ and $X=H^*$ respectively. If $\alpha\in Chom^n(V,W)$ and $\beta\in Chom^m(U,V)$, then $\beta_xu=0$ (resp. $\alpha_yv=0$) for any $x\in F_kX$ (resp. $y\in F_lX$) for any $u\in U$ (resp. $v\in V$), where $k\geq m$ (resp. $l\geq n$).
Thus $(\alpha_x\beta)_yu=\alpha_{x_{(2)}}(\beta_{yx_{(-1)}}u)=0$ if $y\in F_mX$. Hence $\alpha_x\beta\in Chom^m(U,W)$.
Therefore  $\alpha*\beta=\sum\limits_iSh_i\otimes 1\otimes_H(\alpha_{x_i}\beta)\in H^{\otimes2}\otimes_HChom^{\bf r}(U,W)$. Similar to [Wu2, Proposition 5.2], we can prove that $(\alpha*\beta)*v=\alpha*(\beta*v)$ for any $\alpha\in Chom^r(V,W),\beta\in  Chom^r(U,V)$  and  $v\in U$.
 \end{proof}

Next, we introduce another submodule of $Chom(V,W)$ for any two left $H$-modules $V,W$. At first, let us define  almost zero pseudolinear maps.
A pseudolinear map $\phi:V\to W$ is said to be { almost zero} if  $\phi(V_1)=0$, where $V_1$ is a submodule of the $H$-module $V$ such that $V=V_1+V_2$ for some finitely generated submodule $V_2$. Let $Chom^{\circ}(V,W):=\{\phi\in Chom(V,W)|\phi \text{ is almost zero} \}$  and $Cend^\circ (V):=Chom^\circ (V,V)$.
Similar to Proposition \ref{prop21}, we can prove the following 

\begin{proposition} For any three  left $H$-modules $U, V,W $,  there  is a unique $H^{\otimes 2}$ linear map
$\mu: Chom(V, W)\boxtimes Chom^\circ(U, V)\to H^{\otimes 2}\otimes _H Chom^\circ (U,W)$, denoted as $\mu(\alpha,\beta)=\alpha*\beta$, such that
$$(\alpha*\beta)*u=\alpha*(\beta*u)$$ in $H^{\otimes 3}\otimes_H W$ for any $\alpha\in Chom(V, W), \beta\in Chom^\circ(U, V)$ and $ u\in U.$ In particular, $Cend^\circ(U)$ is an  associative $H$-pseudoalgebra. 
\end{proposition}

Let $V$ be a free left $H$-module with a basis $\{e_i|i\in I\}$. Define $e_i^*\in Chom (V,{\bf k})\simeq Hom_H(V,H)$ by
$$e^*_i*(e_j)=(1\otimes 1)\otimes_H\delta_{ij}.$$ It is easy to prove that $\{e^*_i|i\in I\}$ are $H$-linearly  independent, and every almost zero pseudolinear map in $Chom (V,{\bf k})$ is an $H$-linear combination of some ${e_i^* } 's$.  Hence  $Chom^{\circ}(V,{\bf k})$ is a free left $H$-module with a basis $\{e_i^*|i\in I\}$. A pseudolinear map
$f\otimes g\otimes \alpha \in H^{\otimes 2}\otimes End(V_0)$ is a almost zero pseuolinear map of $Cend(H\otimes V_0)$ if and only if the codimension of the kernel of $\alpha $ is finite, where the action of $f\otimes g\otimes \alpha $ on $h\otimes v\in H\otimes V$ is given by $$(f\otimes g\otimes \alpha )*(h\otimes v):=(f\otimes hg)\otimes_H (1\otimes \alpha(v)).$$
Under this action, the pseudoproduct of $H\otimes H\otimes End(V_0)$ is given by 
\begin{eqnarray} \label{eq1}(f\otimes a\otimes \alpha)*(g\otimes b\otimes \beta):=(f\otimes ga_{(1)})\otimes_H(1\otimes ba_{(2)}\otimes \alpha\beta)\end{eqnarray}
by Proposition \ref{prop2.15}.
Let $Hom_{\bf k}^{\circ}(V_0,W_0):=\{f|$ $f$ is ${\bf k}$-linear map and the codimension of $Ker(f)$ is finite$\}.$ Then we have the following result.

\begin{theorem}\label{prop22} Let $V,W$ be two left $H$-modules, where  $V=H\otimes V_0$ is  free. Then $Chom^{\circ}(V,{\bf k})\otimes W\simeq Chom^{\circ}(V,W)$  as left $H$-modules, where the action of $H$ on $Chom(V,{\bf k})\otimes W$ is given by $h(\phi\otimes w):=(h_{(1)}\phi)\otimes (h_{(2)}w)$.  If further $W=H\otimes W_0$ for some vector space $W_0$, then $$Chom^\circ(V,W)\simeq H\otimes H\otimes Hom_{\bf k}^\circ(V_0,W_0),$$ where $H\otimes H\otimes Hom_{\bf k}^\circ(V_0,W_0)$ is a left $H$-module with $h'(h\otimes h^{\prime\prime}\otimes\phi)=(h'h)\otimes h^{\prime\prime}\otimes \phi$. Thus $Chom^{\bf r}(V,W)\supseteq H\otimes H\otimes Hom(V_0,W_0)\supseteq Chom^{\circ}(V,W)$.
Moreover, $Cend^{\circ}(V)\simeq Chom(V,{\bf k})\otimes V\simeq H\otimes H \otimes End^{\circ}(V_0)$. Under this isomorphism, $Chom(V,{\bf k})\otimes V$ can be endowed  with an associative $H$-pseudoalgebra structure, whose pseudoproduct  is given by
\begin{eqnarray}\label{Eq2*}(\alpha\otimes u)*(\beta\otimes v)=(g_{\alpha,v(1)}\otimes 1)\otimes_H(\beta\otimes g_{\alpha,v(-2)}u),\end{eqnarray}
where $\alpha(v)=g_{\alpha,v}\otimes 1\otimes_H1$.
\end{theorem}

\begin{proof} Let  $\{e_i|i\in I\}$ be a basis of $V_0$. Define $\Phi:Chom^{\circ}(V,{\bf k})\otimes W\to Chom^{\circ}(V,W)$ via 
\begin{eqnarray}\label{p1}\Phi(f\otimes w)(v):=(1\otimes S(g_{f,v}))\otimes_Hw\end{eqnarray} for $w\in W$, $v\in V$ and $f\in Chom^\circ(V,{\bf k})$, where $f(v)=g_{f,v}\otimes1\otimes_H1\in H\otimes H\otimes_H{\bf k}$.  Since $(h_1\alpha)*(h_2v)=h_1g_{\alpha,h_2v}\otimes 1\otimes_H1$, we have $\Phi((h_1\alpha)\otimes (h_2v))*w=1\otimes S(h_1g_{\alpha,w})\otimes_H h_2v.$ Thus $\Phi(h_{(1)}f\otimes h_{(2)}w)*v=1\otimes Sg_{\alpha,v}h_{(-1)}\otimes_Hh_{(2)}w=h\otimes Sg_{\alpha,v}\otimes_Hw.$ Consequently, $\Phi(h(\alpha\otimes w))=h\Phi(\alpha\otimes w)$ and $\Phi$ is a homomorphism of left $H$-modules. Let $V_0^\circ$ be the subspace of $V_0^*$ spanned by $\{e^*_i|i\in I\}$.
Since $Chom^\circ(V,{\bf k})\simeq H\otimes V_0^\circ$, the image $im(\Phi)$ of $\Phi$ is contained in $ Chom^{\circ }(V,W)$. Let $g\in Chom^{\circ}(V,W)$. Then there are $i_1,\cdots,i_t\in I$ such that $g(e_i)=0$ if $i\in I\setminus \{i_1,\cdots,i_t\}$. Assume that $g(e_{i_k})=\sum\limits_j(1\otimes h_{i_k,j})\otimes_H w_{i_k,j}$ for $k=1,2,\cdots, t$,
where $h_{i_k,j}\in H$ and $w_{i_k,j}\in W$. Define $f_{i_k,j}\in Chom^{\circ}(V,{\bf k})$ by
$f_{i_k,j}(e_l):=(S(h_{i_k,j})\otimes 1)\otimes_H\delta_{i_k,l}$. Then $g=\Phi(\sum\limits_{k,j}f_{i_k,j}\otimes w_{i_k,j})$ and $\Phi$ is surjective. By [BL, Proposition 4.2], $\Phi$ is an isomorphism of left $H$-modules.

Furthermore, let $V=H\otimes V_0$ and $W=H\otimes W_0$ for some vectors spaces $V_0$ and $W_0$.  From the above arguments, we get that $Chom^\circ(V,{\bf k})$ is spanned by $\{e_i^*|i\in I\}$ over $H$. Let $V_0^{\circ}$ be the vector space over ${\bf k}$ by  $\{e_i^*|i\in I\}$. Then $Chom^\circ(V,{\bf k})\simeq H\otimes V_0^\circ$. Hence $Chom^\circ(V,W)\simeq H\otimes V_0^\circ \otimes H\otimes W_0\simeq H\otimes H\otimes Hom_{\bf k}^\circ (V_0,W_0).$ By Lemma \ref{l2.1}, $Chom^\circ(V,{W})\simeq H\otimes H\otimes Hom_{\bf k}^\circ (V_0,W_0)$ as left $H$-modules.

We simply denote $\Phi(f\otimes w)(v) $ by $(f\otimes v)*v$. According to the action given by (\ref{p1}), we obtain $(\alpha\otimes u)*((\beta\otimes v)*w)=(\alpha\otimes u)*((1\otimes Sg_{\beta,w})\otimes_Hv)=1\otimes Sg_{\alpha,v(1)}\otimes
Sg_{\beta,w}Sg_{\alpha,v(2)}\otimes_Hu=((g_{\alpha,v(1)}\otimes 1)\otimes_H(\beta\otimes Sg_{\alpha,v(2)}u))*w$. Thus $(\alpha\otimes u)*(\beta\otimes v)=(g_{\alpha,v}{}_{(1)}\otimes 1)\otimes_H(\beta\otimes g_{\alpha,v}{}_{(-2)}u)$.
\end{proof}

\begin{remark} (1) Note that  $\alpha*((g_{\beta,w})_{(-2)}v)=g_{\alpha,v}\otimes (g_{\beta,w})_{(-2)}\otimes_H 1=g_{\alpha,v}(g_{\beta,w})_{(2)}\otimes 1\otimes_H 1$ for any $\alpha\in Chom(V,{\bf k})$. Using this we can check that $Chom(V,{\bf k})\otimes V$ is an associative $H$-pseudoalgebra with the pseudoproduct given by (\ref{Eq2*}) even if $V$ is not free. 

(2) If $V=H\otimes V_0$ and $W=H\otimes W_0$ for some vector spaces $V_0,W_0$ over ${\bf k}$, then $Chom^\circ(V,W)\subseteq H\otimes H\otimes Hom_{\bf k}(V_0,W_0)\subseteq Chom^r(V,W)\subseteq Chom(V,W)$.
\end{remark}

\subsection{Notations} If  $a\in \mathcal {A}$ for an associative $H$-pseudoalgebra $ \mathcal {A}$, then $r_a:\mathcal A\to H^{\otimes 2}\otimes_H\mathcal A$, $b\mapsto ((12)\otimes_H1)(b*a)$ is an element in $Cend(\mathcal A)$. To get rid of the notation $((12)\otimes_H1)$ in this definition, we will put the pseudolinear map $r_a$ on the right hand side, namely, $b*r_a=b*a$. Then $(hb)*r_a=(h\otimes 1\otimes_H1)b*r_a$ and $b*(ha)=(1\otimes h\otimes_H1)b*r_a$.
Thus we use $HomC(M,N)$ to denote the set of all linear maps $\phi:M\to H^{\otimes 2}\otimes_HN$ such that $ (hm)\phi=(h\otimes 1\otimes_H1)(m)\phi$ for any $m\in M$ and $ h\in H$. For any $h\in H$ and $\phi\in HomC(M,N)$, define $h\phi$ by $(m)(h\phi):=(1\otimes h\otimes_H1)(m)\phi$. Then $HomC(M,N)$ is also a left $H$-module. Similarly, we can define $HomC^{\bf r}(M,N)$ and $HomC^{\circ}(M,N)$.
Let $EndC(M):=HomC(M,M)$, $EndC^{\bf r}(M):=HomC^{\bf r}(M,M)$ and $EndC^{\circ}(M):=HomC^{\circ}(M,M)$. For $m\in M$ and $\phi\in HomC(M,N)$, $(m)\phi$ is also denoted by $m*\phi$.

\section{Hybrid bimodules}

\subsection{Basic properties of hybrid modules}
Let $\mathcal A$ be an associative $H$-pseudoalgebra and $R$ be an ordinary associative algebra over the field ${\bf k}$.    Assume that  $M$ is both a  left $R$ module and  a left $H$ module. Then $M$ is called a  $R$-$\mathcal A$ {\bf hybrid  bimodule} if $M$ is also a right $\mathcal A$ pseudomodule such that  $$(1\otimes 1\otimes_Hr)(m*a)=(rm)*a$$ for any $a\in\mathcal A$, $m\in M$ and $r\in R$. We write $M={}_RM^{\mathcal A}$ to indicate that $M$ is a $R$-$\mathcal A$ hybrid bimodule.
For any $a\in \mathcal A$, $s\in R$ and $m\in M$, since $(sm)*a=\sum\limits_i Sh_i\otimes 1\otimes_H(sm)_{x_i}a=(1^{\otimes 2}\otimes_Hs)\sum\limits_iSh_i\otimes1\otimes_Hm_{x_i}a$, we have $s(m_{x_i}a)=(sm)_{x_i}a$. Thus $(hsm)*a=\sum\limits_ihSh_i\otimes 1\otimes_H(sm)_{x_i}a= (1^{\otimes 2}\otimes_Hs)\sum\limits_ihSh_i\otimes 1\otimes _Hm_{x_i}a=(1^{\otimes2}\otimes_Hs)(hm)*a=(shm)*a$.  If $M$ is a unital pseudomodule, then $(sh)m=(hs)m$ for any $s\in R$, $m\in M$ and $h\in H$. Thus $M$ is a left $H\otimes R$ module. In the sequel, we always assume that any hybird bimodule $_RM^\mathcal A$ is a left $H\otimes R$ module whether it is unital or not.
 
 A  left $H$ module and right $R$ module $M$ is called can  an $\mathcal A$-$R$ hybrid bimodule if it is  a left $\mathcal A$ pseudomodule satisfying
$$a*(mr)=(a*m)(1\otimes 1\otimes_Hr)$$ for $a\in \mathcal A,m\in M$ and $r\in R$. An $\mathcal A$-$R$ hybrid bimodule $M$ is usually denoted by
${}^{\mathcal A}M_R$. We can prove that $a*((hm)*r)=a*(h(mr))$ for any $h\in H$, $a\in \mathcal A$, $r\in R$  and $m\in {}^{\mathcal A}M_R$. We also assume that any  $\mathcal A$-$R$ hybird bimodule $M$ is a left $H$ and right $R$ bimodule.

The endomorphism ring of a right (resp.  left)  $\mathcal A$  pseudomodule $M$ is denoted by $End(M^{\mathcal A})$ (resp. $End({}^{\mathcal A}M)$).
\begin{example}Let $Y$ be a  commutative left and right $H$ differential algebra. If $Y$ is also an $H$ bimodule and $\mathcal A$ is an associative $H$-pseudoalgebra, then $\mathbb{A}_Y(\mathcal A):=Y\otimes \mathcal A$ is an $\mathbb{A}_Y( \mathcal A)$-$\mathbb{A}_Y(\mathcal A)$ hybrid bimodule by Proposition \ref{prop2.9}, where $Y\otimes \mathcal A$  is regarded as  a left module over the  ordinary algebra $\mathbb{A}_Y(\mathcal A)$.
\end{example}

\begin{lemma} \label{lem31}Let $M$ be a  $R$-$\mathcal A$ hybrid bimodule. Then $l_s:M\to M, m\mapsto sm$ is an endomorphism of the right $\mathcal A$-pseudomodule $M^{\mathcal A}$ for any $s\in R$. Conversely, if $M$ is a right $\mathcal A$-pseudomodule and $R$ is a  subring of $End(M^{\mathcal A})$, then $M$ is a $R$-$\mathcal A$ hybrid bimodule  with the action given by $r\cdot m:=r(m)$.
\end{lemma}
\begin{proof}If $m*a=\sum\limits_if_i\otimes 1\otimes_Hc_i$ for $m\in M$ and $a\in \mathcal A$, then $(sm)*a=\sum\limits_if_i\otimes1\otimes_H(sc_i)=(1\otimes 1\otimes_Hl_s)(m*a)$. Since $h(sm)=s(hm)$ for all $h\in H$, $m\in M$ and $s\in R,$ $l_s$ is an endomorphism of the left $H$ module $M$. Thus $l_s$ is an endomorphism of the right $\mathcal A$ pseudomodule $M$. 

Conversely, since every $s\in R$ is an endomorphism of left $H$-module, we have $s(hm)=hsm$ for any $h\in H$ and $m\in M$. Thus $M$ is a left $H\otimes R$ module. Moreover, since 
$(1\otimes 1\otimes_Hs)(m*a)=(s(m))*a$ for any $s\in End(M^{\mathcal A})$ and $a\in \mathcal A$, $M$ is a $R$-$\mathcal A$ hybrid bimodule.
\end{proof}

\begin{example} Let  $\mathcal A$ be an associative $H$-pseudoalgebra with identity $1$. Then $End(\mathcal A^\mathcal A)=\{a\in \mathcal A|a*1=1\otimes 1\otimes_Ha\}$.
For any $a,b\in End(\mathcal A^\mathcal A)$, the multiplication $ab$ is equal to $a_1b$ for $1\in H^*$. Moreover, $\mathcal A$ becomes a left $End(\mathcal A^\mathcal A)$ module via
$a\cdot m:=a_1m$ for all $m\in\mathcal A$.
\end{example}
 
 For any two left $H\otimes R$ modules $M$ and $N$, let $$Chom_R(M,N):=\{\varphi\in Chom(M,N)|\varphi(rm)=(1^{\otimes 2}\otimes_Hr)\varphi(m), m\in M, r\in R\}.$$ Usually $Chom_R(M,M)$ is denoted by $Cend_R(M)$.   Similarly, define $$HomC_R(M,N):=\{\phi\in HomC(M,N)|(rm)*\phi=(1^{\otimes 2}\otimes_Hr)(m*\phi),m\in M, r\in R\}$$ and $EndC_R(M):=HomC_R(M,M)$.
 
We have proved that $Chom(M,N)\simeq Hom_H(M,H\otimes N)$ in \cite{Wu2}. Furthermore, we have the following result.

\begin{lemma}\label{lem33}For any two left $H\otimes R$ modules $V,W$,  $Chom_R(V,W)\simeq Hom_{H\otimes R}(V,H\otimes W)$ as left $H$-modules, where 
$H\otimes W$ is a left $H\otimes R$ module with $$(h\otimes r)(h'\otimes w):= h_{(1)}h'\otimes h_{(2)}(rw)$$
and
$Hom_{H\otimes R}(V,H\otimes W)$ is a left $H$-module with  $(h\varphi)(v)=\varphi(v)(Sh\otimes 1)$.

Similarly, $HomC(V,W)\simeq Hom_{H\otimes R}(V,H\otimes W)$ as left $H$ modules,  where $H\otimes W$ is a left $H\otimes R$ module with $$(h\otimes r)*(h'\otimes w):=hh'\otimes rw$$ and $Hom_{H\otimes R}(V,H\otimes W)$ is a left $H$ module with $(h\phi)(v):=\sum\limits_i f_ih_{(-1)}\otimes h_{(2)}w_i$ for $\phi(v)=\sum\limits_i f_i\otimes w_i$.
\end{lemma}

\begin{proof}Assume that $\phi*v=\sum\limits_if_i\otimes1\otimes_Hw_i$ for  $\phi\in Chom_R(V,W)$ and $v\in V$. Define $\hat{\phi}(v)=\sum\limits_iS(f_i)\otimes w_i\in H\otimes W$. Since $\phi*(hv)=\sum\limits_if_i\otimes h\otimes_H w_i=\sum\limits_if_ih_{(-1)}\otimes 1\otimes_Hh_{(2)}w_i$, $\hat{\phi}(hv)=\sum\limits h_{(1)}Sf_i\otimes h_{(2)}w_i=h\hat{\phi}(v)$. For any  $r\in R$ and $v\in V$, $\phi*(rv)=\sum\limits_i f_i\otimes 1\otimes_Hrw_i$. Hence $\hat{\phi}(rv)=\sum\limits_{i}S(f_i)\otimes rw_i=r\hat{\phi}(v)$.
This implies  that $\hat{\phi}\in Hom_{H\otimes R}(V,H\otimes M)$, where $H\otimes M$ is a left $H\otimes R$-module with $$(h\otimes r)(h'\otimes m):=h_{(1)}h'\otimes h_{(2)}rm.$$ 
Conversely, for any $\phi\in Hom_{H\otimes R}(V,H\otimes W)$, define $\bar{\phi}*v=\sum\limits_i Sf_i\otimes 1\otimes_Hw_i$ if  $\phi(v)=\sum\limits_i f_i\otimes w_i$.
It is easy to check that $\bar{\phi}\in Chom_R(V,W)$.

Since $S^2=1$, we have $\bar{\hat{\phi}}=\phi$ and $\hat{\bar{\psi}}=\psi$. Thus $Chom_R(V,W)\simeq Hom_{H\otimes  R}(V,H\otimes W)$ as vector spaces. For any $\phi\in Chom_R(V,W)$ and $v\in V$, if $\phi*v=\sum\limits_i f_i\otimes 1\otimes_Hw_i$, then $(h\phi)*v=\sum\limits_ihf_i\otimes1\otimes_Hw_i$. Hence $\widehat{h\phi}(v)=\sum\limits_iS(hf_i)\otimes w_i=(h\hat{\phi})(v),$ and  $Chom_R(V,W)\simeq Hom_{H\otimes R}(V,H\otimes W)$ as left $H$-modules.

Similarly, we can prove that $HomC(V,W)\simeq Hom_{H\otimes R}(V,H\otimes W)$ as left $H$ modules.
\end{proof}

From Lemma \ref{lem33}, we get the following corollary inmediately.
\begin{corollary}For any left $H\otimes R$ module $M$ and $N$, $Chom_R(M^p,N^q)\simeq M_{q\times p}(Chom_R(M,N))$ and $HomC(M^p,N^q)\simeq M_{p\times q}(HomC(M,N))$ for any postive integers $p,q$.
\end{corollary}

For any algebra $B$, we use $_B\mathcal {M}od$ to denote the category of all left $B$-modules. 
Since  $Chom_R(M,$ $N)$ and $HomC_R(M,N)$ are isomorphic to $Hom_{H\otimes R}(M,H\otimes N)$ as left $H$-modules, we have $Chom_R(M,$ $-)$ and $HomC_R(M,-)$ are left exact covariant functors from the category $_{H\otimes R}\mathcal{M}od$  to the category $_H\mathcal{M}od$. Meanwhile $Chom_R(-,N)$ and $HomC_R(-,N)$ are  left exact contravariant functors from $_R\mathcal{M}od$ to $_H\mathcal{M}od$. If $R$ is commutative, then $Chom_R(M,N)$ (resp. $HomC_R(M,N)$) is also a left $R$-module with  $(r\phi)*m:=\phi(rm)$ (resp. $(m)*(r\phi):=(rm)*\phi$) for all $r\in R$, $\phi\in Chom_R(M,N)$ and $m\in M$. Moreover, $(h(r\phi))*m=(h\otimes 1\otimes_H1)\phi(rm)$ and $(r(h\phi))*m=(h\phi)(rm)=(h\otimes 1\otimes_H1)\phi(rm)$. Thus $Chom_R(M,N)$ is also a left $H\otimes R$-module. Similarly, we can prove that $HomC_R(M,N)$ is also a left $H\otimes R$ module.  In this case both $Chom_R(M,-)$ and $HomC_R(M,-)$  are  left covariant  endofunctors of the category $_{H\otimes R}\mathcal{M}od$. Similarly, both  $Chom_R(-,N)$ and $HomC_R(-,N)$ are  left exact contravariant endofunctors  of $_{H\otimes R}\mathcal{M}od$.

\begin{remark}
 Since $Chom_R(M,-)\simeq Hom_{H\otimes R}(M,H\otimes-)$, and $H\otimes-$ is always an exact functor, we can use the functor $Chom_R(M,-)$  to study homological properties of  ${M}$. For example, if $M$ is a projective $H\otimes R$ module, then  $Chom_R(M,-)$ and $HomC_R(M,-)$  are exact functors. Similarly, the functor $Chom_R(-,N)$ is an exact  functor if $H\otimes N$ is an injective $H\otimes R$-module. 
\end{remark}

For any left $H\otimes R$ modules  $M$ and $N$, we have $ Chom_R(M,N)\simeq Hom_{H\otimes R}(M,H\otimes N)=Hom_{H\otimes R}(M, \cup_{n\geq 0}F^nH\otimes N)\supseteq \cup_{n\geq 0}Hom_{H\otimes R}(M, F^nH\otimes N)$. Define  $$Chom_R^{\bf r}(M,N):=\cup_{n\geq 0} Hom_{H\otimes R}(M,F^nH\otimes N)$$ and $Chom^{\circ}_R(M,N)=\{\phi\in Chom(M,N)|$ $\phi(M_1)=0$ for some submodule $M_1$ of $_{H\otimes R}M$ and there is  a finitely generated submodule  $M_2$ of $_{H\otimes R}M$ such that  $M=M_1+M_2\}$.
If $M=H\otimes R^{(J)}$ for some index set $J$, then $Hom_{H\otimes R}(H\otimes R^{(J)}, F^nH\otimes N)\simeq Hom_R(R^{(J)}, F^nH\otimes N)\simeq F^nH\otimes Hom_R(R^{(J)}, N)$. Assume further that $N=H\otimes R^{(J')}$. Then  $Chom_R^{\bf r}(H\otimes R^{(J)},H\otimes R^{(J')})\supseteq H\otimes H\otimes Hom_R(R^{(J)},R^{(J')})$.

Let $Cend^{\bf r}_R(M):=Chom^{\bf r}_R(M,M)$ and $Cend_R^{\circ}(M):= Chom^{\circ}(M,M)$. Then either $Cend^{\bf r}_R(M)$ or $Cend^{\circ}_R(M)$ is an associative 
$H$-pseudoalgebra. Similarly, one can define $HomC^\circ_R(V,W)$, $HomC_R^{\bf r}(V,$ $W)$, $EndC_R^\circ (V)$, $EndC^{\bf r}_R(V)$ and so on.  Moreover, we have the following lemma.

\begin{lemma}\label{lem34}Let $U,V,W$  be three left $H\otimes R$ modules and $T$ is either $r$ or $\circ$.  Then there is a unique $H^{\otimes 2}$-linear map
$\mu: HomC_R^T(U,V)\boxtimes HomC_R^T(V, W)\to H^{\otimes 2}\otimes_H  HomC_R^T(U,$ $ W)$, denoted as $\mu(\alpha,\beta)=\alpha*\beta$, such that
$$u*(\alpha*\beta)=(u*\alpha)*\beta$$ in $H^{\otimes 3}\otimes_H W$ for any $\alpha\in HomC_R^T(U, V), \beta\in HomC_R^T(V, W)$ and $ u\in U.$ In particular, $EndC_R^T(U)$ is an  associative $H$-pseudoalgebra. Moreover, $U$ is a $R$-$EndC_R^T(U)$  hybrid bimodule.

Further assume that $U=H\otimes R^{(J)}$ for some index set $J$. Then $EndC_R^\circ(U)\leq H\otimes H\otimes End_R(R^{(J)})\leq EndC_R^{\bf r}(U)$.
\end{lemma}

\begin{proof} For the sake of convenience,  let us assume that $U$ is a finitely generated $H\otimes R$-module.  The proof is similar  to that of \cite[Lemma 10.1]{BDK},  or that of  Proposition  \ref{prop21}. As an exmple, we only prove this lemma in the case when  $T=\circ$. Then $HomC^{\circ}_R(U,V)=HomC_R(U, V)$. If $\mu$ is well-define, then $\mu (\alpha,\beta)=\sum\limits_i Sh_i\otimes 1\otimes_H\alpha_{x_i}\beta$ for any $\alpha\in HomC_R(U,V)$ and $\beta\in HomC_R(V,W)$, where $\{h_i\}$ and $\{x_i\}$ are dual bases in $H$ and $H^*$ respectively.
Thus one must define an element $\alpha_x\beta$ in $HomC_R(U,W)$  for any $x\in H^*$. Let $u*\alpha=\sum\limits_{\eta}f^{\eta}\otimes 1\otimes_Hv_{\eta}$ for  $u\in U$ and $v_{\eta}*\beta=\sum\limits_{\xi}g_{\eta}^{\xi}\otimes 1\otimes_Hw_{\eta\xi}$.
Since $u_y(\alpha_x\beta)=\sum\limits_i(u_{Sh_iy}\alpha)_{xx_i}\beta=\sum\limits_{\eta,\xi}\langle y, Sg^\xi_{\eta(1)}Sf^\eta \rangle \langle x,Sg^{\eta}_{\xi(2)}\rangle w_{\eta\xi}$, we only need to define $\alpha_x\beta$ via 
$u*(\alpha_x\beta)=\sum\limits_{\eta,\xi}f^{\eta}g^{\xi}_{\eta(1)}\otimes 1\otimes_H\langle x,Sg^{\xi}_{\eta(2)}\rangle w_{\eta\xi} $. Next we prove that $\alpha_x\beta\in HomC_R(U,V)$. For any $r\in R$ and $u\in U$, 
since $(ru)*\alpha=\sum\limits_\eta f^{\eta}\otimes 1\otimes_Hrv_{\eta}$ and $(rv_{\eta})*\beta=\sum\limits_{\xi}g_{\eta}^{\xi}\otimes1\otimes_Hrw_{\eta\xi}$, we have
$(ru)*(\alpha_x\beta)=((1\otimes^{\otimes 2}\otimes_Hr)(u*(\alpha_x\beta))$. Hence $\alpha_x\beta\in HomC_R(U,W)$ for any $x\in H^*$.

Now assume that $U$ is generated by $u_1,u_2,\cdots,u_n$ as left $H\otimes R$-module, and $u_i*\alpha=\sum\limits_i f_i^{\eta}\otimes 1\otimes_Hv_{i\eta}$ and $v_{i\eta}*\beta=\sum\limits_{\xi}g^{\xi}_{i\eta}\otimes 1\otimes_Hw_{i\eta\xi}$. Let $B$ be the subspace of $H$ spanned by $f_i^{\eta}g_{i\eta(1)}^{\xi}$. Then $B$ is  finite-dimensional and $\{x\in H^*|\langle x,B\rangle=0\}$ is a finite-codimensional subspace of $H^*$. Thus there are only finitely many $i$ such that $u_j*(\alpha_{x_i}\beta)\neq 0$ for $1\leq j\leq n$. Hence there are only finitely many $i$ such that $\alpha_{x_i}\beta\neq 0$. Therefore $\alpha*\beta=\sum\limits_i Sh_i\otimes 1\otimes_H\alpha_{x_i}\beta\in H^{\otimes 2}\otimes_H HomC_R(U,W)$.

From the definition of $HomC^{\bf r}(V,W)$, we know that $EndC^{\bf r}(H\otimes R^{(J)})\supseteq H\otimes H\otimes End_R(R^{(J)})$ for any index set $J$.
If $M=R^{(J)}$ is a free left $R$-module, then  $\Phi: Hom_R(R^{(J)},R)\otimes_RR^{(J)}\to End_R(R^{(J)})$ is an injective map, where 
$\Phi(\sum\limits_jf_j\otimes_R m_j)(m):=\sum\limits_jf_j(m)m_j$. We use $End^\circ_R(R^{(J)})$ to denote the image of $\Phi$. Then $EndC_R^\circ (H\otimes R^{(J)})\simeq H\otimes H\otimes End^\circ_R(R^{(J)})$. Hence $H\otimes H\otimes End_R(R^{(J)})\supseteq EndC^\circ_R(H\otimes R^{(J)})$.
\end{proof}

\begin{remark} 
 From  Lemma \ref{lem34}, we can obtain a hybrid bimodule $_RM^{\mathcal A}$ from any $H\otimes R$ left module $M$, where $\mathcal A$ is an  associative $H$ pseudoalgebra of  either $EndC^{\circ}_R(M)$ or $EndC^{\bf r}_R(M)$. Conversely, if $M^{\mathcal A}$ is a right $\mathcal A$ pseudomodule and $R=End(M^{\mathcal A})$ the endomorphism ring of $M^{\mathcal A}$, then $M$ is a left $R$-$\mathcal A$ hybrid bimodule by Lemma \ref{lem31}. 
\end{remark}

Suppose that $M$ is a free left $H\otimes R$ module with a basis $\{e_i|i\in I\}$, where $R$ is an algebra over the field ${\bf k}$. Regard $R$ as a left $H$ module via $hr=\varepsilon (h)r$. Then $R$ becomes a left $H\otimes R$ module via $(h\otimes r)r':=\varepsilon(h)rr'$.  Define
$e_j^*\in HomC_R(M,R)\simeq Hom_{H\otimes R}(M,H\otimes R)$ by $e_i*e_j^*=1\otimes 1\otimes _H\delta_{i,j}$. Note that $HomC_R(M,R)$ is a right $R$-module with
$(m)(\phi r):=(m)\phi (1^{\otimes 2}\otimes_Hr)$.
Similar to Theorem \ref{prop22}, we have the following

\begin{corollary} \label{cor36} Let $V,W$ be two left $H\otimes R$-modules, where  $V=H\otimes R^{(I)}$ for some index set $I$. Then $HomC_R^{\circ}(V,{R})\otimes_RW\simeq HomC_R^{\circ}(V,W)$  as left $H$-modules, where the action of $H$ on $HomC_R(V,{R})\otimes_R W$ is given by $h(\phi\otimes w):=(h_{(1)}\phi)\otimes (h_{(2)}w)$.  In particular, $EndC_R^{\circ}(V)\simeq HomC_R(V,{R})\otimes_R V$. For any left $H\otimes R$ module $V$, $HomC_R(V,{R})\otimes_R V$ can be endowed with an associative $H$-pseudoalgebra structure, whose  pseudoproduct  is given by
\begin{eqnarray}\label{Eq2}(\alpha\otimes_R u)*(\beta\otimes_R v)=(1\otimes g_{\beta,v(-1)})\otimes_H(\alpha\otimes_R g_{\alpha,v(2)}u),\end{eqnarray}
where $\alpha(v)=g_{\alpha,v}\otimes 1\otimes_H1$.
\end{corollary}

\begin{corollary}Let $W$ be a left $R$-module and $V$ be a finitely generated left $R$-module. Then $HomC_R(H\otimes V,H\otimes W)\simeq H\otimes H\otimes Hom_R(V,W)$, where $H\otimes H\otimes Hom_R(V,W)$ is a left $H$-module with the action given by $h(h'\otimes a\otimes\phi)=hh'\otimes a\otimes \phi$ for $h\in H$ and $h'\otimes a\otimes \phi\in H\otimes H\otimes Hom_R(V,W)$. Under this isomorphism,
\begin{eqnarray}\label{b6}(f\otimes a\otimes \phi)*(g\otimes b\otimes \psi)=(f\otimes ga_{(-1)})\otimes_H(1\otimes a_{(2)}b\otimes \psi\phi),\\
\label{b7}(1\otimes v)*(f\otimes a\otimes\phi)=(1\otimes f)\otimes_H(a\otimes\phi(v)),\end{eqnarray} 
for any $f\otimes a\otimes \phi,g\otimes b\otimes \psi\in EndC_R(H\otimes V)$ and $1\otimes v\in H\otimes V$.
\end{corollary}

\begin{proof}By Lemma \ref{lem33}, $HomC_R(H\otimes V, H\otimes W)\simeq Hom_{H\otimes R}(H\otimes V, H\otimes H\otimes  W)\simeq Hom_R(V,H\otimes H\otimes W)\simeq H\otimes H\otimes Hom_R(V,W)$ as $V$ is finitely generated. Explicitly, 
let $v_1,v_2,\cdots,v_m$ be generators of the $R$-module $V$. We can further assume that $v_1,v_2,\cdots,v_m$ are linearly independent over ${\bf k}$. For any  $f\in Hom_{H\otimes R}(H\otimes V,H\otimes H\otimes W)$, if  $(1\otimes v_i)f=\sum\limits_j f_{ij}\otimes g_{ij}\otimes w_{ij}$, then  $f=\sum\limits_{ij}f_{ij}\otimes g_{ij}\otimes \psi_{ij}$, where $\psi_{ij}$ are  maps of $R$ modules defined by $\psi_{ij}(1\otimes v_k)=\delta_{ik}w_{ij}$ and $(1\otimes v)(f_{ij}\otimes g_{ij}\otimes\psi_{ij})=f_{ij}\otimes g_{ij}\otimes \psi_{ij}(v)$.
Define $$\Phi: Hom_{H\otimes R}(H\otimes V,H\otimes H\otimes W)\to H\otimes H\otimes Hom_R(V,W)$$ via $\Phi(f)=\sum\limits_{i,j} f_{ij}\otimes g_{ij}\otimes \psi_{ij}$.
Then $(1\otimes v)\Phi^{-1}(f\otimes a\otimes\phi)=f\otimes a\otimes \phi(v)$. Since $(1\otimes v)\Phi^{-1}(h(f\otimes a\otimes\phi))=hf\otimes a\otimes \phi(v)=(1\otimes v)(h\Phi^{-1}(f\otimes a\otimes \phi))$, $\Phi^{-1}$ is an $H$-module homomorphism. Thus $\Phi$ is an $H$-module isomorphism.
Since  $(1\otimes v)\overline{\Phi^{-1}(f\otimes a\otimes\phi)}=(f\otimes 1)\otimes_H(a\otimes \phi(v))$, we simply rewrite it as $(1\otimes v)*(f\otimes a\otimes \phi)=(1\otimes f)\otimes _H(a\otimes \phi(v))$.
Then  $((1\otimes v)*(f\otimes a\otimes\phi))*(g\otimes b\otimes\psi)=(a_{(1)}\otimes fa_{(2)}\otimes g)\otimes_H(b\otimes \psi\phi(v))=(1\otimes v)*((f\otimes ga_{(-1)})\otimes_H(1\otimes a_{(2)}b\otimes \psi\phi)).$ Hence $(f\otimes a\otimes\phi)*(g\otimes b\otimes\psi)=(f\otimes ga_{(-1)})\otimes_H(1\otimes a_{(2)}b\otimes \psi\phi)$.
\end{proof}

\begin{remark}In the sequel,    $H\otimes H\otimes End_R(V)$ denotes  an associative $H$-pseudoalgebra with the pseudoproduct given by (\ref{b6}). In addition, $H\otimes V$ is a right $H\otimes H\otimes End_R(V)$ with the action given by (\ref{b7}).
\end{remark}
It is well-known that 
$H\otimes V$ is a left $(H\otimes H\otimes End_R(V))^{op}$ pseudomodule. Thus, we have 
 
\begin{corollary}Let $W$ be an $H$-$R$ bimodule and $V$ be a finitely generated right $R$-module. Then $Chom_R(H\otimes V,H\otimes W)\simeq H\otimes H\otimes Hom_R(V,W)$, where $H\otimes H\otimes Hom_R(V,W)$ is a left $H$-module with the action given by $h(h'\otimes a\otimes\phi)=hh'\otimes a\otimes \phi$ for $h\in H$ and $h'\otimes a\otimes \phi\in H\otimes H\otimes Hom_R(V,W)$. Under this isomorphism,
\begin{eqnarray}\label{b9}(f\otimes a\otimes \phi)*(g\otimes b\otimes \psi)=(fb_{(-1)}\otimes g)\otimes_H(1\otimes b_{(2)}a\otimes \phi\psi),\\
(f\otimes a\otimes\phi)*(1\otimes v)=(f\otimes 1)\otimes_H(a\otimes\phi(v)),\end{eqnarray} 
for any $f\otimes a\otimes \phi,g\otimes b\otimes \psi\in Cend_R(H\otimes V)$ and $1\otimes v\in H\otimes V$.
\end{corollary}

\begin{corollary} \label{cor310}For any index set $J$ and a division  algebra  $R$  over ${\bf k}$, let $\mathcal A\subseteq H\otimes H\otimes End_R(R^{(J)})$ be an associative $H$-pseudoalgebra such that 
 $H\otimes R^{(I)}$ is a  right $\mathcal A$  pseudomodule with the action  given by (\ref{b7}). If $\mathcal A\supseteq H\otimes H\otimes End_R^\circ (R^{(J)})$, then 
 $H\otimes R^{(J)}$ is a simple $\mathcal A$ pseudomodule and $End((H\otimes R^{(J)})^\mathcal A)=R$.
\end{corollary}

\begin{proof} Let $N$ be a nonzero $\mathcal A$ sub-pseudomodule of $H\otimes R^{(I)}$. Then there is a nonzero element $\sum\limits_{i\in I}f_i\otimes r_ie_i\in N$.
Consequently, $0\neq (f_{i_0}\otimes 1)\otimes_H(1\otimes  e_{i_0})\in H^{\otimes 2}\otimes_HN$ for some $i_0\in I$.
 Hence $1\otimes e_{i_0}\in N$. Consequently, $N=H\otimes R^{(I)}$. 

Assume that $\varphi$ is an $\mathcal A$ pseudomodule endomorphism of $H\otimes R^{(I)}$ and  $\varphi(1\otimes e_i)=\sum\limits_{j\in I}h_{ij}\otimes r_{ij}e_j$ for some $h_{ij}\in H$ and $r_{ij}\in R$. Then $(1^{\otimes 2}\otimes _H\varphi)((1\otimes e_i)*(f\otimes 1\otimes \phi))=(1\otimes f)\otimes_H\varphi(1\otimes \phi(e_i))=\sum\limits_j(h_{ij}\otimes r_{ij}e_j)*(f\otimes 1\otimes \phi)=\sum\limits_j(h_{ij}\otimes f)\otimes_H(1\otimes r_{ij} \phi(e_j))$  for any $\phi\in End_R(R^{(I)})$ and any $f\in H$. Hence $h_{ij}\in{\bf k}$. Thus $\varphi=a\in R$.
\end{proof}

\subsection{Schur-Weyl duality}
For any  $R$-$\mathcal A$ hybrid bimodule  $M$, which is a finitely generated left $H$-module, 
let $\mathcal A^{(0)}=\mathcal A$, $R^{(0)}=R$ and  $\mathcal A^{(n)}=EndC_{R^{(n-1)}}(M)$, $R^{(n)}=End(M^{\mathcal A^{(n-1)}})$ for $n\geq 1$. Then $\mathcal A^{(n)}$ are associative $H$-pseudoalgebras and $R^{(n)}$ are ordinary associative algebras. It is easy to check that $\mathcal A^{(2)}=\mathcal A^{(1)}$ and $R^{(2)}=R^{(1)}$.

A partition $\lambda$  of $n$ is a representation of $n$ in the form $n =\lambda_1+\lambda_2+\cdots+\lambda_r$, where $\lambda_i$  are positive integers, and $\lambda_i\geq \lambda_{i+1}$. Since the characteristic of ${\bf k}$ is zero, ${\bf k}[S_n]$ is semisimple and every partition $\lambda$ of $n$ determines a unique irreducible representation $V_\lambda$ of ${\bf k}[S_n]$ and every irreducible representation of ${\bf k}[S_n]$ is isomorphic to $V_\lambda$ for some partition $\lambda$ of $n$. Let $D_\lambda$ be the endomorphism ring of $V_\lambda$. Then $H\otimes H\otimes D_\lambda\simeq EndC_{{\bf k}[S_n]}(H\otimes V_\lambda)$.

\begin{theorem} For any finite-dimensional vector space $V$,  let $M=H\otimes V^{\otimes n}$. Then $M$ is a left $H\otimes {\bf k}[S_n]$ module with the action given by 
 $$(h'\otimes \sigma) (h\otimes(v_1\otimes v_2\otimes\cdots\otimes v_n)):=h'h\otimes v_{\sigma^{-1}(1)}\otimes v_{\sigma^{-1}(2)}\otimes\cdots\otimes v_{\sigma^{-1}(n)}$$ 
 for $h'\otimes\sigma\in H\otimes{\bf k}[S_n]$ and $h\otimes v_1\otimes v_2\otimes\cdots\otimes v_n\in H\otimes v^{\otimes n}$.
In addition,  $EndC_{{\bf k}[S_n]}(H\otimes V^{\otimes n})=H\otimes H\otimes S^nEnd(V)^{\otimes n}=H\otimes H\otimes R'\simeq \oplus_\lambda H\otimes H\otimes  M_{n_\lambda}(D_\lambda)$, where $R'$ is the image of the universal enveloping algebra
$U(\mathfrak{gl}(V))$ under its natural action on $V^{\otimes n}$.
\end{theorem}

\begin{proof}
 Note that $EndC(H\otimes V^{\otimes n})=H\otimes H\otimes End(V^{\otimes n})\simeq H\otimes H\otimes End(V)^{\otimes n}$.
 For any  $(h\otimes v_1\otimes v_2\otimes \cdots\otimes v_n)\in H\otimes V^{\otimes n}$,  $(f\otimes a\otimes \varphi_1\otimes\varphi_2\otimes\cdots\otimes \varphi_n) \in EndC(H\otimes V^{\otimes n})$ and $\sigma\in S_n$,  we have $$\begin{array}{l}(h\otimes v_{\sigma^{-1}(1)}\otimes v_{\sigma^{-1}(2)}\otimes \cdots\otimes v_{\sigma^{-1}(n)})*(f\otimes a\otimes \varphi_1\otimes\varphi_2\otimes\cdots\otimes \varphi_n)\\
 =(h\otimes f)\otimes_H(a\otimes \varphi_{\sigma^{-1}(1)}(v_{\sigma^{-1}(1)})\otimes\varphi_{\sigma^{-1}(2)}(v_{\sigma^{-1}(2)})\otimes\cdots\otimes \varphi_{\sigma^{-1}(n)}(v_{\sigma^{-1}(n)}))\end{array}$$ if and only if $\varphi_1\otimes\varphi_2\otimes\cdots\otimes \varphi_n\in S^n(End(V)^{\otimes n})$.  Thus $EndC_{{\bf k}[S_n]}(H\otimes V^{\otimes n})=H\otimes H\otimes S^n(End(V)^{\otimes n})=H\otimes H\otimes R'.$

 Since ${\bf k}[S_n]$ is a semismple ring,  $V^{\otimes n}=\oplus_{\lambda}V_{\lambda}^{n_\lambda}$ is direct sum of simple ${\bf k}[S_n]$ modules $V_{\lambda}$. Thus $H\otimes V^{\otimes n}=\oplus_\lambda (H\otimes V_\lambda)^{n_\lambda}$.  Moreover, $Hom_{{\bf k}[S_n]}(V_\lambda, V_\mu)=0$ if $\lambda\neq \mu$. Hence $Hom_{{\bf k}[S_n]}(H\otimes V_\lambda,H\otimes H \otimes V_\mu)=0$ if $\lambda\neq \mu$. Consequently, $HomC_{{\bf k}[S_n]}(H\otimes V_{\lambda},H\otimes V_\mu)\simeq Hom_{H\otimes{\bf k}[S_n]}(H\otimes V_\lambda,H\otimes H\otimes V_\mu)\subseteq Hom_{{\bf k}[S_n]}(H\otimes V_\lambda,H\otimes H\otimes V_\mu)=0$ if $\lambda\neq \mu$.
Therefore, $EndC_{{\bf k}[S_n]}(H\otimes V^{\otimes n})\simeq \oplus_\lambda M_{n_\lambda}(Cend_{{\bf k}[S_n]}(H\otimes V_\lambda))\simeq \oplus_\lambda H\otimes H\otimes M_{n_\lambda}(D_\lambda)$. 
\end{proof}

\begin{theorem}(Schur-Weyl) For any finite-dimensional vector space $V$  over an algebraically closed field ${\bf k}$ and any integer $n\geq 2$, the left $H\otimes {\bf k}[S_n]$ module 
$H\otimes V^{\otimes n}$ can be decomposed as $H\otimes V^{\otimes n}=\oplus_\lambda L_\lambda\otimes V_\lambda$, where $V_\lambda$ is an irreducible representation of $S_n$, $L_\lambda$ is either zero or a simple right $H\otimes H\otimes End(V)$ pseudomodule.
\end{theorem}

\begin{proof}If ${\bf k}$ is algebraically closed, then $D_\lambda={\bf k}$ and $V_\lambda \simeq {\bf k}^{m_\lambda}$ as left ${\bf k}[S_n]$ modules, where $m_\lambda$ is the dimension of the simple module $V_\lambda$ over ${\bf k}$. Thus 
$H\otimes V^{\otimes n}\simeq \oplus_{\lambda} H\otimes V_{\lambda}^{n_{\lambda}}\simeq \oplus_\lambda H\otimes {\bf k}^{n_\lambda}\otimes {\bf k}^{m_\lambda}$. Since $EndC_{{\bf k}[S_n]}(H\otimes V^{\otimes n})=\oplus_\lambda H\otimes H\otimes M_{n_\lambda}({\bf k})$, $H\otimes {\bf k}^{n_\lambda}$ is a simple $EndC_{{\bf k}[S_n]}(H\otimes V^{\otimes n})$ pseudomodule by Corollary \ref{cor310}.
\end{proof}

Similarly, we can define an action of $S_n$ on $(H\otimes V)^{\otimes n}$ by $\sigma((h_1\otimes v_1)\otimes(h_2\otimes v_2)\otimes\cdots\otimes(h_n\otimes v_n)):=(h_{\sigma^{-1}(1)}\otimes v_1)\otimes(h_{\sigma^{-1}(2)}\otimes v_2)\otimes\cdots\otimes(h_{\sigma^{-1}(n)}\otimes v_n)$. Under this action, $(H\otimes V)^{\otimes n}$ becomes  a left ${\bf k}[S_n]$-module. Let $(H\otimes V)^{\odot n}$ be the invariant subspace of $(H\otimes V)^{\otimes n}$ under this  action. Then $(H\otimes V)^{\odot n}$ is also a submodule of $(H\otimes V)^{\otimes n}$ as left $H$-modules.
Define $\alpha:H\otimes V^{\otimes n}\to (H\otimes V)^{\odot n}$ by $\alpha(h\otimes v_1\otimes v_2\otimes\cdots\otimes v_n)=(h_{(1)}\otimes v_1)\otimes (h_{(2)}\otimes v_2)\otimes\cdots\otimes (h_{(n)}\otimes v_n)$ and $\beta:(H\otimes V)^{\odot n}\to H\otimes V^{\otimes n}$ by $\beta((h_1\otimes v_1)\otimes (h_2\otimes v_2)\otimes\cdots\otimes (h_n\otimes v_n))=\prod_{i=2}^n\varepsilon(h_i)h_1\otimes (v_1\otimes v_2\otimes\cdots\otimes v_n)$. Then $\beta\alpha=1$. Since $H$ is cocommutative, $Im\alpha$ is a submodule of the right ${\bf k}[S_n]$-module $(H\otimes V)^{\otimes n}$.

Let $\mathcal A=EndC_R((H\otimes V)^{\odot n})$. For any $1\otimes a\otimes \varphi \in EndC(H\otimes V)$ and $(h_1\otimes v_1)\otimes (h_2\otimes v_2)\otimes \cdots\otimes (h_n\otimes v_n)\in (H\otimes V)^{\otimes n}$, define  \begin{eqnarray}\label{b8}\begin{array}{l}((h_1\otimes v_1)\otimes (h_2\otimes v_2)\otimes \cdots\otimes (h_n\otimes v_n))*(1\otimes a\otimes\varphi)\\
=\sum\limits_{i=1}^n\sum\limits_{j=1}^{m_i}(h_i\otimes1)\otimes_H(h_1\otimes v_1)\otimes (h_2\otimes v_2)\otimes\cdots \otimes(h_{i-1}\otimes v_{i-1})\otimes\\
 \qquad \qquad (a\otimes \varphi(v_i))\otimes (h_{i+1}\otimes v_{i+1})\otimes\cdots\otimes(h_n\otimes v_n).\end{array}\end{eqnarray}
 Then $$\begin{array}{l} ((h_{\sigma^{-1}(1)}\otimes v_{\sigma^{-1}(1)})\otimes (h_{\sigma^{-1}(2)}\otimes v_{\sigma^{-1}(2)})\cdots(h_{\sigma^{-1}(n)}\otimes v_{\sigma^{-1}(n)}))*(1\otimes a\otimes\varphi)\\
 =\sum\limits_{i=1}^n\sum\limits_{j=1}^{m_i}(1\otimes h_{\sigma^{-1}(i)})\otimes_H(h_{\sigma^{-1}(1)}\otimes v_{\sigma^{-1}(1)}\otimes (h_{\sigma^{-1}(2)}\otimes v_{\sigma^{-1}(2)})\otimes
 \cdots \otimes(h_{\sigma^{-1}(i-1)}\otimes v_{\sigma^{-1}(i-1)})\otimes\\
 \qquad\qquad (a\otimes \varphi(v_{\sigma^{-1}(i)}))\otimes (h_{\sigma^{-1}(i+1)}\otimes v_{\sigma^{-1}(i+1)})\otimes\cdots\otimes(h_{\sigma^{-1}(n)}\otimes v_{\sigma^{-1}(n)}))\\
 =(1^{\otimes 2}\otimes_H\sigma)((h_1\otimes v_1)\otimes (h_2\otimes v_2)\cdots(h_n\otimes v_n))*(1\otimes a\otimes\varphi))\end{array}$$ provided that $h_1\otimes h_2\otimes \cdots\otimes h_n\in (H^{\otimes n})_{S_n}$, where $(H^{\otimes n})_{S_n}$ is the invariant subspace of $H^{\otimes n}$ under the action of $S_n$. Similarly, we can prove that (\ref{b8}) defines a pseudolinear map in $\mathcal A$.  Let $\mathcal A'$ be the subalgebra generated by this pseudolinear map.
Then $\mathcal A'\subseteq S^nEndC_R((H\otimes V)^{\odot n})=\mathcal A$.

\subsection{ Density Theorem}
A right $\mathcal A$ pseudomodule $M^{\mathcal A}$ is said to be nontrivial if $M^{\mathcal A}*{\mathcal A}\neq 0$.
A  nontrivial right $\mathcal A$ pseudomoule $M^{\mathcal A}$ is said to be simple if it has only tow sub-pseudomodules, namely $0$ and $M^{\mathcal A}$.
Suppose that $\varphi:M^{\mathcal A}\to N^{\mathcal A}$ is a homomorphism of  right $\mathcal A$ pseudomodules. Then $Ker(\varphi):=\{m\in M^{\mathcal A}|\varphi(m)=0\}$ is a submodule of $M^{\mathcal A}$ as $(1^{\otimes 2}\otimes_H\varphi)(m*a)=\varphi (m)*a$.
Let $Im(\varphi)=span_H(\{\varphi(m_xa)|a\in \mathcal{A},m\in M^{\mathcal A},x\in H^*\})$. Since $(1^{\otimes 2}\otimes_H\varphi)(m*a)=\varphi(m)*a$, $\varphi(m_xa)=\varphi(m)_xa$.
For any $x,y\in H^*$, $a,b\in \mathcal A$ and $m\in M^{\mathcal A}$, we have $(\varphi(m)_xa)_yb=\varphi(m)_{x_{(2)}}(a_{yx_{(-1)}}b)\in Im(\varphi)$. Hence $Im(\varphi)$ is a submodule of $N^{\mathcal A}$.

If $\varphi:M^{\mathcal A}\to N^{\mathcal A}$ is an isomorphism from $\mathcal A$ pseudomodule $M^{\mathcal A}$ to $\mathcal A$ pseudomodule $N^{\mathcal A}$, then $\varphi^{-1}:N^{\mathcal A}\to M^{\mathcal A}$ is also an isomorphism. In fact, for any $a\in \mathcal A$ and $n\in N^{\mathcal A}$, $(1^{\otimes 2}\otimes_H\varphi^{-1})(n*a)=(1^{\otimes 2}\otimes_H\varphi^{-1})(\varphi\varphi^{-1}(n)*a)=(1^{\otimes 2}\otimes_H\varphi^{-1}\varphi)(\varphi^{-1}(n)*a)=\varphi^{-1}(n)*a$. These arguments show that  the following result holds.

\begin{lemma} (Schur's Lemma) If $M^{\mathcal A}$ is a simple right $\mathcal A$ pseudomodule, then $D:=End(M^{\mathcal A}):=\{\varphi|\varphi$ is an endomorphism of the right $\mathcal A$-pseudomodule $M\}$ is a division ring.
\end{lemma}

A right $\mathcal A$ pseudomodule $M$ is said to be faithful if $M*a=0$ only if $a=0$. Every pseudomodule can be regarded as a faithful pseudomodule over a proper associative $H$-pseudoalgebra.

\begin{proposition}  Let $N$ be a submodule of a right $\mathcal A$-pseudomodule $M$. Then  (1) $r_M(m):=\{a\in \mathcal A|m*a=0\}$ is a right ideal of $\mathcal A$ for any $m\in M$.
(2) $\mathfrak{A}_N(M):=\{a\in \mathcal A|m*a\in H^{\otimes 2}\otimes_HN,\forall m\in M\}$ is an ideal  of $\mathcal A$. In particular $\mathfrak{A}(M):=\mathfrak{A}_0(M)$ is an ideal of $\mathcal A$. Furthermore,  $M$ is a faithful $\mathcal A/\mathfrak{A}_0(M)$ pseudomodule.
\end{proposition}

\begin{proof}  It is obvious.
\end{proof}

 A right $\mathcal A$ pseudomodule is said to be semisimple if it is a sum of simple $\mathcal A$ pseudomodules. It is easy to prove that any unital right semisimple $\mathcal A$ pseudomodule is a direct sum of simple $\mathcal A$ pseudomodule, and any sub-pseudomodule of a unital semisimple pseudomodule is a direct summand.

\begin{lemma} \label{lem27}Let $M$ be a semisimple unital $\mathcal A$ pseudomodule and  $R=End(M^{\mathcal A})$. If  $M$ is a finitely generated left $H\otimes R$ module and $\mathcal A'=Cend_R(M)$, then any  $\mathcal A$ sub-pseudomodule $N$ of $M$ is also an $\mathcal A'$ pseudomodule. For any $n\geq 1$, let $S=M_n(R)$. Then $End(M{}^{(n)\mathcal A})$ is the set of maps $(m_1,m_2,\cdots, m_n)\mapsto (v_1,v_2,\cdots, v_n)$, where $v_i=\sum\limits_{j=1}^nr_{ij}m_j$, $r_{ij}\in R$. Moreover, for any $a'\in \mathcal A'$, the map $(m_1, m_2,\cdots,m_n)\mapsto (m_1*a',m_2*a',\cdots, m_n*a')$ belongs to  $Cend_S(M^n)$.
\end{lemma}
\begin{proof} For any   submodule  $N$ of $M$,  there is an idempotent $e\in R$ such that $N=eM$. Let $a'\in \mathcal A'$. Then $N*a'=(1^{\otimes 2}\otimes_He)M*a'$. This implies that $N$ is an $\mathcal A'$ pseudomodule.

Let $f\in End(M^{(n)\mathcal A})$. Then $f(0,\cdots,0, m_i,0,\cdots,0)=(m_{1i},m_{2i},\cdots,m_{ni})$ for any $m_i\in M$. This gives   maps  $m_i\mapsto m_{ji}$ in $R$ for $1\leq j\leq n$. Denote these maps  by $r_{ji}$. Then $f(m_1,m_2,\cdots, m_n)=(\sum\limits_{j=1}^nr_{1j}m_j, \sum\limits_{j=1}^nr_{2j}m_j,\cdots,\sum\limits_{j=1}^nr_{nj}m_j)$. Since for any  $a'\in \mathcal A'$, we have 
$(f(m_1,m_2,\cdots,m_n))*a'=(1^{\otimes 2}\otimes_Hf)(m_1*a',m_2*a',\cdots,m_n*a')$. Thus $a'$ induce a pseudolinear map in  $Cend_S(M^{n})$
\end{proof}

For any subset $N$ of   a right $\mathcal A$ pseudomodule  $M$, let  $N\circ\mathcal A$ be the $H$-submodule of $M$ generated by 
the set $\{n_{x}a|a\in\mathcal A,x\in H^*,n\in N\}$. It is easy to check that $N\circ\mathcal A$ is a sub-pseudomodule of the $\mathcal A$ pseudomodule $M$.

\begin{proposition}\label{prop316}Let $M$ be a semisimple unital right $\mathcal A$ pseudomodule and $R=End(M^{\mathcal A})$.  If  $M$ is a finitely generated left $H\otimes R$ module and $\mathcal A'=Cend_R(M)$, then for any $a'\in \mathcal A'$, any $m_1,m_2,\cdots,m_k\in M$ and $x\in H^*$, there are $a_i\in \mathcal A$ and $x_i\in H^*$  such that $(m_1,m_2,\cdots, m_k)_xa'=\sum\limits_i (m_1,m_2,\cdots,m_k)_{x_i}a_i$.
\end{proposition}

\begin{proof} Since $M$ is a semisimple unital $\mathcal A$ pseudomodule, $M^{(k)}$ is a semisimple unital $\mathcal A$-pseudomodule. So $N=(m_1,m_2,\cdots,m_k)\circ \mathcal A$ is a submodule of $M^{(k)}$.
Hence $N$ is also an $\mathcal A'$ pseudomodule by Lemma \ref{lem27}. Since $M^{(k)}$ is unital, $(m_1,m_2,\cdots,m_k)\in N$. Thus $(m_1,m_2,\cdots,m_k)_xa'\in N$. Consequently, there exist $a_i\in\mathcal  A, y_i\in H^*$ and $h_i\in H$ such that $(m_1,m_2,\cdots,m_k)_xa'=\sum\limits_ih_i((m_1,m_2,\cdots,m_k)_{y_i}a_i)=\sum\limits_i(m_1,m_2,\cdots,m_k)_{y_ih_{i(1)}}h_{i(2)}a_i$.
\end{proof}

From Proposition \ref{prop316}, we get $(m_1,m_2,\cdots,m_k)*a'\in H^{\otimes 2}\otimes_H(m_1,m_2,\cdots,m_k)\circ\mathcal A $ and $(m_1,m_2,\cdots,m_k)\circ \mathcal A'\subseteq (m_1,m_2,\cdots,m_k)\circ\mathcal  A$. If $M$ is a faithful $\mathcal A$ pseudomodule, that is $r(M):=\{a\in \mathcal A|M*a=0\}=0$, then $\mathcal A\subseteq\mathcal  A'$ and $(m_1,m_2,\cdots,m_k)\circ \mathcal A'=(m_1,m_2,\cdots,m_k)\circ \mathcal A$.  

\begin{theorem}(Density theorem) Let $M$ be a simple unital  right $\mathcal A$ pseudomodule. Then  $M=H\otimes D^{(I)}$ for some index set $I$, where $D=End(M^{\mathcal A})$ is a division ring. For any linearly independent  $m_1,m_2,\cdots,m_k$ over $D$,  $M^{(k)}=(m_1,m_2,\cdots,m_k)\circ \mathcal A$.

Further assume that $M$ is a simple faithful unital $\mathcal A$ pseudomodule and $M=H\otimes D^n$. Then $\mathcal A\subseteq H\otimes H\otimes M_n(D)$,  $M^{(k)}=(m_1,m_2,\cdots,m_k)\circ (H\otimes H\otimes M_n(D))=(m_1,m_2,\cdots,m_k)\circ \mathcal A$ for any $1\leq k\leq n$, where $m_1,m_2,\cdots,m_k$ are linearly independent  over $D$. 
\end{theorem}

\begin{proof} First, let us prove that $M$ is a free left $H$-module. Let $1$ be an identity of $\mathcal A$ and $M_1=\{m_x1|m\in M,x\in X\}$. Note that  $(m_y1)_x1=m_{y_{(2)}}(1_{xy_{(-1)}}1)=\sum\limits_i m_{Sh_i y}(1_{xSx_i}1)$ for any $x,y\in X$.  Since $1_{xSx_i}1=\langle xSx_i,1\rangle=\langle x,1\rangle \langle x_i,1\rangle=\delta_{1,i}\delta_{1,x}$, $(m_y1)_x1=\delta_{1,x}m_y1$.
 If $m*1=\sum\limits_iSh_i\otimes 1\otimes_Hm_{x_i}1$ for
 $m\in M$, then $(m-\sum\limits_{i}h_i(m_{x_i}1))*1=\sum\limits_{j}Sh_j\otimes 1\otimes_Hm_{x_j}1-
\sum\limits_{i,j}h_iSh_j\otimes1\otimes_H(m_{x_i}1)_{x_j}1=0$. Hence $m=\sum\limits_ih_i(m_{x_i}1)$ and $M=HM_1$. If $\sum\limits_{i=1}^ng_i(m_i)_{y_i}1=0$ for some linearly independent elements $g_i\in H$, then $0=\sum\limits_{i=1}^n(g_i(m_i)_{y_i}1)*1=\sum\limits_{i=1}^n\sum\limits_jg_iSh_j\otimes1\otimes_H((m_i)_{y_i}1)_{x_j}1=
\sum\limits_{i=1}^ng_i\otimes1\otimes_H(m_i)_{y_i}1$. Thus $(m_i)_{y_i}1=0$ and $M$ is a free left $H$-module.

Next, we prove that $M_1$ is a left $End(M^{\mathcal A})$-module. Observe that $M_1=M^0$, where $M^0=\{m\in M|m*1=1\otimes 1\otimes_Hm\}$. In fact, if $m\in M^0$, then
$m_11=m$ and $m\in M_1$. Conversely, if $m\in M$, then $(m_11)_x1=m_1(1_x1)=\delta_{1,x}m_11$ for any $x\in X$. Thus $(m_11)*1=1\otimes 1\otimes _Hm_11$. Consequently $M_1=M^0$. Hence $\phi(m)\in M^0$ for any $m\in M^0$ and $\phi\in End(M^{\mathcal A})$.

Since $M^0$ can be viewed as a left vector space over the division ring $End(M^{\mathcal A})$, we have $M\simeq H\otimes D^{(I)}$ for some index set $I$, where $D=End(M^{\mathcal A})$.

Let $\varphi:D^{(I)}\to D^{(I)}$ be an endomorphism of the left vector space $_DM$ such that $\varphi(m_i)=m_i'$ for $1\leq i\leq k$. Then $(1\otimes m_i)*(1\otimes1\otimes_H\varphi)=(1\otimes 1)\otimes_H(1\otimes m_i')$. Hence $(m_1',m_2',\cdots,m_k')=\sum\limits_ih_i((m_1,m_2,\cdots,m_k)_{x_i}a_i)$
for some $a_i\in \mathcal A, x_i\in H^*$ and $h_i\in H$. Therefore $M^{k}=(m_1,m_2,\cdots m_k)\circ (H\otimes H\otimes End_D(D^{(I)}))\subseteq (m_1,m_2,\cdots,m_k)\circ \mathcal A'\subseteq (m_1,m_2,\cdots, m_k)\circ \mathcal A\subseteq M^k$ by Proposition \ref{prop316}. Consequently, $M^k=(m_1,m_2,\cdots,m_k)\circ (H\otimes H\otimes End_D(D^{(I)}))=(m_1,m_2,\cdots,m_k)\circ \mathcal A$.

If $M=H\otimes D^n$, then $\mathcal A' =Cend_D(M)=H\otimes H\otimes M_n(D)$ by Corollary  \ref{cor310}. Thus $\mathcal A\subseteq H\otimes H\otimes M_n(D)$.
\end{proof}

\begin{corollary}If $\mathcal A$ is a unital associative $H$-pseudoalgebra such that $\mathcal A$ is a right simple pseudomodule over itself, then $\mathcal  A=H\otimes D$ for some division algebra $D$.
\end{corollary}

\subsection{Functors  among categories of hybrid modules} 
Let $M,N$ be two $R$-$\mathcal A$ hybrid bimodules and  $\varphi:M\to N$  be a left $H\otimes R$ module homomorphism. Then $\varphi$ is said to be a homomorphism from the $R$-$\mathcal A$ hybrid bimodule $M$ to the $R$-$\mathcal A$ hybrid bimodule $N$ if $$(1\otimes 1\otimes_H\varphi)(m*a)=\varphi(m)*a$$ for any $m\in M$ and $a\in\mathcal A$. All $R$-$\mathcal A$ hybrid bimodules form a category ${}_R\mathcal{M}od^\mathcal A$ with the homomorphisms of $R$-$\mathcal A$ hybrid bimodules. Similarly, we can define the category ${}^\mathcal A\mathcal{M}od_R$ of all $\mathcal A$-$R$ hybrid bimodules. 

The category $_R\mathcal{M}od^{\mathcal A}$  is a subcategory $_{H\otimes R}\mathcal{M}od$ of all left $H\otimes R$ modules.  Let $_{H\otimes R'}N_R$ be an $H\otimes R'$-$R$ bimodule and $M$ be a  $R$-$\mathcal A$ hybrid bimodule. Then $N\otimes_RM$ is a right $\mathcal A$ pseudomodule with the action given by
$$(n\otimes_Rm)*a=\sum\limits_j1\otimes g_j\otimes_H(n\otimes _Rm_j),$$ where $m*a=\sum\limits_j1\otimes g_j\otimes_Hm_j$. It is easy to check that $N\otimes_RM$ is a $R'$-$\mathcal A$ hybrid bimodule. Thus $\ N\otimes_R- \ $ determines a functor from $_R\mathcal{M}od^{\mathcal A}$ to $_{R'}\mathcal{M}od^{\mathcal A}$. If there are two bimodules  $_{H\otimes R'}N_R$  and $_{H\otimes R}N'_{R'}$ such that $N_R$, $_{R'}N,$ $_RN'$ and $N'_{R'}$ are progenerators, then $_R\mathcal{M}od^{\mathcal A}$ is equivalent to $_{R'}\mathcal{M}od^{\mathcal A}$. Since any vector space can be regard as a left $H$-module with the trivial action,  and the functor  $N\otimes_R-$ is a tensor product over algebra $R$ and is not over $H\otimes R$, we have $_R\mathcal{M}od^{\mathcal A}$ is equivalent to $_{R'}\mathcal{M}od^{\mathcal A}$ provided that $R$ and $R'$ are Morita equivalent. Consequently, $\mathcal{M}od^{\mathcal A}$ is equivalent to $_{M_n({\bf k})}\mathcal{M}^{\mathcal A}$ for any $n\geq 2$.

\begin{proposition}\label{prop318}
Let $_{H\otimes R'}M_R$ be a $H\otimes R'$-$R$  bimodule and $_{H\otimes R}N_T$ be a $H\otimes R$-$T$ bimodule. If $_{R'}M^{\mathcal A}$ is a  $R'$-$\mathcal A$ hybrid bimodule such that $(mr)*a=(m*a)(1^{\otimes 2}\otimes_Hr)$  for any $m\in M$, $a\in \mathcal A$ and $r\in R$, then the  $R'$-$T$ bimodule $M\otimes _RN$ is a $R'$-$\mathcal A$ hybrid bimodule. Moreover, $(m\otimes_Rnt)*a=((m\otimes_Rn)*a)(1^{\otimes 2}\otimes_Ht)$ for any $m\in M$, $n\in N$, $a\in \mathcal A$ and $t\in T$.
\end{proposition}

\begin{proof}Since $r(hn)=h(rn)$ for any $h\in H$, $r\in R$ and $n\in N$,  we have $M\otimes_RN$ is a left $H$-module with the action given by $h(m\otimes_Rn):=h_{(1)}m\otimes_Rh_{(2)}n$. Further, if $r'\in R'$, then $r'(h(m\otimes_Rn))=(r'h_{(1)}m)\otimes_Rh_{(2)}n=h((r'm)\otimes_Rn)=r'h(m\otimes_Rn)$.  Hence $M\otimes_RN$ is a left $H\otimes R'$ module.
Let $a\in \mathcal A$ and $m\otimes_R n\in M\otimes_RN$. Define $(m\otimes_Rn)*a:=\sum\limits_j1\otimes g_j\otimes_H(m_j\otimes_Rn)$, where  $m*a=\sum\limits_j1\otimes g_j\otimes_Hm_j$. Then $((hr')(m\otimes _Rn))*a=\sum\limits_j1\otimes g_jh_{(-2)}\otimes_Hr'(h_{(1)}m)\otimes_Rh_{(3)}n=(h\otimes 1\otimes_Hr')((m\otimes_Rn)*a)$ for any $h\in H$, $r'\in R'$, $m\otimes n\in M\otimes _RN$ and $a\in \mathcal A$. Thus $M\otimes_RN$ is a  $R'$-$\mathcal A$ hybrid bimodule.
It is easy to check that $(mr\otimes_R n)*a=(m\otimes_R rn)*a$ and $(m\otimes_Rnt)*a=((m\otimes_Rn)*a)(1^{\otimes 2}\otimes_Ht)$. \end{proof}

Let $_{R'}\mathcal{M}od^{\mathcal A}_R$ be a subcategory of $_{R'}\mathcal{M}od^{\mathcal A}$ consisting of all  $R'$-$\mathcal A$ hybrid bimodules $M$, which are also right $R$-modules satisfying
$$(mr)*a=(m*a)(1^{\otimes 2}\otimes_Hr)$$ for all $m\in M$, $r\in R$ and $a\in\mathcal A$. Then an $H\otimes R$-$T$ bimodule $_{H\otimes R}N_T$ determines a functor $_-\otimes_RN$ from $_{R'}\mathcal{M}od^{\mathcal A}_R$ to $_{R'}\mathcal{M}od^{\mathcal A}_T$ by Proposition \ref{prop318}.

\begin{proposition}Let $N$ be a finitely generated free left $R$-module and  $_{H\otimes R}N_{R'}$ be a  $H\otimes R$-$R'$ bimodule. If  $M$ is a $R$-$\mathcal A$ hybrid bimodule, then $Hom_R(N,M)$ is a $R'$-$\mathcal A$ hybrid bimodule, where $(h\phi)(n)=h_{(1)}\phi(h_{(-2)}n)$, $(r\phi)(n)=\phi(nr)$ and $(\phi*a)(n)=\phi(n)*a.$
\end{proposition}

\begin{proof}Since $N$ is a free $R$ module of finite rank, the map $\Phi: Hom_R(N,R)\otimes_R M\to Hom_R(N,M)$, $\sum\limits_i \phi_i\otimes_R m_i\mapsto \Phi(\sum\limits_i \phi_i\otimes_R m_i)$ is an isomorphism of vector spaces over ${\bf k}$,  where $\Phi(\sum\limits_i \phi_i\otimes_R m_i)(n):=\sum\limits_i\phi_i(n)m_i$ for any $\sum\limits_i \phi_i\otimes_R m_i\in Hom_R(N,R)\otimes_RM$ and $n\in N$.  For any $r'\in R'$, $\Phi(r'\sum\limits_i\phi_i\otimes_R m_i)(n)=\sum\limits_i(r'\phi_i)(n)m_i=\sum\limits_i\phi_i(nr')m_i=\Phi(\sum\limits_i\phi_i\otimes_R m_i)(nr')=(r'\Phi(\sum\limits_i\phi_i\otimes_R m_i))(n)$. Hence $\Phi$ is an isomorphism of left $R'$ modules. Let $h\in H$ and $\sum\limits_i \phi_i\otimes_R m_i\in Hom_R(N,R)\otimes _RM$. Then $\Phi(\sum\limits_ih_{(1)}\phi_i\otimes_R h_{(2)}m_i)(n)=\sum\limits_i(h_{(1)}\phi_i)(n)(h_{(2)}m_i)=\sum\limits_i\phi_i(h_{(-1)}n)(h_{(2)}m_i)=h_{(1)}(\sum\limits_i\phi(h_{(-2)}n)m_i)=(h\Phi(\sum\limits_i\phi_i\otimes_Rm_i))(n).$ Hence $\Phi$ is also an isomorphism of left $H$ modules.

Define $(\sum\limits_j \phi_j\otimes_R m_j)*a:=\sum\limits_{i,j}1\otimes g_{ij}\otimes_H (\phi_i\otimes _Rm_{ij})$, where $m_i*a=\sum\limits_j1\otimes g_{ij}\otimes_Hm_{ij}$. Then $Hom_R(N,R)\otimes_RM$ becomes a right $\mathcal A$ pseudomodule. For $r'\in R$, $(r'\sum\limits_i \phi_i\otimes _Rm_i)*a=\sum\limits_j1\otimes g_{ij}\otimes_H(r'\phi_i\otimes _Rm_{ij})=(1^{\otimes 2}\otimes_Hr')((\sum\limits_i\phi_i\otimes _Rm_i)*a).$ Thus $Hom_R(N,R)\otimes_RM$ is a left $R'$ and right $\mathcal A$ hybrid bimodule.

Now let $\{e_1,e_2,\cdots,e_n\}$ is a basis of $N$ as a left $R$ module and $\varphi\in Hom_R(N,M)$. Then $\Phi(\sum\limits_{i=1}^ne_i^*\otimes_R \varphi(e_i))=\varphi.$
Assume that $\varphi(e_i)*a=\sum\limits_j1\otimes g_{ij}\otimes_Hm_{ij}$. Then $(1^{\otimes 2}\otimes_H\Phi)((\sum\limits_{i=1}^ne^*_i\otimes_R \varphi(e_i))*a)=\sum\limits_{i,j}1\otimes g_{ij}\otimes_H\Phi(e_i^*\otimes _Rm_{ij}),$ where $\sum\limits_j1\otimes g_{ij}\otimes_H \Phi(e_i^*\otimes _Rm_{ij})(e_k)=\sum\limits_{ij}\delta_{ik}\otimes g_{ij}\otimes_Hm_{ij}=\delta_{ik}\varphi(e_i)*a.$ Thus $(1^{\otimes 2}\otimes_H\Phi)((\sum\limits_{i=1}^ne^*_i\otimes_R \varphi(e_i))*a)=\varphi*a=\Phi(\sum\limits_ie_i^*\otimes _R\varphi(e_i))*a.$  Since $(\varphi*a)(e_i)=(\Phi(\sum\limits_k\phi_k\otimes_Rm_k)*a)(e_i)= \sum\limits_j1\otimes g_{ij}\otimes_H m_{ij}=\varphi(e_i)*a$ for each $i$, $(\varphi*a)(n)=\varphi(n)*a$ for any $n\in N$.
Hence $Hom_R(N,M)$ is a  $R'$-$\mathcal A$ hybrid bimodule and $\Phi$ is an isomorphism of $R'$-$\mathcal A$ hybrid bimodules.
\end{proof}

\section{bipseudomodule}
In Section 3, we give a criterion to detemine that two categories $_R\mathcal{M}od^{\mathcal A}$ and $_{R'}\mathcal{M}od^{\mathcal A}$ are equivalent for a fixed associative $H$-pseudoalgebra $\mathcal A$. In this section, we study when $_R\mathcal{M}od^{\mathcal A}$ and $_R\mathcal{M}od^{\mathcal A'}$ are equivalent for the same $R$. 
For any two right  $\mathcal A$ pseudomodules  $M$ and $N$, let $$Chom(M^{\mathcal A},N^{\mathcal A}):=\{\varphi\in Chom(M,N)|\varphi*(m*a)=(\varphi*m)*a,\forall a\in \mathcal A,m\in M\}.$$ It is easy to check that $h\varphi\in Chom(M^{\mathcal A},N^{\mathcal A})$ for any $h\in H$ and $\varphi\in Chom(M^{\mathcal A},N^{\mathcal A})$, that is, $Chom(M^{\mathcal A},N^{\mathcal A})$ is a submodule of $Chom(M,N)$ as left $H$-modules.

\begin{proposition}\label{pro41} If $\phi:M'\to M$ and $\psi:N\to N'$ are homomorphisms of $\mathcal A$ pseudomodules. Then $\phi\times \psi: Chom(M^{\mathcal A}, N^{\mathcal A})\to Chom(M'{}^{\mathcal  A},N'{}^{\mathcal A})$, $\alpha\mapsto (1^{\otimes 2}\otimes_H\psi)\alpha\phi$ is a homomorphism of left $H$-modules.  Hence $Chom(\oplus_{i=1}^nM_i^{\mathcal A},\oplus_{j=1}^mN_j^{\mathcal A})\simeq \oplus_{i=1,j=1}^{n,m}Chom(M_i^{\mathcal A},N_j^{\mathcal A})$.
\end{proposition}

\begin{proof}
For any $\alpha \in Chom(M^{\mathcal A},N^{\mathcal A})$, $h\in H$ and $m\in M$, we have $(1^{\otimes 2}\otimes_H \psi)\alpha\phi(hm)=(1^{\otimes 2}\otimes_H\psi) \alpha (h\phi(m))=(1\otimes h\otimes_H1)(1^{\otimes 2}\otimes_H\psi)\alpha(\phi(m)).$ This means that $(1^{\otimes 2}\otimes_H\psi)\alpha\phi\in Chom(M',N')$.
For any $m\in M'$ and $a\in \mathcal A$, $(1^{\otimes 2}\otimes_H\psi)(\alpha\phi)(m*a)=(1^{\otimes 2}\otimes_H\psi)(\alpha(\phi(m)*a))=\sum\limits_iSh_j\otimes Sh_i\otimes 1\otimes_H\psi(\alpha_{x_j}(\phi(m)_{x_i}a))=\sum\limits_{i,j}Sh_jSh_{i(1)}\otimes Sh_{i(2)}\otimes 1\otimes_H\psi((\alpha_{x_j}\phi(m))_{x_i}a)=\sum\limits_{i,j}Sh_jSh_{i(1)}\otimes Sh_{i(2)}\otimes 1\otimes_H\psi((\alpha_{x_j}\phi(m)))_{x_i}a=(((1\otimes 1\otimes_H\psi)\alpha\phi)(m))*a$.  Thus $(1^{\otimes 2}\otimes_H\psi)\alpha\phi\in Chom(M'{}^{\mathcal A},N'{}^{\mathcal A})$.

Let $\lambda_j:N_j\to \oplus_{j=1}^m N_j$ (resp.  $\pi_j:\oplus_{j=1}^mN_j\to N_j$) be the canonical embedding (resp. projective) homomorphisms of right $\mathcal A$ pseudomodules of  the direct sum $\oplus_{j=1}^m N_j$. Assume that $\lambda'_i:M_i\to \oplus_{i=1}^n M_i $ (resp.  $\pi_i':\oplus_{i=1}^nM_i\to M_i$) be the canonical embedding (resp. projective) homomorphisms of right $\mathcal A$ pseudomodules of  the direct sum $\oplus_{i=1}^n M_i$. Then $\pi_{i,j}:Chom(\oplus_{i=1}^n M_i^\mathcal A,\oplus_{j=1}^m N_j^\mathcal A)\to Chom(M_i^\mathcal A,N_j^\mathcal A)$, $\alpha \mapsto (1^{\otimes2}\otimes_H\pi_j)\alpha\lambda_i'$, are homomorphisms of left $H$-modules. Similarly, $\lambda_{i,j}:Chom(M_i^\mathcal A,N_j^\mathcal A)\to Chom(\oplus_{i=1}^nM_i^\mathcal A,\oplus_{j=1}^mN_j^\mathcal A)$, $\alpha\mapsto (1^{\otimes 2}\otimes_H\lambda_j)\alpha\pi_i'$ are homomorphisms of left $H$-modules. It is easy to check that $\pi_{i,j}\lambda_{i',j'}=\delta_{i,i'}\delta_{j,j'}$ and $\sum\limits_{i,j}\lambda_{i,j}\pi_{i,j}=1$. Hence $Chom(\oplus_{i=1}^nM_i^{\mathcal A},$ $\oplus_{j=1}^mN_j^{\mathcal A})\simeq \oplus_{i=1,j=1}^{n,m}Chom(M_i^{\mathcal A},N_j^{\mathcal A})$.
\end{proof}

For any  two left $\mathcal A$ pseudomodules $M$ and $N$, let  $HomC(^{\mathcal A}M,{}^{\mathcal A}N):=\{\phi\in HomC(M,N)|(a*m)*\phi=a*(m*\phi)\}$, $EndC(^{\mathcal A}M):=HomC(^{\mathcal A}M,{}^{\mathcal A}M)$.

\begin{lemma}\label{lem41}Let $M$ be a right $\mathcal A$ pseudomodule. Then  $H\otimes M$ is a left $H$ module with the action given by $h'(h\otimes m)=(h'h)\otimes m$ for $h'\in H$ and $h\otimes m\in H\otimes M$.  For any $t=\sum\limits_{i}Sh_i\otimes m_i\in H\otimes M$, define $l_t(a):=\sum\limits_{i,j}Sh_iSh_j\otimes 1\otimes _Hm_{ix_j}a$ for any $a\in\mathcal A$. Then $l_t\in Chom(\mathcal A^\mathcal A,M^\mathcal A)$ and $l_t(a)=0$ if and only if $(\sum\limits_iSh_im_i)*a=0$.

Assume further that $\mathcal A$ is a unital associative $H$-pseudoalgebra  with identity $1$. Then 
$$Chom({\mathcal A}^{\mathcal A},M^{\mathcal A})\simeq  \{\sum\limits_i h_{i(1)}\otimes h_{i(2)}m_i|\sum\limits_j1_{x_i}(Sh_jm_j)=m_i, \forall i\},$$ which is a submodule of the left $H$-module $H\otimes M$. If  $M$ is a unital right $\mathcal A$ pseudomodule, then $Chom(\mathcal A^\mathcal A,M^\mathcal A)\simeq H\otimes M$.

In particular, $H\otimes \mathcal A$ is an associative $H$-pseudoalgebra with pseudoproduct 
\begin{eqnarray}\label{equ31}(h\otimes a)*(h'\otimes b)=\sum\limits_ihSh_i\otimes h'\otimes_H1\otimes  a_{x_i}b\end{eqnarray}
for $h,h'\in H$ and $a,b\in\mathcal A$.  

Let $Y$ be an $H$ bimodule commutative associative  algebra such that it is both left and right $H$-differential algebra. Then the annihilator algebra $\mathcal A_Y(H\otimes \mathcal A)\simeq Y\otimes\mathcal A$ as associative algebras, $Y\otimes \mathcal A$ is the ordinary associative algebra defined in  Proposition \ref{prop2.9}.
\end{lemma}

\begin{proof} For any $t=\sum\limits_iSh_i\otimes m_i\in H\otimes \mathcal A$, let $l_t(a):=\sum\limits_{i,j}Sh_iSh_j\otimes 1\otimes_H (m_i)_{x_j}a$ for any $a\in\mathcal A$.
If $h\in H$ and $a\in\mathcal A$, then $l_t(ha)=\sum\limits_{i}(Sh_i\otimes 1\otimes_H1)\sum\limits_j(Sh_j\otimes 1\otimes_Hm_{ix_j}(ha))=\sum\limits_{i,j}Sh_iSh_j\otimes h\otimes_Hm_{ix_j}a.$ Thus $l_t\in Chom(\mathcal A,M)$. It is easy to check that $l_t\in Chom(\mathcal A^\mathcal A,M^\mathcal A)$ and $l_{ht}=hl_t$ for any $h\in H$. If $l_t(a)=0$, then $\sum\limits_{i,j}Sh_iSh_j\otimes 1\otimes_H (m_i)_{x_j}a=0.$ Thus $l_{tx_k}(a)=\sum\limits_{i,j}\langle x_k, h_jh_i\rangle m_{ix_j}a=\sum\limits_im_{i x_kSh_i}a=(\sum\limits_iSh_im_i)_{x_k}a.$
Consequently, $l_t(a)=0$ if and only if $(\sum\limits_iSh_im_i)*a=0$.

Next, let us assume that $\mathcal A$ is a unital associative $H$-pseudoalgebra with indenity $1$.
Define $\Phi:Chom(\mathcal A^{\mathcal A},M^{\mathcal A})\to H\otimes M$ by $\Phi(\alpha):=\sum\limits_i Sh_i\otimes \alpha_{x_i}{1}$.  Then $\sum\limits_j(\alpha_{x_j}1)_{x_iSh_j}1=\sum\limits_j\alpha_{x_{j(2)}}(1_{(x_iSh_j)x_{j(-1)}}1)$ $=\sum\limits_{j,p}\alpha_{x_p}(1_{(x_iSh_j)(S(x_jSh_p))}1)$. For any $h\in H$, we have $$\sum\limits_j\langle (x_iSh_j)(S(x_jSh_p)),h\rangle=\sum\limits_j\langle x_i, h_{(1)}h_j\rangle\langle x_j ,h_{(-2)}h_p\rangle=\langle x_i,h_{(1)}h_{(-2)}h_p\rangle=\varepsilon(h)\langle x_i,h_p\rangle.$$
 Hence $\sum\limits_j(\alpha_{x_j}1)_{x_iSh_j}1=\alpha_{x_i}1.$  If $m_j=m_{j1}1$ for each $j$, then  $\sum\limits_j(m_j)_{x_iSh_j}1=\sum\limits_j(m_j{}_11)_{x_iSh_j}1=\sum\limits_jm_j{}_1(1_{x_iSh_j}1)$. Since $1*1=1\otimes 1\otimes_H1$, $1_{x_iSh_j}1=\langle x_iSh_j,1\rangle=\langle x_i, h_j\rangle =\delta_{ij}1$. Thus $\sum\limits_j(m_j)_{x_iSh_j}1=m_i$.  
As $(h\alpha)*1=\sum\limits_ihSh_i\otimes1\otimes_H\alpha_{x_i}{1}$,  we have 
 $\Phi(h\alpha)=h\Phi(\alpha)$.

If $\Phi(\alpha)=\Phi(\beta)$,  then $\sum\limits_iSh_i\otimes 1\otimes _H\alpha_{x_i}1=\sum\limits_i Sh_i\otimes1\otimes_H\beta_{x_i}1$.
 Thus $\alpha({1}_1a)=\sum\limits_{j}Sh_j\otimes 1\otimes_H\alpha_{x_j}({1}_1a)=\sum\limits_{j}Sh_j\otimes 1\otimes_H(\alpha_{x_{j(2)}}1)_{x_{j(1)}}a=\sum\limits_{i,j}Sh_j\otimes 1\otimes_H(\alpha_{x_{i}}1)_{Sh_i x_{j}}a
 =\sum\limits_{j}Sh_j\otimes 1\otimes_H\beta_{x_j}({1}_1a)=\beta({1}_1a).$  Since $\mathcal A$ is unital, $\alpha({1}_1a)=\alpha({a})= \beta({a})$ and $\alpha=\beta$.  Thus 
 $\Phi$ is injective.

Conversely, let  $\sum\limits_iSh_i\otimes m_i\in H\otimes M$ satisfying $\sum\limits_j(m_j)_{x_iSh_j}1=m_i$ for each $i$ and $t=\sum\limits_iSh_im_i$. Then  
$l_{tx_k}1=\sum\limits_{i,j}\langle x_k,h_jh_i\rangle (m_i)_{x_j}1=\sum\limits_i (m_i)_{x_kSh_i}1=m_k.$ Thus
$\Phi(l_t)=\sum\limits_{i,j}Sh_i\otimes (m_j)_{x_iSh_j}1=\sum\limits_i Sh_i\otimes m_i$. Thus  $\Phi$ is an isomorphism of left $H$ modules.
 
 For $h\otimes a,h'\otimes b\in H\otimes \mathcal A$ and $c\in\mathcal A$, $(h\otimes a)*((h'\otimes b)*c)=\sum\limits_{i,j}hSh_i\otimes h'Sh_j\otimes 1\otimes_Ha_{x_i}(b_{x_j}c)=\sum\limits_{i,j}hSh_iSh_{j(1)}\otimes h'Sh_{j(2)}\otimes_H(a_{x_i}b)_{x_j}c=(\sum\limits_{i}hSh_i\otimes h'\otimes_H1\otimes a_{x_i}b)*c$. Hence $(h\otimes a)*(h\otimes b)=\sum\limits_ihSh_i\otimes h'\otimes_H1\otimes a_{x_i}b$.
 
 Note that $\xi: \mathcal {A}_Y(H\otimes \mathcal A)=Y\otimes_HH\otimes \mathcal A\simeq Y\otimes \mathcal A$, $y\otimes_H h\otimes a\mapsto yh\otimes a$, is an isomorphism of vector spaces. For any $y\otimes h\otimes_Ha,y'\otimes h'\otimes_Hb\in Y\otimes_HH\otimes \mathcal A$, $(y\otimes h\otimes_Ha)(y'\otimes h'\otimes _Hb)=\sum\limits_i(yhSh_i)(y'h')\otimes_H 1\otimes a_{x_i}b$. Thus $\xi$ is also an isomorphism of associative algebras.
\end{proof}

\begin{lemma}\label{le41}Assume that  $H$ is an enveloping algebra of a finite-dimensional Lie algebra.
Let $M$ be a left unital $\mathcal A$ pseudomodule, where $\mathcal A$ is a unital associative $H$-pseudoalgebra  with identity $1$ such that $a_11=a$ for any $a\in\mathcal A$. Then $HomC({}^{\mathcal A}{\mathcal A},{}^{\mathcal A}M)\simeq  H\otimes M$,  where $H\otimes M$ is a left $H$ module with the action given by $h'(h\otimes m)=(h'h)\otimes m$ for $h'\in H$ and $h\otimes m\in H\otimes M$. 

Espicially, $EndC({}^\mathcal A\mathcal A)$ is a subalgebra of the $H$-pseudoalgebra $H\otimes \mathcal A$ with the pseudoproduct given by (\ref{equ31}).
\end{lemma}

\begin{proof}For any $t=\sum\limits_{i}h_i\otimes m_i\in H\otimes M$, let $r_t:\mathcal A\to H^{\otimes 2}\otimes _HM$, $a\mapsto \sum\limits_{i,j}Sh_j\otimes h_i\otimes_Ha_{x_j}m_i$. Then $r_t\in HomC({}^\mathcal A\mathcal A,{}^\mathcal AM)$.
If $H\otimes M$ is a left $H$ module with the action given by $h'(h\otimes m)=h'h\otimes m$, then $r_{ht}=hr_t$ for any $h\in H$. 

 If $M$ is a unital left $\mathcal A$ pseudomodule and $r_t=0$, then
$0=1*r_t=\sum\limits_{i,j}Sh_j\otimes h_i\otimes _H1_{x_j}m_i$ and $\sum\limits_{i,j}h_i\otimes Sh_j\otimes _H1_{x_j}m_i=0$. Applying $1\otimes (\langle 1,\ \rangle \otimes_H1)$  to the previous equation, we get $\sum\limits_ih_i\otimes1\otimes_Hm_i=0$ and $t=0$. Therefore, the map $t\mapsto r_t$ is a left $H$ module monomorphism from $H\otimes\mathcal A$ to $HomC({}^\mathcal A\mathcal A, {}^\mathcal AM)$.

Furthermore,  define $\Psi:HomC({}^\mathcal A\mathcal A,{}^\mathcal AM)\to H\otimes M$ by $\Psi(\alpha)=\sum\limits_ih_{i(1)}\otimes h_{i(-2)}(1_{x_i}\alpha)$, where $1*\alpha=\sum\limits_iSh_i\otimes 1\otimes_H1_{x_i}\alpha$.  Then $1*r_{\Psi(\alpha)}=\sum\limits_jSh_j\otimes h_{i(1)}\otimes_H1_{x_j}(h_{i(-2)}(1_{x_i}\alpha))=
\sum\limits_i(1*(h_{i(-2)}1_{x_i}))(1\otimes h_{i(1)}\otimes_H\otimes 1)=\sum\limits_{i,j}Sh_j\otimes \varepsilon(h_i)\otimes_H1_{x_j}(1_{x_i}\alpha)=\sum\limits_{j}Sh_j\otimes 1\otimes_H1_{x_j}(1_1\alpha)=\sum\limits_jSh_j\otimes 1\otimes_H1_{x_j}\alpha.$ Thus $\Psi(r_{\Psi(\alpha)})=\Psi(\alpha)$.
If $\Psi(\alpha)=\Psi(\beta)$, then $\sum\limits_ih_{i(1)}\otimes h_{i(-2)}(1_{x_i}\alpha)=\sum\limits_ih_{i(1)}\otimes h_{i(-2)}(1_{x_i}\beta)$. Thus $a*\alpha=\sum\limits_iSh_i\otimes 1\otimes_H(a_11)_{x_i}\alpha=\sum\limits_i1\otimes h_{i(1)}\otimes_Hh_{i(-2)}(a_1(1_{x_i}\alpha))=\sum\limits_i1\otimes h_{i(1)}\otimes_Hh_{i(-2)}(a_1(1_{x_i}\beta))=a*\beta$ for any $a\in\mathcal A$. Hence $\alpha=\beta$ and  $\Psi$ is injective. Since $\Psi(r_{\Psi(\alpha)})=\Psi(\alpha)$,
$r_{\Psi(\alpha)}=\alpha$ and $r:H\otimes M\to HomC({}^\mathcal A\mathcal A,{}^\mathcal AM)$, $t\mapsto r_t$  is an isomorphism.

For any $h\otimes a,h'\otimes b\in H\otimes\mathcal A$ and $c\in\mathcal A$, $(c*(h\otimes a))*(h'\otimes b)=\sum\limits_{i,j}Sh_iSh_{j(1)}\otimes hSh_{j(2)}\otimes h'\otimes_H(c_{x_i}a)_{x_j}b=\sum\limits_{i,j}Sh_i\otimes hSh_j\otimes h'\otimes_Hc_{x_i}(a_{x_j}b)=c*(\sum\limits_jhSh_j\otimes h'\otimes_H1\otimes a_{x_j}b)$.
Hence $(h\otimes a)*(h'\otimes b)=\sum\limits_ihSh_i\otimes h'\otimes_Ha_{x_i}b$ and $EndC({}^\mathcal A\mathcal A)$ is a subalgebra of the $H$-pseudoalgebra of $H\otimes\mathcal A$ with pseudoproduct given by (\ref{equ31}).
\end{proof}

\begin{theorem} \label{cor43}Let $\mathcal A$ be a unital associative $H$-pseudoalgebra with identity $1$ such that $a_11=a$ for any $a\in \mathcal A$. 
 Then
 $\mathcal {A}$ is a left unital $H\otimes \mathcal A$ pseudomodule and $EndC(^{H\otimes \mathcal{A}} \mathcal A)\simeq \{\sum\limits_iSh_i\otimes a_i|\sum\limits_iSh_i\otimes a_i=\sum\limits_{i,j}Sh_iSh_{j(-1)}\otimes Sh_{j(2)}(1_{x_j}a_i)\}\leq H\otimes \mathcal A$, where $H\otimes \mathcal A$ is a left $H$ module with the action given by
 $h(h'\otimes a):=h'h_{(-1)}\otimes h_{(2)}a$ for $h\in H$ and $h'\otimes a\in H\otimes\mathcal A$.
 \end{theorem}

\begin{proof} From (\ref{equ31}), we can obtain that  $(1\otimes 1)*(1\otimes 1)=1\otimes 1\otimes_H1\otimes 1$ and $(1\otimes 1)*(h\otimes a)=\sum\limits_iSh_ih_{(-1)}\otimes 1\otimes_Hh_{(2)}\otimes 1_{x_i}a$ for any $h\otimes a\in H\otimes \mathcal A$. Thus $(1\otimes 1)_1(h\otimes a)=\sum\limits_i\langle 1,h_{(1)}h_i\rangle h_{(2)}\otimes 1_{x_i}a=h\otimes a$. Hence  $1\otimes 1$ is an identity of $H\otimes \mathcal A$.
It is easy to check that $\mathcal A$ is a left $H\otimes \mathcal{A}$ pseudomodule.  For $1\otimes 1\in H\otimes \mathcal A$ and $a\in \mathcal A$, $(1\otimes 1)*a=\sum\limits_iSh_i\otimes 1\otimes _H1_{x_i}a$.  Hence $(1\otimes 1)_1a=1_1a=a$ for any $a\in\mathcal A$. Thus $\mathcal {A}$ is a unital left $H\otimes \mathcal A$ pseudomodule. If $\alpha\in EndC(^{H\otimes \mathcal A}\mathcal A)$, then $((Sh_i\otimes a)*b)*\alpha=\sum\limits_{j,k} Sh_iSh_jSh_{k(1)}\otimes Sh_{k(2)}\otimes 1\otimes_H(a_{x_j}b)_{x_k}\alpha=(Sh_i\otimes a)*(b*\alpha)=\sum\limits_{j,k}Sh_iSh_k\otimes Sh_j\otimes 1\otimes_Ha_{x_k}(b_{x_j}\alpha)$ for any $a,b\in \mathcal A$. In particular, $\sum\limits_{j,k} Sh_jSh_{k(1)}\otimes Sh_{k(2)}\otimes 1\otimes_H(a_{x_j}1)_{x_k}\alpha=\sum\limits_{j,k}Sh_j\otimes Sh_k\otimes 1\otimes_Ha_{x_j}(1_{x_k}\alpha)$. Thus \begin{eqnarray}\label{eq31}\begin{array}{lll}\sum\limits_kSh_{k(1)}\otimes Sh_{k(2)}\otimes 1\otimes_H1_{x_k}\alpha&=&\sum\limits_{j,k} Sh_jSh_{k(1)}\otimes Sh_{k(2)}\otimes 1\otimes_H(1_{x_j}1)_{x_k}\alpha\\
&=&\sum\limits_{j,k}Sh_j\otimes Sh_k\otimes 1\otimes_H1_{x_j}(1_{x_k}\alpha)\\
&=&\sum\limits_{j,k}1\otimes Sh_kSh_{j(-1)}\otimes Sh_{j(-2)}\otimes_HSh_{j(3)}(1_{x_j}(1_{x_k}\alpha)).\end{array}\end{eqnarray}Applying the functor $\varepsilon\otimes 1\otimes\varepsilon \otimes_H1$ to (\ref{eq31}), we get $\sum\limits_kSh_k\otimes 1_{x_k}\alpha=\sum\limits_{j,k}Sh_jSh_{k(-1)}\otimes Sh_{j(2)}(1_{x_j}(1_{x_k}\alpha))$.
Define $\Lambda:EndC(^{H\otimes \mathcal A}\mathcal A)\to H\otimes  \mathcal A$ by $\Lambda(\alpha)=\sum\limits_i Sh_i\otimes 1_{x_i} \alpha$. If $\alpha\neq \beta$, then there is $a\in\mathcal A$ and $x_j\in X$ such that $a_{x_j}\alpha\neq a_{x_j}\beta$.
As $a_{x_j}\alpha=(a_11)_{x_j}\alpha=a_{1}(1_{x_j}\alpha)$, we have $\Lambda$ is injective. Since $\sum\limits_iSh_i\otimes 1\otimes_H1_{x_i}(h\alpha)=\sum\limits_iSh_i\otimes h\otimes_H1_{x_i}\alpha$, $\Lambda(h\alpha)=\sum\limits_iSh_ih_{(-1)}\otimes h_{(2)}(1_{x_i}\alpha)$ and $\Lambda$ is a left $H$ module homomorphism.

For any  $\sum\limits_i Sh_i\otimes a_i\in H\otimes \mathcal A$ satisfying $\sum\limits_iSh_i\otimes a_i=\sum\limits_{i,j}Sh_iSh_{j(-1)}\otimes Sh_{j(2)}(1_{x_j}a_i)$, let  $$b*\alpha=\sum\limits_{i,j}Sh_j\otimes Sh_i\otimes_Hb_{x_j}a_i$$ for any $b\in \mathcal A$. Then $\alpha\in EndC(\mathcal A)$ and $((1\otimes m)*a)*\alpha=(\sum\limits_jSh_j\otimes 1\otimes_Hm_{x_j}a)*\alpha=\sum\limits_{j,k}Sh_jSh_{k(1)}\otimes Sh_{k(2)}\otimes Sh_i\otimes_H(m_{x_j}a)_{x_k}a_i$ and $(1\otimes m)*(a*\alpha)=(1\otimes m)\sum\limits_jSh_j\otimes Sh_i\otimes_Ha_{x_j}a_i=\sum\limits_{j,k}Sh_k\otimes Sh_j\otimes Sh_i\otimes_Hm_{x_k}(a_{x_j}a_i)$ for any $1\otimes m\in H\otimes\mathcal A$. Thus $\alpha\in EndC(^{H\otimes \mathcal A}\mathcal A)$. Moreover,  $\Lambda (\alpha)=\sum\limits_{i,j}Sh_jSh_{i(-1)}\otimes Sh_{i(2)}(1_{x_j}a_i)=\sum\limits_iSh_i\otimes a_i$. Hence $\Lambda$ is an isomorphism.
 \end{proof}

\begin{example}\label{ex44} Let $\mathcal A=H\otimes A$ be a current $H$-pseudoalgebra, where $A$ is an ordinary associative algebra over ${\bf k}$ with identity. Then $(Sh_i\otimes (Sh_j\otimes a))*(1\otimes b)=Sh_iSh_j\otimes 1\otimes_H(1\otimes ab)$ for any $a,b\in A$. Thus $(Sh_i\otimes (Sh_p\otimes a))*((Sh_j\otimes(Sh_q\otimes  b))*(1\otimes c))=Sh_iSh_p\otimes  Sh_jSh_q\otimes 1\otimes_H1\otimes abc=((Sh_iSh_ph_{q(1)}h_{j(1)}\otimes 1)\otimes_H(Sh_{j(2)}Sh_{q(2)}\otimes (1\otimes ab)))*(1\otimes c)$ for any $a,b,c\in A$. Hence $H\otimes \mathcal A$ is an associative $H$-pseudoalgebra with pseudoproduct $(Sh_i\otimes (Sh_p\otimes a))*(Sh_j\otimes(Sh_q\otimes  b))=(Sh_iSh_ph_{q(1)}h_{j(1)}\otimes 1)\otimes_H(Sh_{j(2)}Sh_{q(2)}\otimes (1\otimes ab))=Sh_iSh_p\otimes Sh_jSh_q\otimes _H(1\otimes 1\otimes ab)$.

Let $n\geq 1$ be an integer and $A$ be an ordinary associative algebra. Then $H^{\otimes n}\otimes A$ is an associative $H$-pseudoalgebra with the following operation
$$(h_1\otimes h_2\otimes\cdots\otimes h_n\otimes a)*(g_1\otimes g_2\otimes\cdots\otimes g_n\otimes b):=(h_1h_2\cdots h_n)\otimes(g_1g_2\cdots g_n)\otimes_H(1^{\otimes n-1}\otimes ab)$$
and $$g_1(h_1\otimes h_2\otimes\cdots\otimes h_n\otimes a)=g_1h_1\otimes h_2\otimes\cdots\otimes h_n\otimes a$$
for all $h_i,g_i\in H$ and $a,b\in A$. The annihilator algebra $\mathcal{A}_Y(H^{\otimes n}\otimes A)\simeq Y\otimes H^{\otimes n-1}\otimes A$  with product
$(y_1\otimes h_1\otimes h_2\otimes \cdots\otimes h_{n-1}\otimes a)*(y_2\otimes g_1\otimes  g_2\otimes \cdots\otimes g_{n-1}\otimes b)=(y_1h_1h_2\cdots h_{n-1})(y_2g_1g_2\cdots g_{n-1})\otimes  1^{\otimes n-1}\otimes ab$.
\end{example}

\begin{remark}\label{rem44}From Theorem  \ref{cor43} and Proposition \ref{pro41}, we get that $Cend((\mathcal A^\mathcal A)^n)\simeq H\otimes M_n(\mathcal A)\simeq M_n(H)\otimes \mathcal A\simeq M_n(H\otimes \mathcal A)$, which  is an associative $H$-pseudoalgebra for any associative unital $H$-pseudoalgebra $\mathcal A$ and $n\geq 1$, where $(a_{ij})*(b_{ij})=(c_{ij})$ and $c_{ij}=\sum\limits_{k=1}^na_{ik}*b_{kj}$. 
\end{remark}

\begin{definition} A right $\mathcal A$ pseudomodule $M$ is said to be finitely generated if there are $m_1,m_2,\cdots,m_s\in M$ such that $M$ is equal to its $H$ submodule generated by $\{m_1{}_{x_{i_1}}a_1,m_2{}_{x_{i_2}}a_2,\cdots,$ $m_s{}_{x_{i_s}}a_s|i_1,i_2,\cdots,i_s\in I, a_i\in \mathcal A\}$. In this case, we say that $M$ is generated by $m_1,m_2,\cdots,m_s$ as a right $\mathcal A$ pseudomodule.
\end{definition}

Since $h(m_{x_i}a)=m_{h_{(1)}x_i}h_{(2)}a$, $h(m_{x_i}a)$ belongs to the vector space spanned by $\{m_{x_j}a|j\in I,a\in A\}$ over ${\bf k}$ by Proposition \ref{pr21}.
Thus the left $H$ submodule generated $\{m_1{}_{x_{i_1}}a_1,m_2{}_{x_{i_2}}a_2,\cdots,m_s{}_{x_{i_s}}a_s|i_1,$ $i_2,\cdots,i_s\in I, a_i\in \mathcal A\}$ is equal to the subspace 
spanned by $\{m_1{}_{x_{i_1}}a_1,m_2{}_{x_{i_2}}a_2,\cdots,m_s{}_{x_{i_s}}a_s|i_1,i_2,\cdots,$ $i_s\in I, a_i\in \mathcal A\}$ over ${\bf k}.$

For  any  unital associative $H$-pseudoalgebra  $\mathcal A$  with a unit $1$,  $\mathcal A^\mathcal A$ is generated by $1$ as a right $\mathcal A$ pseudomodule since $1_1a=a$ for any $a\in \mathcal A$.  Any direct sum of finite many finitely generated $\mathcal A$ pseudomodules is  again a finitely generated $\mathcal A$ pseudomodule.
Similarly to   \cite[ Lemma 10.1]{BDK} or Lemma \ref{lem34} , we get the following result.

\begin{lemma}
Let $U,V,W$  be three right $\mathcal A$ pseudomodules, and assume that $U$ is a finitely generated $\mathcal A$ pseudomodule. Then there is a unique polylinear map
$\mu: Chom(V^\mathcal A, W^\mathcal A)\boxtimes Chom(U^\mathcal A, V^\mathcal A)\to H^{\otimes 2}\otimes_H Chom(U^\mathcal A, W^\mathcal A)$,  denoted as $\mu(\alpha,\beta)=\alpha*\beta$, such that
$$(\alpha*\beta)*u=\alpha*(\beta*u)$$ in $H^{\otimes 3}\otimes_H W$ for any $\alpha\in Chom(V^\mathcal A,  W^\mathcal A), \beta\in Chom(U^\mathcal A, V^\mathcal A)$ and $ u\in U.$ In particular, $Chom(U^\mathcal A,$ $U^\mathcal A)$, denoted by $Cend(U^\mathcal A)$,  is a subalgebra of $Cend(U)$ provided that $Cend(U)$ is an associative $H$-pseudoalgebra. 
\end{lemma}

\begin{proof} If $\mu$ is well-defined, then $\mu (\phi,\psi)=\sum\limits_i Sh_i\otimes 1\otimes_H\phi_{x_i}\psi$ for $\phi\in Chom(V^\mathcal A,W^\mathcal A)$ and $\psi\in Chom(U^\mathcal A,W^\mathcal A)$, where $\{h_i\}$ and $\{x_i\}$ are dual bases in $H$ and $H^*$ respectively.
Thus one must define $\phi_x\psi\in Chom(U^\mathcal A,W^\mathcal A)$ for any $x\in H^*$. Assume that $\psi*u=\sum\limits_{\alpha}f^{\alpha}\otimes 1\otimes_Hv_{\alpha}$ for  $u\in U$ and $\phi*v_{\alpha}=\sum\limits_{\beta}g_{\alpha}^{\beta}\otimes 1\otimes_Hw_{\alpha\beta}$.
Then $\sum\limits_i\phi_{x_i}(\psi_{y(h_iSx)}u)=\sum\limits_{\alpha,\beta}\langle y,Sf^{\alpha}_{(1)}\rangle \langle x,f^{\alpha}_{(2)}Sg_{\alpha}^{\beta}\rangle w_{\alpha\beta}$.
Thus $(\phi_x\psi)*u=\sum\limits_iSh_i \otimes 1\otimes (\phi_{x}\psi)_{x_i}u=\sum\limits_{i,j}Sh_i\otimes 1\otimes_H\phi_{x_j}(\psi_{x_ih_jSx_i}u)=\sum\limits_{i,\alpha,\beta}Sh_i\otimes 1\otimes_H\langle x_i, Sf_{(1)}^\alpha\rangle\langle x_j,f^\alpha_{(2)}Sg_\alpha^\beta \rangle=\sum\limits_{\alpha,\beta}f^{\alpha}_{(1)}\otimes 1\otimes \langle x,f^{\alpha}_{(2)}Sg_{\alpha}^{\beta}\rangle w_{\alpha\beta}$. For $y,z\in H^*$ and $a\in \mathcal A$, since $$\begin{array}{lll}((\phi_x\psi)_yu)_za&=&(\phi_{x_{(2)}}(\psi_{yx_{(-1)}}u))_za\\
&=&\phi_{x_{(3)}}((\psi_{yx_{(-1)}}u)_{zx_{(-2)}}a)\\
&=&
\phi_{x_{(4)}}(\psi_{yx_{(-1)}}(u_{zx_{(-3)}x_{(2)}y_{(-1)}}a))\\
&=&\phi_{x_{(2)}}(\psi_{y_{(2)}x_{(-1)}}(u_{zy_{(-1)}}a))\\
&=&(\phi_x\psi)_{y_{(2)}}(u_{zy_{(-1)}}a),\end{array}$$  $(\phi_{x}\psi)*(u*a)=((\phi*\psi)*u)*a$ and $\phi_x\psi\in Chom(U^\mathcal A,$ $W^\mathcal A)$.

Now let $U$ be  generated by $u_1,u_2,\cdots,u_n$ as a right $\mathcal A$ pseudomodule and $\psi*u_i=\sum\limits_i f_i^{\alpha}\otimes 1\otimes_Hv_{i\alpha}$ and $\phi*v_{i\alpha}=\sum\limits_{\beta}g^{\beta}_{i\alpha}\otimes 1\otimes_Hw_{i\alpha\beta}$. If $T$ is the subspace of $H$ spanned by $f_{i(2)}^{\alpha}Sg_{i\alpha}^{\beta}$, then $T$ is  finite-dimensional and $\{x\in H^*|\langle x,T\rangle=0\}$ is a finite-codimensional subspace of $H^*$. Thus there are only finitely many $i$ such that $(\phi_{x_i}\psi)*u_j\neq 0$ for $1\leq j\leq n$. 
For any  $u=\sum\limits_pg_p(u_p{}_{y_p}a_p)\in U$,  $(\phi_{x_i}\psi)*u=\sum\limits_{p,j}Sh_j\otimes g_p\otimes_H((\phi_{x_i}\psi)_{x_{j(2)}}u_p)_{y_px_{j(1)}}a_p=\sum\limits_{p,j,q}Sh_j\otimes g_p\otimes_H((\phi_{x_i}\psi)_{x_q}u_p)_{y_pSh_qx_{j}}a_p.$
Hence there are only finitely many $i$ such that $\phi_{x_i}\psi\neq 0$. Therefore $\phi*\psi=\sum\limits_i Sh_i\otimes 1\otimes_H\phi_{x_i}\psi\in H^{\otimes 2}\otimes_H Chom(U^\mathcal A,W^\mathcal A)$.
\end{proof}

Let $Chom^{\bf r}(M^{\mathcal A},N^{\mathcal A}):=Chom^{\bf r}(M,N)\cap Chom(M^{\mathcal A},N^{\mathcal A})$ and  $Chom^{\circ}(M^{\mathcal A},N^{\mathcal A}):=
\{\phi \in Chom(M,$ $N)|\phi(M_1)=0$ for some sub-pseudomodule  $M_1$ of $M^\mathcal A $ and there exists a finitely generated $\mathcal A$ sub-pseudomodule $M_2$ of $M^\mathcal A$ such that $M=M_1+M_2\}$ for any two right $\mathcal A$ pseudomodules $M$ and $N$. Meanswhile,  let us define $Cend^{\bf r}(M^{\mathcal A}):=Chom^{\bf r}(M^{\mathcal A},$ $M^{\mathcal A})$  and $Cend^{\circ}(M^{\mathcal A}):=Chom^{\circ}(M^{\mathcal A},M^{\mathcal A})$. Similarly, we can define $HomC(^\mathcal AM,{}^\mathcal AN)$, $HomC^{\bf r}(^\mathcal AM,$ ${}^\mathcal AN)$, $HomC^\circ(^\mathcal AM,{}^\mathcal AN)$, $EndC(^\mathcal AM)$,  $EndC^{\bf r}(^\mathcal AM)$ and  $EndC^\circ(^\mathcal AM)$.

\begin{remark}Let $M^\mathcal A$ be a finitely generated right $\mathcal A$ pseudomodule. Then $Cend(M^\mathcal A)$ is an associative $H$-pseudoalgebra and $(\phi*\psi)*m=\phi*(\psi*m)$ for any $\phi,\psi\in Cend(M^\mathcal A)$ and $m\in M$. Moreover, $M$ is a $Cend(M^\mathcal A)$-$\mathcal A$ bi-pseudomodule.
\end{remark}

For example, if $\mathcal A$ is an associative $H$-pseudoalgebra with identity $1$ such that $a_11=a$ for any $a\in\mathcal A$, then $Cend(\mathcal A^\mathcal A)\simeq H\otimes \mathcal A$ and $\mathcal A$ is an $H\otimes \mathcal A$-$\mathcal A$ bipseudomodule, where $(\sum\limits Sh_i\otimes a_i)*a=\sum\limits_{i,j}Sh_iSh_j\otimes 1\otimes_H a_{ix_j}a$ for any $a\in\mathcal A$ and $\sum\limits_i Sh_i\otimes a_i\in H\otimes \mathcal A$. Furthermore, 
for any right $\mathcal A$ pseudomodule $M$,  let  $\mathcal B$ be either $Cend^{\circ}(M^\mathcal A)$ or $Cend^{\bf r}(M^\mathcal A)$. Then $M$ is a $\mathcal B$-$\mathcal A$ bi-pseudomodule.

\begin{proposition} \label{pr410}(1)
Let $M$ is a $R$-$\mathcal A$ hybrid  bimodule and $N$ is a  $\mathcal B$-$R$ hybrid bimodule. Then $M\otimes_R N$ is a  $\mathcal B$-$\mathcal A$ bi-pseudomodule with
 actions given by
$$(m\otimes_R n)*a=\sum\limits_i 1\otimes g_i\otimes_Hm_i\otimes_R n$$
and $$b*(m\otimes_R n)=\sum\limits_iSh_i\otimes 1\otimes_H(m\otimes _R (b_{x_i}n))$$
where $m*a=\sum\limits_i1\otimes g_i\otimes_Hm_i$.  In particular, $\mathcal A\otimes \mathcal B$ is an $\mathcal A$-$\mathcal B$ bi-pseudomodule.

(2) If $R$ is a commutative algebra over ${\bf k}$, then $N\otimes_R M$ is a $\mathcal B$-$\mathcal A$ bipseudomdule with the actions given by 
$$b*(n\otimes_R m)=\sum\limits_iSh_i\otimes 1\otimes_H(b_{x_i}n)\otimes_R m$$
and $$(n\otimes_R m)*a=\sum\limits_j1\otimes g_j\otimes_Hn\otimes_R m_j$$
where $m*a=\sum\limits_j1\otimes_H g_j\otimes m_j$.

(3) In the case when $R$ is a commutative ${\bf k}$ algebra, the flipping map $\sigma: M\otimes_R N\to N\otimes_R M,m\otimes_R n\mapsto n\otimes _Rm$ is an isomorphism of  $\mathcal B$-$\mathcal A$ bi-pseudomodules.
\end{proposition}

\begin{proof}Similar to Proposition \ref{prop318}, we can prove that $M\otimes_RN$ is a left $\mathcal B$ pseudomodule and a right $\mathcal A$ pseudommodule.
For $a\in \mathcal B,$  $m'\in M$, $m\in N$ and $a'\in \mathcal A$, we have $(a*(m'\otimes m))*a'=(\sum\limits_iSh_i\otimes 1\otimes_H((a_{x_i}m'\otimes m ))*a'=\sum\limits_i((Sh_i\otimes 1)\Delta\otimes 1\otimes_H1)( (a_{x_i}m')\otimes m)*a'=\sum\limits_{i,j}Sh_i\otimes 1\otimes g_j\otimes_H  (a_{x_i}m)\otimes n_j$ and
$a*((m\otimes n)*a')=\sum\limits_j(1\otimes (1\otimes g_j)\Delta\otimes_H1)a*(m'\otimes n_j)=\sum\limits_{i,j}Sh_i\otimes 1\otimes g_j\otimes_H(a_{x_i}m')\otimes n_j$. 
Thus $M\otimes_R N$ is a $\mathcal B$-$\mathcal A$ bi-pseudomodule.

The rest claims can be checked easily.
\end{proof}

\begin{proposition}\label{propo48h} Let 
 $M$ be an $\mathcal A$-$\mathcal B$ bi-pseudomodule,  $N$ be a left $\mathcal A$ pseudomodule and $P$ be a right $\mathcal B$ pseudomodule. Then

(1) $HomC({}^\mathcal AM^\mathcal B,{}^\mathcal AN)$ is a left $\mathcal B$ pseudomodule determined by $m*(b*\phi):=(m*b)*\phi$, for any $m\in M, b\in\mathcal B$ and $\phi\in Hom({}^\mathcal AM^\mathcal B,{}^\mathcal AN)$;

(2) $Chom({}^\mathcal AM^\mathcal B,P^\mathcal B)$ is a right $\mathcal A$ pseeudomodule determined by $(\phi*a)*m:=\phi*(a*m)$, for any $m\in M, a\in\mathcal A$ and $\phi\in Chom({}^\mathcal AM^\mathcal B,P^\mathcal B)$;

(3) $HomC({}^\mathcal AN,{}^\mathcal AM^\mathcal B)$ is a right $\mathcal B$ pseudomodule determined by $n*(\phi*b)=(n*\phi)*b$,  for any $n\in N, b\in \mathcal B$ and $\phi\in HomC({}^\mathcal AN,{}^\mathcal AM^\mathcal B)$; 

(4) $Chom(P^\mathcal B,{}^\mathcal AM^\mathcal B)$ is a left $\mathcal A$ pseudomodule determined by $(a*\phi)*p:=a*(\phi*p)$ for any $a\in \mathcal A,p\in P$ and $\phi\in Chom(P^\mathcal B,{}^\mathcal AM^\mathcal B)$.
\end{proposition}

\begin{proof} The proofs of (1)--(4) are similar. We only prove (2) as an example. For any   $h,h'\in H$, $\phi\in Chom(^\mathcal AM^\mathcal B,P^\mathcal B)$ and $m\in M$, we have 
$((h\phi)*(h'a))*m=(h\phi)*((h'a)*m)=\sum\limits_{i,j}hSh_i\otimes h'Sh_j\otimes 1\otimes_H\phi_{x_i}(a_{x_j}m)=(h\otimes h'\otimes 1 _H\otimes 1)((\phi*a)*(m))$ and 
$(\phi*a)*(hm)=(1\otimes 1\otimes h\otimes _H1)(\phi*a)*m$. Thus $\phi*a\in H^{\otimes 2}\otimes _HChom(M,P)$ and $*:Chom(M,P)\boxtimes_H\mathcal A\to H^{\otimes 2}\otimes_HChom(M,P)$ is a left $H^{\otimes 2}$ linear map. Moreover, $(\phi*a)*(m*b)=\phi*(a*m*b)=\phi(a*m)*b=((\phi*a)*m)*b$. Hence $\phi*a\in H^{\otimes 2}\otimes_HChom(M^\mathcal B,P^\mathcal B)$.
\end{proof}

From Proposition \ref{pro41} and Proposition \ref{propo48h}, we know that $Chom(^\mathcal AM^\mathcal B,-)$ is a functor from $\mathcal{M}od^\mathcal B$ to $\mathcal{M}od^\mathcal A$ for any $\mathcal A$-$\mathcal B$  bi-pseudomodule $M$.
In particular, $Chom(^\mathcal A\mathcal A^\mathcal A,-)$ is an endofunctor of $\mathcal{M}od^\mathcal A$.  If $M$ is a unital right  pseudomodule over a unital associative $H$-pseudoalgebra $\mathcal A$, and 
 $\Phi:Chom(\mathcal A^\mathcal A,M^\mathcal A)\to H\otimes M$, $\alpha\mapsto \sum\limits_iSh_i\otimes \alpha_{x_i}1$ is a left $H$-module isomorphism given by the proof of Lemma \ref{lem41}, then
$$\begin{array}{lll}(1^{\otimes 2}\otimes_H\Phi)(\alpha*a)&=&\sum\limits_iSh_i\otimes 1\otimes_H\Phi(\alpha_{x_i}a)\\
&=&\sum\limits_{i,j}Sh_i\otimes 1\otimes_H(Sh_j\otimes (\alpha_{x_i} a)_{x_j}1) \end{array}$$ 
and $\Phi(\alpha)*a= (\sum\limits_i Sh_i\otimes \alpha_{x_i}1)*a$. 
For any associative $H$-pseudomodule $\mathcal A$  with identity $1$ such that $a_11=a$ for any $a\in\mathcal A$,  $Cend(\mathcal {A}^\mathcal{A})\simeq H\otimes \mathcal A$ is also an associative $H$-pseudoalgebra with identity $1\otimes 1$ by  Theorem \ref{cor43}.  

\begin{definition}Let  $M$ be an $\mathcal A'$-$\mathcal A$ bi-pseudomodule and $N$ be an  $\mathcal A$-$\mathcal B$ bi-pseudomodule. Assume that  $\mathcal{T}$ is an $\mathcal A'$-$\mathcal B$ sub-bipseudomodule of $M\otimes N$ generated by the following subset $$\left\{\begin{array}{ll}(m_xa)\otimes n-m\otimes (a_{x}n),&a'_{x}(m_{y}a)\otimes n-a'_xm\otimes a_yn\\
m_xa\otimes n_yb-m\otimes (a_xn)_yb,& a'_xm\otimes (a_yn)_zb-a'_x(m_ya)\otimes n_zb\end{array}\mid \begin{array}{l}m\in M,n\in N,x,y,z\in X,\\  a\in \mathcal A,a'\in \mathcal A',b\in \mathcal B\end{array}\right\}.$$  The  quotient bi-pseudomodule $(M\otimes N)/\mathcal T$, denoted by $M\otimes^\mathcal AN$, is called the tensor product of $M$ and $N$ over $\mathcal A$.  We  use 
$m\otimes^\mathcal A n$ to denote the image $m\otimes n+\mathcal{T}$ of $m\otimes n$ in $(M\otimes N)/\mathcal{T}$.  
 \end{definition}

For any associative $H$-pseudoalgebra $\mathcal A$, one can define $M\otimes^\mathcal AN$ for any left $\mathcal A$ pseudomodule $N$ and any right $\mathcal A$ pseudomodule $M$ over $\mathcal A$ by $M\otimes N/\mathcal{T}'$, where $\mathcal{T}'$ is a left $H$-module generated by  $\{(m_xa)\otimes n-m\otimes (a_{x}n)
\mid m\in M,n\in N,x\in X,a\in \mathcal A\}.$ 

\begin{proposition}\label{pro414}Let  $M$ be an $\mathcal A'$-$\mathcal A$ bi-pseudomodule and $N$ an $\mathcal A$-$\mathcal B$ bi-pseudomodule.  Then $M\otimes ^\mathcal AN$ is an $\mathcal A'$-$\mathcal B$
bipseudomodule satisfying  $m*a\otimes^\mathcal  A n=\sum\limits_iS h_{i}\otimes1\otimes _Hm_{x_i}a\otimes n=\sum\limits_i Sh_{i}\otimes1\otimes_H m \otimes a_{x_i}n=m\otimes^\mathcal A a*n$ for any $m\in M$, $a\in \mathcal A$ and $n\in N$.
\end{proposition}

\begin{proof} Let $a',b' \in \mathcal A'$, $a\in\mathcal  A$, $b\in \mathcal B$ and $m\otimes n\in M\otimes N$. Then 
$$\left\{\begin{array}{l} a'*(m\otimes a_yn)=\sum\limits_{i}Sh_i\otimes 1\otimes_H (a'_{x_i}m)\otimes a_yn,\\
a'*((m_ya)\otimes n)=\sum\limits_iSh_i\otimes1\otimes_H a'_{x_{i}}(m_{y}a)\otimes n,\end{array}\right.$$
 $$\left\{\begin{array}{l}a'*(m_xa\otimes n_yb)=\sum\limits_iSh_i\otimes 1\otimes_H a'_{x_i}(m_xa)\otimes n_yb,\\
a'*(m\otimes (a_xn)_yb)=\sum\limits_i Sh_i\otimes 1\otimes_Ha'_{x_i}m\otimes (a_xn)_yb,\end{array}\right.$$
$$\left\{\begin{array}{l}b'*(a'_x(m_ya)\otimes n)=\sum\limits_i Sh_i\otimes 1\otimes_Hb'_{x_i}(a'_x(m_ya))\otimes n=\sum\limits_iSh_i\otimes 1\otimes_H (b'_{x_{i(2)}}a')_{xx_{i(1)}}(m_ya)\otimes n,\\
b'*(a'_xm\otimes a_yn)=\sum\limits_iSh_i\otimes 1\otimes_H b'_{x_i}(a'_xm)\otimes a_yn=\sum\limits_iSh_i\otimes 1\otimes_H(b'_{x_{i(2)}}a')_{xx_{i(2)}}m\otimes a_yn,\end{array}\right.$$
and
$$\left\{\begin{array}{l}\begin{array}{lll} b'*(a'_x(m_ya)\otimes n_zb)&=&\sum\limits_j Sh_j\otimes 1\otimes_Hb'_{x_j}(a'_{x}(m_ya))\otimes n_zb\\
&=&\sum\limits_iSh_i\otimes 1\otimes_H(b'_{x_{i(2)}}a')_{xx_{i(1)}}(m_ya)\otimes n_zb,\end{array}\\
\begin{array}{lll}b'*(a'_xm\otimes (a_yn)_zb)&=&\sum\limits_iSh_i\otimes 1\otimes_H b'_{x_i}(a'_xm)\otimes (a_yn)_zb\\
&=&\sum\limits_i Sh_i\otimes 1\otimes_H(a'_{x_{i(1)}}b')_{xx_{i(1)}}m\otimes (a_yn)_zb.\end{array}\end{array}\right.$$
Define $a'*(m\otimes^\mathcal A n):=\sum\limits_i Sh_i\otimes 1\otimes_H(a'_{x_i}m\otimes^\mathcal A n)$ for $m\otimes^\mathcal A n\in M\otimes^\mathcal AN$, Then $a'*(m\otimes^\mathcal An)$ is well-defined. The rest claims can be routinely  verified.
\end{proof}

\begin{lemma} Let $\mathcal A$ be a unital associative $H$-pseudoalgebra  with identity $1$ such that $a=a_11$ for any $a\in \mathcal A$ and $M$ be a $\mathcal B$-$\mathcal A$ bi-pseudomodule.  If $M$ is a unital right $\mathcal A$ pseudomodule, then $M$ is isomorphic to $M\otimes ^\mathcal A\mathcal A$ as  $\mathcal B$-$\mathcal A$ bi-pseudomodules.\end{lemma}

\begin{proof}
Since $m\otimes^\mathcal A a=m\otimes^\mathcal A a_11=m_1a\otimes^\mathcal A 1$, $M\otimes ^\mathcal A\mathcal A$ is spanned by $m\otimes^\mathcal A 1$. Then the map $\alpha: M\to M\otimes ^\mathcal A\mathcal A$, $m\mapsto m\otimes^\mathcal A 1$ is surjective. Let $\beta:M\otimes^\mathcal A\mathcal A\to M,m\otimes^\mathcal A a\mapsto m_1a$. Then $\alpha \beta(m\otimes^\mathcal A a)=m_1a\otimes ^\mathcal A1=m\otimes ^\mathcal Aa_11=m\otimes^\mathcal A a$. Thus $\alpha\beta=id_{M\otimes^\mathcal AN}$.  Since $M$ is a unital $\mathcal A$-pseudomodule, $\beta\alpha(m)=m_11=m$ and $\beta\alpha=id_{M}$.

Since $m*(ha)= \sum\limits_i Sh_i (h_{(-1)})\otimes 1\otimes_H h_{(2)}(m_{x_i}a)$, $m_1(ha)=h(m_1a)$.  Similarly, as $(h_{(1)}m)*(h_{(2)}a)=\sum\limits_ih_{(1)}Sh_i\otimes 1\otimes_Hm_{x_i}(h_{(2)}a)$,  we get that  $(h_{(1)}m)_1(h_{(2)}a)=m_1(ha)$. Thus $\beta(h(m\otimes a))=m_1(ha)=h(m_1a)=h\beta(m\otimes a).$

For $b\in \mathcal B$, $(1^{\otimes 2}\otimes_H\beta)(b*(m\otimes^\mathcal A a))=(1^{\otimes 2}\otimes_H\beta)(b*(m_1a\otimes^\mathcal A 1))=\sum\limits_i Sh_i\otimes 1\otimes_H(b_{x_i}(m_1a))_11=\sum\limits_i Sh_i\otimes 1\otimes _H b_{x_i}(m_1a)=b*\beta(m\otimes^\mathcal A a)$. Similarly, $(1^{\otimes 2}\otimes_H\beta)((m\otimes^\mathcal A a)*a')=\sum\limits_j1\otimes g_j\otimes_H\beta(m\otimes^\mathcal A b_j)$ and $(m_1a)*a'=\sum\limits_i Sh_i\otimes 1\otimes_H(m_1a)_{x_i}a'=\sum\limits_i Sh_i\otimes 1\otimes_Hm_1(a_{x_i}a')=\sum\limits_j1\otimes g_j\otimes_H m_1b_j$, where $a*a'=\sum\limits_j1\otimes g_j\otimes b_j$. Hence $\beta$ is a homomorphism of $\mathcal B$-$\mathcal A$ bi-pseudomodules.
\end{proof}

\begin{lemma}Let $N$ be an $\mathcal A$-$\mathcal B$ bi-pseudomodule, $M$ be a right $\mathcal A$-pseudomodule and $L$ be a left $\mathcal B$-pseudomodule. Then
$(M\otimes ^\mathcal AN)\otimes^\mathcal BL$ is isomorphic to $M\otimes^\mathcal A(N\otimes ^\mathcal BL)$.
\end{lemma}
\begin{proof}For $m\in M$, $n\in N$, $l\in L$, $a\in\mathcal A$, $b\in \mathcal B$ and $x,y\in X$, we have $(m_xa\otimes^\mathcal A n)_yb\otimes^\mathcal B l=(m_xa\otimes^{\mathcal A} n_yb)\otimes^\mathcal B l=m_xa\otimes^\mathcal A (n_yb\otimes^\mathcal B l)=m_xa\otimes^{\mathcal A }(n\otimes^\mathcal B b_yl)=m\otimes^\mathcal A a_x(n\otimes^\mathcal B b_yl)=m\otimes^{\mathcal A }((a_xn)\otimes^\mathcal B b_yl)$.
Thus the map $(m\otimes ^\mathcal An)\otimes^\mathcal B l\mapsto m\otimes^\mathcal A (n\otimes^\mathcal B l)$ is an isomorphism from $(M\otimes ^\mathcal AN)\otimes ^\mathcal BL$ to $M\otimes^\mathcal A(N\otimes^\mathcal BL)$.
\end{proof}

For any  finitely generated  right $\mathcal A$-pseudomodule $M$,  let $M^{c*}:=Chom(M^\mathcal A,\mathcal A^\mathcal A)$, and $\mathcal A'=Cend(M^\mathcal A)$. Note that $Chom(M^\mathcal A,\mathcal A^\mathcal A)\subseteq Chom(M,\mathcal A)\simeq Hom_H(M,H\otimes \mathcal A)$. We can regarded $M$ as a left $\mathcal A'$ pseudomodule such that $\alpha*(m*a)=(\alpha*m)*a$ for any $\alpha\in \mathcal A'$, $a\in\mathcal A$ and $m\in M$ by Proposition \ref{propo48h}. If   $M$  is a  finitely generated $\mathcal A$  pseudomodule, then we can define $\mathcal A$ and $\mathcal A'$ actions on $M^{c*}$ as follows. For any $a\in \mathcal A$, $a'\in \mathcal A'$, $m\in M$ and $\phi\in M^{c*}$,
$$(a*\phi)*m:=a*(\phi*m)\in H^{\otimes 3}\otimes_H \mathcal A,\qquad (\phi*a')*m=\phi*(a'*m)\in H^{\otimes 3}\otimes_H\mathcal A.$$
 Since $((a*\phi)*a')*m=a*(\phi*(a'*m))=(a*(\phi*a'))*m$ for any $m\in M$, $(a*\phi)*a'=a*(\phi*a')$. This means that $M^{c*}$ is an $\mathcal A$-$\mathcal A'$ bi-pseudomodule.  Moreover,  $M^{c*}\otimes^{\mathcal A'}M$ is an $\mathcal A$-$\mathcal A$ bi-pseudomodule by Proposition \ref{pro414}, where 
\begin{eqnarray}h(m'\otimes^{\mathcal  A'} m)=h_{(1)}m'\otimes^{\mathcal  A'} h_{(2)}m,
\end{eqnarray}
\begin{eqnarray}a*(m'\otimes ^{\mathcal  A'}m)=\sum\limits_iSh_i\otimes1\otimes_H((a_{x_i}m')\otimes^{\mathcal  A'} m)
\end{eqnarray}
and \begin{eqnarray}(m'\otimes^{\mathcal  A'} m)*a=\sum\limits_j1\otimes g_j\otimes_H(m'\otimes^{\mathcal  A'} n_j)
\end{eqnarray}
provided that $m*a=\sum\limits_j1\otimes g_j\otimes_Hn_j$.

Similarly, $M\otimes^\mathcal AM^{c*}$ is an $\mathcal A'$-$\mathcal A'$ bi-pseudomodule, where $(\phi*a')*m:=\phi*(a'*m)$ for $\phi\in M^{c*}$, $a'\in Cend(^\mathcal AM)$ and $m\in M$.
Since $(a*\phi)*m=a*(\phi*m)$ for $a\in \mathcal A, $ $\phi\in M^{c*}$ and $m\in M$, $((\phi*m)*\psi)*n=(\phi*m)*(\psi*n)=\phi*(m*(\psi*n))$ for $m\in M$ and $\phi,\psi\in M^{c*}$.
For any $\phi\in M^{c*}$ and $m,n\in M$, define $\mu:M\otimes^\mathcal AM^{c*}\to H^{\otimes 2} \otimes_H\mathcal A'$ by $$\mu(m\otimes^{\mathcal  A} \phi)*n:=m*(\phi*n)$$ and $\tau:M^{c*}\otimes^{\mathcal A'}M\to H^{\otimes 2}\otimes_H\mathcal A$ by
$$\tau(\phi\otimes{}^{\mathcal  A'} m):=\phi*m.$$
Then $n*\tau(\phi\otimes^{\mathcal  A'} m)=n*(\phi*m)=\mu(n\otimes^{\mathcal  A}  \phi)*m$. Since $(\phi*\mu (m\otimes^{\mathcal  A} \psi))*n=\phi*(\mu(m\otimes^{\mathcal  A} \psi)*n)=\phi*(m*(\psi*n))$,  $(\phi*\mu (m\otimes^{\mathcal  A} \psi))*n=((\phi*m)*\psi)*n$  and $\phi*\mu (m\otimes^{\mathcal  A} \psi)=(\phi*m)*\psi.$ Hence 
$\tau(\phi\otimes^{\mathcal  A'} m)*\psi=(\phi*m)*\psi=\phi*\mu(m\otimes^{\mathcal  A} \psi)$.

From the above arguments, we can give a more general definition as follows. 

\begin{definition}A Morita context is a set $(\mathcal A,\mathcal A',M,M',\tau,\mu)$, where $\mathcal A$ and $\mathcal A'$ are $H$-pseudoalgebras, $M$ is an left $\mathcal A'$ and right $\mathcal A$ bi-pseudomodule,
$M'$ is a left $\mathcal A$ and right $\mathcal A'$ bi-pseudomodule, $\tau: M'\otimes ^{\mathcal A'}M\to H^{\otimes 2}\otimes_H\mathcal A$ is an  $\mathcal A$-bilinear $H^{\otimes 2}$-linear map and $\mu:M\otimes^\mathcal AM'\to H^{\otimes 2}\otimes_H\mathcal A'$ is an $\mathcal A'$-bilinear and $H^{\otimes 2}$-linear map such that $\mu(n\otimes^\mathcal A m')*m=n*\tau(m'\otimes^{\mathcal  A'} m)$ and $\tau(m'\otimes^{\mathcal  A'} m)*n'=m'*\mu(m\otimes^\mathcal A n')$.

$\mu$ (resp. $\tau$) is said to be surjective if there are $m_i\in M$ and $m_i'\in M'$ such that $\sum\limits_i \mu(m_i\otimes^\mathcal A m_i')=1\otimes 1\otimes_Ha'$ (resp.  $\sum\limits_i\tau(m_i'\otimes^{\mathcal  A'} m_i)=1\otimes 1\otimes_Ha$) for any $a'\in \mathcal A'$ (resp.  $a\in \mathcal A$).
\end{definition}

\begin{example}\label{ex419} Let $\mathcal A$ be an associative $H$-pseudoalgebra with an identity $1$ such that $a_11=a$ for any $a\in\mathcal A$.  Assume that $M:= \mathcal A^n$ is a direct sum of $\mathcal A$ as a right pseudomodule over itself.  We write the element $m$ in $M$ as a column vector, that is, $m=\left(\begin{array}{c}a_1\\ a_2\\ \vdots\\ a_n\end{array}\right)$, denoted by $(a_1,a_2,\cdots,a_n)^T$, where $a_i\in \mathcal A$.  
Then $M^{c*}=Chom(M^\mathcal A,\mathcal A^\mathcal A)\simeq (H\otimes \mathcal A)^n$ by Proposition \ref{pro41} and Lemma \ref{lem41} and $\mathcal A'=Cend(M^\mathcal A)\simeq  M_n(H\otimes\mathcal A)$ by Remark \ref{rem44}. We write the element in $M^{c*}$ as a row vector. Let $(h_{ij}\otimes a_{ij})\in M_n(H\otimes \mathcal A)$ and $(a_1,a_2,\cdots, a_n)^T\in M$, where $h_{ij}\otimes a_{ij}\in H\otimes \mathcal A$. Then $(h_{ij}\otimes a_{ij})(a_1,a_2,\cdots,a_n)^T=(\sum\limits_{j=1}^n(h_{1j}\otimes a_{1j})*a_j, \sum\limits_{j=1}^n(h_{2j}\otimes a_{2j})*a_j,\cdots,\sum\limits_{j=1}^n(h_{nj}\otimes a_{ij})*a_j)^T$, where $(h_{ij}\otimes a_{ij})*a_j=\sum\limits_kh_{ij}Sh_k\otimes1\otimes_H1(a_{ij})_{x_k}a_j$ by Theorem \ref{cor43}. Thus $M$ is an $\mathcal A'$-$\mathcal A$ bi-pseudomodule. In addition, $M^{c*}$ is a left $\mathcal A$ pseudomodule with the action $a*(h_1\otimes a_1,h_2\otimes a_2,\cdots, h_n\otimes a_n)=(a*(h_1\otimes a_1),a*(h_2\otimes a_2),\cdots, a*(h_n\otimes a_n))$, where 
$a*(h_i\otimes a_i)=\sum\limits_jSh_j\otimes h_i\otimes_H(a)_{x_j}a_i$ by Lemma \ref{le41}. Consequently, $M^{c*}$ is an $\mathcal A$-$\mathcal A'$ bi-pseudomodule.

Moreover, $\tau((h_1\otimes b_1,h_2\otimes b_2,\cdots, h_n\otimes b_n)\otimes^\mathcal A(a_1,a_2,\cdots, a_n)^T)=\sum\limits_j(h_j\otimes b_j)*a_j$ and 
$\mu((a_1,a_2,\cdots,$ $ a_n)^T\otimes^\mathcal A(h_1\otimes b_1,h_2\otimes b_2,\cdots, h_n\otimes b_n))=( a_i*(h_j\otimes b_j))$, where $a_i*(h_j\otimes b_j)=\sum\limits_kSh_k\otimes h_j\otimes (b_j)_{x_k}a_i$. Therefore, $(\mathcal A, M_n(H\otimes \mathcal A), \mathcal A^n, (H\otimes \mathcal A)^n,\tau,\mu)$ is a Mortia context. 
\end{example}

\begin{remark}Let $(\mathcal A,\mathcal A',M,M',\tau,\mu)$ be a Morita context. Then $\left(\begin{array}{cc}\mathcal A'& M\\ M'&\mathcal A\end{array}\right)$ is an associative $H$-pseudoalgebra, where
$$\left(\begin{array}{cc}a'& m\\ m'&a\end{array}\right)*\left(\begin{array}{cc}b'& n\\ n'&b\end{array}\right)=\left(\begin{array}{cc}a'*b'+\mu(m\otimes^{\mathcal  A} n')& a'*n+m*b\\ m'*b'+a*n'&\tau(m'\otimes^{\mathcal  A'} m)+a*b\end{array}\right).$$\end{remark}

Assume that $\mathcal A$ is an associative $H$-pseudoalgebra. Let $\mathcal{PM}_\mathcal A$ be the category of all right $\mathcal A$ pseuomodules, whose morphisms are homomorphisms of right $\mathcal A$ pseudomodules. Recall that a left $H$-module homomorphism $\phi$ from a right $\mathcal A$ pseudomodule $M$ to a right $\mathcal A$ pseudomodule $N$ is said to be homomorphism of right $\mathcal A$ pseudomodules if $\phi(m)*a=(1^{\otimes 2}\otimes_H\phi)(m*a)$ for any $m\in M$ and $a\in \mathcal A$. If $\mathcal A$ is a unital associative $H$-pseudoalgebra, then we use $\mathcal{PM}_\mathcal A^{\circ}$ to denote  the subcategory of $\mathcal{PM}_\mathcal A$ consisting of all unital right $\mathcal A$ pseudomodules.

\begin{example}Let $M$ be an $\mathcal A$-$\mathcal B$ bi-pseudomodule. Then $F(N):=N\otimes^\mathcal AM$ and $F(\phi):N\otimes ^\mathcal AM\to N'\otimes ^\mathcal AM, n\otimes^\mathcal A m\mapsto \phi(n)\otimes^\mathcal A m$ is a functor from $\mathcal{PM}_\mathcal A$ to $\mathcal{PM}_\mathcal B$, where $\phi:N\to N'$ is a homomorphism of right $\mathcal A$ pseudomodules.

For any right $\mathcal A$ pseudomodule $M$,  $H\otimes M$ becomes a right $\mathcal A$ pseudomodule with the following actions
$$(h\otimes m)*a=\sum\limits_iSh_i\otimes 1\otimes_H(h\otimes m_{x_i}a)$$
and $$h'(h\otimes m)=h'h\otimes m.$$
The mappings $M\to H\otimes M$ and $\phi \to id_H\otimes \phi$ is also a functor from $\mathcal{PM}_\mathcal A$ to itself. This functor is denoted by $H\otimes_-$. It is obvious that $H\otimes M$ is also a unital $\mathcal A$ pseudomodule provided that $M$ is a unital pseudomodule. Thus $H\otimes_-$ is an endofunctor of the category $\mathcal{PM}_\mathcal A^{\circ}$. 
\end{example}

By Theorem \ref{cor43} and Proposition \ref{propo48h}, $Chom({}^\mathcal A\mathcal A^\mathcal A,-)$ is also an endofunctor of $\mathcal{PM}_\mathcal A$.

\begin{proposition} For each $M\in \mathcal{PM}^\circ_\mathcal A$, the map $\Phi:Chom(^\mathcal A\mathcal A^\mathcal A, M^\mathcal A)\to H\otimes M$, $\varphi\mapsto \sum\limits_i Sh_i\otimes \varphi_{x_i}1$ defines  a natural equivalence from the  functor  $Chom(^\mathcal A\mathcal A^\mathcal A,-)$ to the functor  $H\otimes_-$.
\end{proposition}

\begin{proof}We only need to verify $\Phi^{-1}$ is an isomorphism of right $\mathcal A$ pseudomodules. For any  $\sum\limits_i Sh_i\otimes m_i\in H\otimes M$ and $a,b\in\mathcal A$,  $(\Phi^{-1}(\sum\limits_i Sh_i\otimes m_i)*a)*b=\Phi^{-1}(\sum\limits_iSh_i\otimes m_i)(a*b)=\sum\limits_j(1\otimes (Sh_j\otimes 1)\Delta\otimes_H1)\Phi^{-1}(\sum\limits_iSh_i\otimes m_i)(a_{x_j}b)=
\sum\limits_{i,j,k}Sh_k\otimes Sh_j\otimes1\otimes_HSh_i\otimes (m_i)_{x_k}(a_{x_j}b)$. Since $m_i*(a*b)=(m_i*a)*b$, we have 
$\sum\limits_{i,j,k}Sh_k\otimes Sh_j\otimes1\otimes_HSh_i\otimes (m_i)_{x_k}(a_{x_j} b)=
\sum\limits_{i,j,k} Sh_kSh_{j(1)}\otimes Sh_{j(2)}\otimes 1\otimes_HSh_i\otimes (m_{ix_k}a)_{x_j}b=(1\otimes 1\otimes_H\Phi^{-1})((\sum\limits_iSh_i\otimes m_i)*a)*b$.
Thus $\Phi^{-1}$ is a homomorphism of right $\mathcal A$ pseudomodules.
\end{proof}

\begin{definition}Two associative unital $H$-pseudoalgebras $\mathcal A$ and $\mathcal B$ are said to be conformal Morita equivalent if there is a functor $F$ from $\mathcal{PM}_\mathcal A^{\circ}$  to $\mathcal{PM}_\mathcal B^{\circ}$  and a functor $G$ from $\mathcal{PM}_\mathcal B^{\circ}$ to $\mathcal{PM}_\mathcal A^{\circ}$ such that $F\circ G\simeq H\otimes_-$ and $G\circ F\simeq H\otimes_-$.
\end{definition}

\begin{theorem} ( {\bf Morita Theorem}) Let $(\mathcal A,\mathcal A',M,M',\tau,\mu)$ be a Morita context such that $\mu$ and $\tau$ are surjective, where both $\mathcal A$ and $\mathcal A'$ are unital with   identities  $1$ and $1'$  respectively. Then 

(1) There are finite elements $n_j\in M$ (resp. $m_j\in M$)  and finite element $\phi_i\in Chom(M^\mathcal A,$ $\mathcal A^\mathcal A)$ (resp. $\psi_j\in HomC(^{\mathcal A'}M,{}^{\mathcal A'}{\mathcal A'})$)  such that $\sum\limits_i Sh_{i(1)}\otimes Sh_{i(2)}\otimes 1\otimes_H1'_{x_i}m=\sum\limits_j n_j*(\phi_j*m)$ (resp. $\sum\limits_i Sh_i\otimes 1\otimes 1\otimes_H m_{x_i}1=\sum\limits_j(m*\psi_j)*m_j$) for any $m\in M$. If $M$ is a unital $\mathcal{A}'$ pseudomodule, then $M$ is generated by $n_j$ as a right $\mathcal A$ pseudomodule.
Similar result for left $\mathcal A$ psedomodule $M'$ (resp. right $\mathcal A'$ pseudomodule $M'$) holds.

(2) If $M$ (resp. $M'$)  is a right unital $\mathcal A$ pseudomodule (resp. left unital $\mathcal A'$ pseudomodule), then $\tau$ (resp. $\mu$) is injective.

(3) If $M$ is a left unital $\mathcal A'$ pseudomodule, then 
the map $l(m'):m\mapsto \tau(m'\otimes^{\mathcal A'} m)$ is an isomorphism from $M'{}^{\mathcal A'}$ to $Chom(M^\mathcal A,\mathcal A^\mathcal A)$. Similarly, if $M'$ is a left unital $\mathcal A$ pseudomodule, then the map $l'(m): m'\mapsto \mu(m\otimes^\mathcal A m')$ is an isomorphism from $M$ to $Chom(^{\mathcal A'}M,{}^{\mathcal A'}\mathcal A')$.

(4) The map $\lambda:\mathcal A'\to Cend(M^\mathcal A); a'\mapsto \lambda(a')$, where $\lambda(a')*m:=a'*m$, is an associative $H$-pseudoalgebra isomorphism from $\mathcal A'$ to $Cend(M^{\mathcal A})$.
Similarly, the map $\lambda':\mathcal A\to Cend(M'{}^{\mathcal A'}); a\mapsto \lambda'(a)$, where $\lambda'(a)*m':=a*m'$, is an associative $H$-pseudoalgebra isomorphism from $\mathcal A$ to $Cend(M'{}^{\mathcal A'})$. 

(4$^\prime$) The map $\Lambda:\mathcal A\to EndC(^{\mathcal A'}M),a\mapsto \Lambda(a)$,  where $m*(\Lambda(a))=m*a$,  is an isomorphism of associative $H$-pseudoalgebras
and the map $\Lambda':\mathcal A'\to EndC(M^\mathcal A), a'\mapsto \Lambda'(a')$, where $m'*(\Lambda'(a'))=m'*a'$, is an isomorphism of associative $H$-pseudoalgebras.

(5) $\mathcal A$ and $\mathcal A'$ are conformal Morita equivalent.

(6) The center $Z(\mathcal A)$ of $\mathcal A$ is isomorphic to the center  $Z(\mathcal A')$ of $\mathcal A'$.
\end{theorem}

\begin{proof}(1) If $\mu$ is surjective, then there are $n_j\in M $ and $n_j'\in M'$ such that $\sum\limits_j\mu(n_j\otimes^\mathcal A n_j')=1\otimes 1\otimes_H1'$. Thus $\sum\limits_iSh_{i(1)}\otimes Sh_{i(2)}\otimes1\otimes_H1'_{x_i}m=(\Delta\otimes 1\otimes_H1) 1'*m=(1\otimes 1\otimes_H1')*m=\sum\limits_j\mu (n_j\otimes ^\mathcal An_j')*m=\sum\limits_jn_j*\tau(n_j'\otimes^{\mathcal A'} m).$ Hence $m=(\varepsilon^{\otimes 3}\otimes_H1)\sum\limits_jn_j*\tau(n_j'\otimes^{\mathcal  A'}m)$ and $M$ is generated by $n_j$. The maps $l(n_j'):M\to \mathcal A$, where $l(n_j')(m)=\tau(n_j'\otimes^{\mathcal A'} m)$, are in $Chom(M^\mathcal A,\mathcal A^\mathcal A)$. 

Similarly, if there are $m_j\in M$ and $m_j'\in M'$ such that $\sum\limits_j\tau(m_j'\otimes^{\mathcal A'} m_j)=1\otimes 1\otimes_H1$, then $\sum\limits_i Sh_i\otimes 1\otimes 1\otimes_Hm_{x_i}1=m*(1\otimes 1\otimes_H1)=\sum\limits_jm*\tau(m_j'\otimes^{\mathcal A'} m_j)=\sum\limits_j\mu(m\otimes^\mathcal A m_j')*m_j$.
The maps $r(m_j')$, where $m*(r(m_j')):=\mu(m\otimes^\mathcal A m_j')$, are in $HomC(^{\mathcal A'}M,{}^{\mathcal A'}\mathcal A')$.

(2) Assume that  $\sum\limits_i\tau (m_i'\otimes^{\mathcal A'} m_i)=0$ and $\sum\limits_j\tau(n_j'\otimes^{\mathcal A'} n_j)=1\otimes 1\otimes_H1\in H^{\otimes 2}\otimes_H\mathcal A$. Then $\sum\limits_{ij}(m_i'\otimes ^{\mathcal A'}m_i)*\tau(n_j'\otimes^{\mathcal A'} n_j)=\sum\limits_i(1\otimes \Delta\otimes_H1)(m_i'\otimes^{\mathcal A'} m_i)*1=(1\otimes \Delta\otimes_H1)\sum\limits_{ij}m_i'\otimes^{\mathcal A'} (m_i*\tau(n_j'\otimes^{\mathcal A'} n_j)) =(1\otimes \Delta\otimes_H1)\sum\limits_{ij}m_i'\otimes^{\mathcal A'} \mu(m_i\otimes^\mathcal A n_j')*n_j= (1\otimes \Delta\otimes_H1)\sum\limits_{ij}m_i'*\mu(m_i\otimes^{\mathcal A} n_j')\otimes^{\mathcal A'} n_j=(1\otimes \Delta\otimes_H1)\sum\limits_{ij}\tau(m_i'\otimes^{\mathcal A'} m_i)*n_j'\otimes^{\mathcal A'} n_j=0.$ Hence $\sum\limits_i(m_i'\otimes^{\mathcal A'} m_i)*1=0$. Let $m_i*1=\sum\limits_j 1\otimes g_{ij}\otimes _Hn_{ij}=\sum\limits_j g_{ij(-1)}\otimes 1\otimes_Hg_{ij(2)}n_{ij}$. Then $m_i{}_11=\sum\limits_jg_{ij}n_{ij}=m_i$ and $\sum\limits_{i,j}1\otimes g_{ij}\otimes_Hm_i' \otimes^{\mathcal A'}n_{ij}=0$. Thus
 $\sum\limits_im_i'\otimes ^{\mathcal A'}m_i=0$ and $\tau$ is injective.

(3) Define $l(m')*m:=\tau(m'\otimes^{\mathcal A'} m)$ for any $m'\in M'$. Then  $l(m')\in Chom (M,\mathcal A)$. Since $l(m')*(m*a)=\sum\limits_i(1\otimes (Sh_i\otimes 1)\Delta\otimes_H1)\tau(m'\otimes^{\mathcal A'}m_{x_i}a)=\tau(m'\otimes^{\mathcal A'}m*a)=
(l(m')*m)*a$, $l(m')\in Chom(M^\mathcal A,\mathcal A^\mathcal A)$. For any $a'\in \mathcal A$, $a\in \mathcal A$ and $m'\in M'$,
we have $l(a*m'*a')*m=\tau((a*m'*a')\otimes^{\mathcal A'} m)=a*\tau(m'\otimes^{\mathcal A'} a'* m)=a*(l(m')*(a'*m))=(a*l(m')*a')*(m).$ Thus $l(a*m'*a')=a*l(m')*a'$ and $l:M'\to Chom(M^{\mathcal A},\mathcal A^\mathcal A)$ is an $\mathcal A$-$\mathcal A'$ bi-pseudomodule homomorphism. If $l(m')=0$, then $\tau(m'\otimes^{\mathcal A'} m)=0$ for any  $m\in M$. Since  there are $n_j\in M$ and $n_j'\in M'$ such that $\sum\limits_j\mu(n_j\otimes^\mathcal A n_j')=1\otimes1\otimes_H1$,  $m'*(1\otimes 1\otimes_H1)=m'*\sum\limits_j\mu(n_j\otimes ^{\mathcal A}n_j')=\sum\limits_j\tau(m'\otimes^{\mathcal A'} n_j)*n_j'=0$.
Hence $m'=0$ and $l$ is injective. Now let $\phi\in Chom(M^\mathcal A,\mathcal A^\mathcal A)$. As before, we write $1\otimes 1\otimes_H1=\sum\limits_i\mu(n_i\otimes^\mathcal A n_i')$. Let
$\sum\limits_kg_k\otimes g_k'\otimes 1\otimes_H m_k'=\sum\limits_i(\phi*n_i)*n_i'$. Then $\sum\limits_k\tau(g_k\otimes g_k'\otimes 1\otimes_Hm_k'\otimes^{\mathcal A'} m)=\sum\limits_{i}(\phi*n_i)*\tau(n_i'\otimes^{\otimes A'} m)=\sum\limits_{i}\phi*(n_i*\tau(n_i'\otimes^{\mathcal A'} m))=\sum\limits_{ki}\phi*(\mu(n_i\otimes^\mathcal A n_i')*m)=\phi*((1\otimes 1\otimes_H1)*m)=\sum\limits_i\phi*(Sh_{i(1)}\otimes Sh_{i(2)}\otimes 1\otimes_H1_{x_i}m)=\sum\limits_{i,j} Sh_j\otimes Sh_{i(1)}\otimes Sh_{i(2)}\otimes1\otimes_H\phi_{x_j}(1_{x_i}m).$ 
Since $\sum\limits_k( 1\otimes (g_k\otimes g_k'\otimes1)\Delta\otimes_H1)\tau(m_k'\otimes^{\mathcal A'} m)=
\sum\limits_{ki}Sh_{i}\otimes g_k\otimes g_k'\otimes 1\otimes_H\tau((m_k')_{x_i}m),$  
$$\sum\limits_{ki}Sh_{i}\otimes g_k\otimes g_k'\otimes 1\otimes_H\tau((m_k')_{x_i}m)=
 \sum\limits_{i,j} Sh_j\otimes Sh_{i(1)}\otimes Sh_{i(2)}\otimes1\otimes_H\phi_{x_j}(1_{x_i}m).$$ 
Thus 
$\sum\limits_k \varepsilon (g_kg_k')\tau(m_k'\otimes^{\mathcal A'}m)= \sum\limits_k  \varepsilon(g_k'g_k)Sh_i \otimes_H\tau((m_k')_{x_i}m)=\sum\limits_jSh_j\otimes 1\otimes_H\phi_{x_j}(m)=\phi(m).$ Hence $\phi(m)=\tau((\sum\limits_k\varepsilon(g_{k}'g_k)m'_{k})\otimes^{\mathcal A'} m)$ and $l$ is surjective.

(4) Since $M$ is an $\mathcal A'$-$\mathcal A$ bi-pseudomodule, $\lambda(a')\in Cend(M^\mathcal A)$. As $\lambda(a'*b')*m=\lambda(a')*(\lambda(b')*m)=(\lambda(a')*\lambda(b'))*m$. This means that 
$\lambda$ is an associative $H$-pseudoalgebra homomorphism. If $\lambda(a')=0$, then $\sum\limits_iSh_i\otimes 1\otimes 1\otimes_Ha'_{x_i}1'=\sum\limits_j a'*\mu(n_j\otimes^\mathcal A n_j')=
\sum\limits_j\mu((a'*n_j)\otimes^\mathcal A n_j')=0$, where $\sum\limits_j\mu(n_j*n_j')=1\otimes 1\otimes _H1'$. Thus $a'=a'_{1}1'=0$ and $\lambda$ is injective. Now let $\phi\in Cend(M^{\mathcal A})$.
Assume that  $\sum\limits_j\mu((\phi*n_j)\otimes^\mathcal A n_j')=\sum\limits_j\sum\limits_{j_i}f^1_{j_i}\otimes f^2_{j_i}\otimes f^3_{j_i}\otimes 1\otimes_Ha_{j_i}'$. Then
$\sum\limits_j\sum\limits_{j_i}((f^1_{j_i}\otimes f^2_{j_i}\otimes f^3_{j_i}\otimes1)\Delta\otimes 1\otimes_H1)a_{j_i}*m=\sum\limits_j\mu((\phi*n_j)\otimes^\mathcal A n_j')*m=\sum\limits_j(\phi*n_j)*\tau(n_j'\otimes^{\mathcal A'} m)=
\sum\limits_j\phi(n_j*\tau(n_j'\otimes^{\mathcal A'} m))=\sum\limits_j\phi((\mu(n_j\otimes^\mathcal A n_j')*m)=\sum\limits_i(1\otimes (Sh_{i(1)}\otimes Sh_{i(2)}\otimes1)\Delta\otimes_H1)\phi(1'_{x_i}m).$
Hence $\phi(m)=\sum\limits_j\sum\limits_{j_i}\varepsilon(f_{j_i}^1f^2_{j_i}f^3_{j_3})a_{j_i}*m$. By now we have proved that $\lambda$ is an isomorphism of associative $H$-pseudoalgebras. Similarly, we can prove that $\lambda'$ is an isomorphism of associative $H$-pseudoalgebras.

(4$^\prime$) This can be proved in the same way as (4).

(5) Let $\tau^{\circ}:M'\otimes^{\mathcal A'}M\to H\otimes A, $ $\sum n_i'\otimes n_i\mapsto \sum\limits_{i,j}Sh_j\otimes \tau({n_i'}_{x_j}n_i)$, where 
$\tau(\sum n_i'\otimes^{\mathcal A'} n_i)=\sum\limits_{i,j}Sh_j\otimes 1\otimes_H\tau({n_i'}_{x_j}n_i)$. Since $\tau$ is surjective, so is
$\tau^{\circ}$. Thus  $\tau^{\circ}$ is an isomorphism.  For any right $\mathcal {A}$ pseudomodule $N$,  define $\gamma:N\otimes^\mathcal A(M'\otimes^{\mathcal A'}M)\to H\otimes N$ via $\gamma (n\otimes ^\mathcal A m'\otimes^{\mathcal A'} m)=\sum\limits_jSh_j\otimes n_1(\tau^{\circ}(m'_{x_j}m))=(1\otimes \beta)((12)(n\otimes \tau^{\circ}(m'\otimes m)))$.  Then $\gamma$ is a bijection.  It is easy to check that 
$_-\otimes^\mathcal AM'\circ _-\otimes^{\mathcal A'}M\simeq H\otimes_-$ and $_-\otimes^{\mathcal A'}M\circ _-\otimes^{\mathcal A}M'\simeq H\otimes_-$ as functors.   Therefore $\mathcal A$ and $\mathcal A'$ are conformal Mortia equivalent.

(6) For any $\phi\in Cend(M^\mathcal A)\cap EndC(^{\mathcal A'}M)$ and any $\psi\in Cend(M^\mathcal A)$, we have  $\psi=\lambda(a')$ for some $a'\in \mathcal A'$ and $\phi=\Lambda(a)$ for some $a\in \mathcal A$. Thus  $\psi *(m*\phi)=a'*(m*a)=(a'*m)*a=(\psi*m)*\phi$ for any $m\in M$. Consequently, $\phi\in Z(Cend(M^\mathcal A))$. Conversely, if $\phi\in Z(Cend(M^\mathcal A))$ and $a'\in\mathcal  A$, then $(a'*m)*\phi=(\lambda(a')*m)*\phi=\lambda(a')(m*\phi)=a'*(m*\phi)$ for any $m\in M$. Thus $\phi\in EndC(^{\mathcal A'}M)$ if we write $\phi(m)$ as $m*\phi$. Hence  $Z(Cend(M^\mathcal A))=Cend(M^{\mathcal A})\cap EndC(^{\mathcal A'}M)$. Similarly, we have $Z(EndC(^{\mathcal A'}M))=Cend(M^\mathcal A)\cap EndC(^{\mathcal A'}M)$.
\end{proof}

\enddocument